\documentclass[preprint,12pt]{elsarticle}
\usepackage{color}
\usepackage{xcolor}
\usepackage[normalem]{ulem}

\usepackage{amssymb}
\usepackage{amsmath}
\usepackage{caption}
\usepackage{subcaption}
\usepackage{lineno,hyperref}
\usepackage{array} 
\usepackage{booktabs}
\usepackage{multirow}
\usepackage{bm}
\usepackage{appendix}
\usepackage{graphicx}
\graphicspath{{Imagesmcmcnet/}{Imagesmcmcnetadded/}}

\usepackage{tikz}
\usepackage{mathrsfs}
\usepackage[left=2cm,right=3cm,top=2cm,bottom=2cm]{geometry}

\newcommand{\Ac}{\mathcal{A}}

\newcommand{\Dc}{\mathcal{D}}
\newcommand{\Ec}{\mathcal{E}}

\newcommand{\Gc}{\mathcal{G}}

\newcommand{\Nb}{\mathbb{N}}

\newcommand{\Rb}{\mathbb{R}}

\newtheorem{theorem}{Theorem}
\usepackage{lipsum} 

\newcommand{\kom}[1]{}
\renewcommand{\kom}[1]{{\bf [#1]}}
\definecolor{deepgreen}{rgb}{0.0, 0.8, 0.0}
\definecolor{deepmagenta}{rgb}{0.5, 0.0, 0.5}
\definecolor{blau}{rgb}{0.1,0.0,0.9}

\newcounter{komcounter}
\numberwithin{komcounter}{section}

\begin{document}
\begin{frontmatter}



\title{MCMC-Net: Accelerating Markov Chain Monte Carlo with Neural Networks for Inverse Problems}
\tnotemark[1,2]


\author[inst1]{Sudeb Majee}
\ead{smajee@uncc.edu}
\affiliation[inst1]{
organization={Department of Mathematics and Statistics},
addressline={UNC Charlotte}, 
city={Charlotte},
postcode={28223}, 
state={NC},
country={USA}
}

\author[inst2]{Anuj Abhishek\corref{cor1}}
\ead{axa1828@case.edu}
\affiliation[inst2]{organization={Department of Mathematics, Applied Mathematics and Statistics},
            addressline={Case Western Reserve University}, 
            city={Cleveland},
            postcode={44106-7058}, 
            state={OH},
            country={USA}}
            
\author[inst3]{Thilo Strauss}
\ead{Thilo.Strauss@xjtlu.edu.cn}
\affiliation[inst3]{organization={Xi'an Jiaotong-Liverpool University},
            addressline={School of AI and Advanced Computing},  
            city={Suzhou},
            postcode={215000}, 
            state={Jiangsu},
            country={China}}

\affiliation[inst4]{
organization={ Center for Trustworthy Artificial Intelligence through Model Risk Management (TAIMing AI)},
addressline={UNC Charlotte}, 
city={Charlotte},
postcode={28223}, 
state={NC},
country={USA}
}

\author[inst1,inst4]{Taufiquar Khan}    
\ead{taufiquar.khan@charlotte.edu}
\cortext[cor1]{Corresponding author}

\begin{abstract}
In many computational problems, using the Markov Chain Monte Carlo (MCMC) can be prohibitively time-consuming. We propose MCMC-Net, a simple yet efficient way to accelerate MCMC via neural networks. The key idea of our approach is to substitute the true likelihood function of the MCMC method with a neural operator based surrogate. We extensively evaluate the accuracy and speedup of our method on three different PDE-based inverse problems where likelihood computations are computationally expensive, namely electrical impedance tomography, diffuse optical tomography, and quantitative photoacoustic tomography. 

MCMC-Net performs similar to the classical likelihood counterpart but with a significant speedup. We conjecture that the method can be applied to any problem with a sufficiently expensive likelihood function. We also analyze MCMC-Net in a theoretical setting for the different use cases. We prove a universal approximation theorem-type result to show that the proposed network can approximate the mapping resulting from forward model evaluations to a desired accuracy. Furthermore, we establish convergence of the surrogate posterior to the true posterior under Hellinger distance.  
\end{abstract}




\end{frontmatter}



\section{Introduction}
\label{Introduction}

\textcolor{black}{Markov chain Monte Carlo (MCMC) methods have been used for almost six decades and have become the standard method for learning Bayesian complex models since the early 1990s \cite{gelfand1990sampling}. These approaches are widely used across many scientific disciplines, including physics, biology, and economics, where they serve as powerful tools for sampling from complex probability distributions. MCMC methods based on the Metropolis-Hastings (MH) algorithm \cite{metropolis1953equation,hastings1970monte} construct Markov chains whose stationary distributions are the Bayesian posterior distributions \cite{robert2005monte,roberts2004general}. 
A key component of MCMC methods is the evaluation of likelihood functions, representing the likelihood of observing data given a particular set of model parameters. In Bayesian inference, for example, the likelihood plays a crucial role in updating prior beliefs to posterior distributions. The accuracy and speed of MCMC algorithms are strongly dependent on the accurate computation and computational cost of the likelihood function. Our proposed method is to substitute the likelihood function with a neural network in order to speed up excessively slow computations. Furthermore, we evaluate this proposed approach on three different kinds of inverse problems where the evaluation of the likelihood function is computationally expensive.} 

\subsection{Background}

Inverse problems are a class of mathematical modeling problems that aim to determine the unknown causes from known consequences. In contrast to direct problems, which predict outcomes based on known parameters and inputs, inverse problems work backward to determine the underlying parameters that generate observable data. These problems arise in many real-world applications, including medical imaging, geophysics, astronomy, oceanography, weather prediction, and non-destructive testing, among others 
\cite{kaipio2006statistical,tarantola2005inverse}. They are often challenging to solve because they can be ill-posed, meaning that solutions may be nonexistent, non-unique, or highly sensitive to data variations. To address these concerns, techniques such as regularization are used, making inverse problems a crucial area of research for obtaining meaningful information from complex systems. Many real-world inverse problems are governed by Partial Differential Equations (PDEs), with the system states described by PDE solutions. The properties of these systems, such as absorption coefficient, permeability, and thermal conductivity, are defined by model parameters that cannot be directly measured. Instead, these parameters are inferred from discrete and noisy observations of the states. Since inverse problems are generally ill-posed, solutions often rely on classical regularization theory \cite{engl1996regularization} or Bayesian inference \cite{stuart2010inverse}. The Bayesian technique provides a flexible framework for solving inverse problems by imposing a prior distribution on the parameters, which allows the incorporation of prior knowledge, which can be thought of as regularization by itself, e.g., see \cite{CE_18}. In recent decades, this technique has gained significant attention for its benefits \cite{kaipio2006statistical,stuart2010inverse}.
A Bayesian approach compared to classical regularization is different in that the Bayesian approach only requires continuity of the solution with respect to parameters, whereas classical regularization methods typically require computing some type of gradient, which can be problematic on its own. Furthermore, the Bayesian approach is simple to implement. 

The Bayesian inversion finds the probability distribution of input parameters using measured data (corrupted by noise) and other available knowledge. Samples from this (posterior) distribution are produced using MCMC methods. However, the formulation of Bayesian inverse problems (BIP) poses various issues, among these one is typically interested in first showing that the BIP of interest is well-posed and subsequently in providing theoretical guarantees for the Bayesian solution (i.e., the posterior density) to converge to the `truth.' From an implementational perspective, dealing with the discretized (finite), albeit very high-dimensional posterior distributions can be difficult due to the expensive-to-solve forward models and high-dimensional parameter spaces. As a result, direct sampling approaches, such as MCMC-based methods \cite{cotter2013mcmc,gelman1997weak,goodman2010ensemble} will incur excessive computation costs. 

Typical approaches to deal with these problems include \textit{(i)}. model reduction methods \cite{cui2015data, cui2016scalable,lieberman2010parameter,marzouk2009dimensionality}, which exploit the intrinsic low dimensionality; \textit{(ii)}. direct posterior approximation methods, such as Laplace approximation and variational inference \cite{bui2013computational,schillings2020convergence}; and \textit{(iii)}. surrogate modeling \cite{li2014adaptive,yan2017convergence,zhou2020adaptive,gao2023adaptive}, which substitutes the expensive model with a low-cost replacement.

Among the strategies listed above, surrogate modeling (for the forward model) offers an easy and principled approach to integrate deep-learning based methods into classical MCMC algorithms for efficiently accelerating the sampling of posterior distributions. Deep neural networks (DNN) have recently gained popularity in science and engineering as surrogate models due to their ability to approximate high-dimensional problems \cite{han2018solving,raissi2019physics,schwab2019deep,tripathy2018deep,zhu2018bayesian}. In general, DNN employs the ability of neural networks to build a quick-to-evaluate surrogate model to approximate the parameter-to-observation maps \cite{zhou2020adaptive,deveney2019deep,yan2021acceleration}.

Operator learning architectures like Fourier Neural Operators (FNOs) \cite{li2020fourier} and DeepONets \cite{lu2021learning} can describe complex models in infinite-dimensional domains as high-dimensional approximations. Therefore, they are potential surrogates, as described in \cite{cao2023residual,genzel2022solving}. Raoni{\'c} et al. \cite{raonic2024convolutional} proposed novel modifications of the Convolutional neural network (CNN) to enforce structure-preserving continuous-discrete equivalence and enable the genuine, alias-free learning of operators. The resulting architecture is called a Convolutional Neural Operator (CNO). However, employing approximate models directly may generate a discrepancy or modeling error, worsening an already ill-posed situation and resulting in a poor outcome. In this article, besides advocating for a fusion of deep-learning based methods with MCMC algorithms, we also provide asymptotic guarantees for the surrogate posterior to converge in an appropriate sense to the true posterior. Finally, we also note that learning the inversion operator directly in a non-Bayesian framework is also possible and has been studied in the literature before, \cite{mish_23,park20,abhishek2025}. We list our main contributions in this article below.

\subsection{Contribution}


\begin{itemize}

\item To our knowledge, this is the first time a neural operator (\textit{NO}) has been used to replace the forward model evaluator in BIP to estimate coefficients of PDEs.

\item From the theoretical point of view, the proposed network architecture replaces a mapping from function space to an operator space (e.g., Neumann-to-Dirichlet operators) in the cases of Electrical Impedance Tomography and Diffuse Optical Tomography. For Quantitative Photoacoustic tomography, we use the \textit{NO} as a surrogate model to replace the forward model between two function spaces. In addition, we analyze a universal approximation theorem-type result to show that the proposed network structure can approximate the respective forward maps in a suitable asymptotic sense. Following this, we also show that the surrogate posterior converges to the true posterior asymptotically in the Hellinger metric.

\item In our numerical experiments, we show that MCMC-Net offers a significant speed-up compared to when the likelihood is evaluated using a typical finite element solver. This is beneficial for practical applications.

\item Even though we use the proposed technique to accelerate MCMC in BIP, it also has broader implications. In particular, in any application where the posterior is explored by an MCMC based method and requires a computationally expensive likelihood function evaluation, one can effectively use deep learning based surrogate models to achieve substantial speed-up.
\end{itemize}

\subsection{Organization}
The rest of the paper is organized as follows. In section \ref{BayesianInverseProblems}, we describe Bayesian inverse problems, and in section \ref{imagingprobs}, we discuss the inverse problems considered for computational investigation.
In section \ref{OperatorLearning}, we delve into deep learning for operator approximation, establish universal approximation theorem-type results for forward operator learning in inverse problems, and present convergence results for the posterior distribution. Section \ref{NumericalExperiments} presents the numerical experiments. Finally, section \ref{conclusion} concludes the paper.

\section{Bayesian Inverse Problems}
\label{BayesianInverseProblems}

\subsection{Theoretical underpinnings of PDE-based inverse problems}

To illustrate the theoretical underpinnings of our approach, we consider a steady-state physical system governed by the following PDE: 
\begin{align}\label{generalpde}
\begin{cases}
    \mathcal{D}(u(x); q(x)) &= 0, \quad x \in \Omega, \\
    \mathcal{B}(u(x)) &= 0, \quad x \in \partial \Omega,
\end{cases}
\end{align}
where \( \mathcal{D} \) represents a general partial differential operator defined in a domain \( \Omega \subset \mathbb{R}^d \), \( \mathcal{B} \) denotes the boundary operator acting on the boundary \( \partial \Omega \), \( q \) signifies the unknown parameter, and \( u \) describes the state field of the system, see e.g. \cite{Nickl_book2}. Here, $\mathcal{X}$ represents the space of the unknown parameter \(q\), and $\mathcal{Y}$ is the space of the observed data \( y \), which consists of measurements collected from the system. We describe the relationship between the unknown parameter $q$ and the observed data $y$ as:
\begin{equation}\label{eq:2}
y = \mathcal{G}(q) + \eta\,,
\end{equation}
where $y$ represents the observed data, $\mathcal{G}: \mathcal{X} \to \mathcal{Y}$ is the forward operator that maps the unknown parameter $q$ to the measurements $y$. {{Typically, $\mathcal{Y}:=\mathbb{R}^k$ is some finite-dimensional space, and $\eta \sim \mathcal{N}(0, \Sigma)$ is the measurement noise, modeled as a Gaussian random vector (RV) with zero mean and covariance matrix $\Sigma$.}} This forward problem models the physical process that generates the observable data $y$ from the unknown parameter $q$. On the other hand, the inverse problem involves recovering $q$ from $y$. {Indeed, for the inverse problem at hand, as the observed data $y$ lies in a finite-dimensional space, while the parameter of interest $q$ lies in some appropriate (infinite-dimensional) function space, this inverse problem is severely ill-posed. Now we will describe the framework for Bayesian inversion developed in \cite{stuart2010inverse}. In the Bayesian approach, the parameter $q$ and the observed data $y$ are modeled as random variables. The Bayesian `solution' to such an inverse problem is then the posterior measure of the random variable $q|y$. To evaluate this posterior measure, we begin by placing a prior probability measure on the space of parameters, i.e., we assume $q \sim \mu_0$ where $\mu_0$ is some probability measure on the space $\mathcal{X}$. Let the noise be independent of $q$ and be distributed according to the Gaussian measure, $\eta \sim \mathbb{Q}_{0}$, whose distribution is given by the multivariate normal $\mathcal{N}(0,\Sigma)$. Assuming that the data is given according to the additive noise model \eqref{eq:2}, we can say that the random variable $y|q$ is distributed according to the measure $\mathbb{Q}_q$, whose distribution function is given by $\mathcal{N}(\mathcal{G}(q),\Sigma)$. Furthermore, in this case,  there exists a positive Radon-Nikdoym density given by $\frac{d \mathbb{Q}_q}{d\mathbb{Q}_0} (y)=\exp{(-\mathit{\Phi}(q; y))} $, where,  $\mathit{\Phi}: X\times Y\to \mathbb{R}$ is the `log-likelihood' function, which is also sometimes referred to as a `potential.' The likelihood term measures how well a given parameter $q$ describes the observed data $y$. Now we consider the following two product measures, $\nu_0=\mu_0\times \mathbb{Q}_0$ and $\nu=\mu_0\times \mathbb{Q}_q$ on the product space $X\times Y$. Then, we have the following analogue of Bayes' theorem on infinite-dimensional spaces:
\begin{theorem}{\cite[Theorem 14.]{dashti17}} \label{th:1.1}
 Let $\mathit{\Phi}:X\times Y\to \mathbb{R}$ be $\nu_{0}$ measurable and let $Z_y$ defined as $\int_{X}\exp(-\mathit{\Phi}(q; y)) d\mu_{0}:=Z_y>0$ for $\mathbb{Q}_0$ a.s. $y$, then the conditional distribution of $q|y$ denoted by $\mu^y$ exists under $\nu$. Furthermore, $\mu^y\ll \mu_0$ and \begin{align}
\frac{d\mu^y}{d\mu_0}(q)=\frac{1}{Z_y}\exp(-\mathit{\Phi}(q; y)).
\end{align}
\end{theorem}
 Theorem \ref{th:1.1} can be interpreted as a statement about the existence of a solution to a Bayesian inverse problem as it establishes conditions for the existence of the posterior density.}

Overall, the Bayesian approach updates the prior about $q$ based on the observed data, and Bayes' theorem helps us to combine the prior and likelihood to obtain the posterior knowledge about $q$.

\subsection{Accelerating MCMC by using Operator Network surrogates}

\begin{figure}[htbp]
    \centering
    \includegraphics[scale=.3]{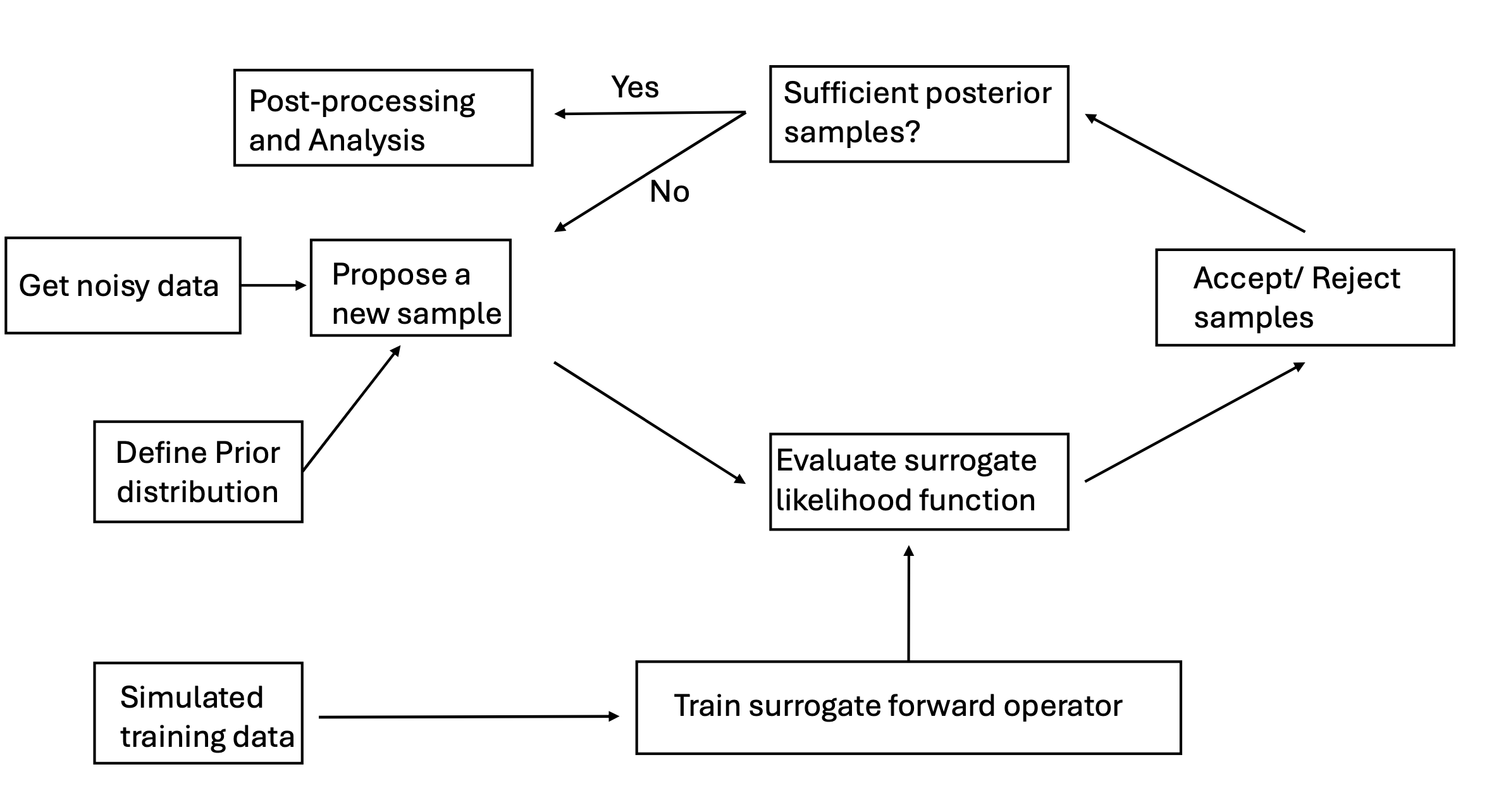}
    \caption{Flowchart illustrating the MCMC-Net workflow}
    \label{flow}
\end{figure}

For complex inverse problems of the kind described above, typically, any discretized representation of the posterior distribution is very high-dimensional and difficult to sample directly. \textit{MCMC} methods provide a solution by generating samples from the posterior distribution. In \textit{MCMC}, a Markov chain is constructed whose stationary distribution is the target posterior $ \mathsf{P}(q | y)$. The process involves iterating between proposing new samples and deciding whether to accept or reject them based on a criterion derived from the posterior distribution. One effective \textit{MCMC} technique for sampling in high-dimensional spaces is the \textit{preconditioned Crank--Nicolson (pCN) method}. This method proposes new samples by perturbing the current sample $q^{(t)}$ according to the following rule:
\begin{equation}
    q^{*} = q^{(t)} + \sqrt{1 - p^2} \, z,
\end{equation}
where $z \sim \mathcal{N}(0, I)$ is a standard Gaussian random variable and $p \in (0, 1)$ controls the amount of perturbation. The term $\sqrt{1 - p^2}$ ensures that the proposal step maintains stability while allowing sufficient exploration of the parameter space. After proposing a new sample $q^{*}$, the \textit{pCN} method evaluates the likelihood of the new sample using the acceptance probability:
\begin{equation}
    \alpha = \min\left( 1, \frac{\mathsf{P}(y | q^{*}) \mathsf{P}(q^{*})}{\mathsf{P}(y | q^{(t)}) \mathsf{P}(q^{(t)})} \right).
\end{equation}
If a uniform random variable $ u^{*} \sim \mathcal{U}(0, 1)$ is less than $\alpha$, the new sample is accepted ($q^{(t+1)} = q^{*}$); otherwise, the current sample is retained ($q^{(t+1)} = q^{(t)}$). By adjusting the parameter $p$, the algorithm can tune the exploration-exploitation trade-off, facilitating faster convergence to the target posterior distribution. It is important to note that each likelihood computation requires an evaluation of the forward model \( \mathcal{G} \), which can be computationally intensive, especially in complex real-world scenarios. To mitigate this, it is essential to replace the expensive forward model with a computationally inexpensive surrogate model. In this paper, we have discussed a deep learning-based surrogate to replace the traditional FEM solver, which {\textit{significantly}} reduces computational costs. We have tested our proposed technique on three inverse problems introduced in  {{section \ref{imagingprobs}}}. {{Here, we describe the essential idea behind replacing the exact forward model with a (learned) surrogate model and the relation of the `true' posterior with the `surrogate' (i.e., approximate) posterior. The idea behind surrogate modeling is to replace the forward model given by \eqref{eq:2} with the following equation (see also Figure \ref{flow}):
\begin{equation}\label{eq:2_theta}
y = \mathcal{G}_{\theta}(q) + \eta\,,
\end{equation}
where $\mathcal{G}_{\theta}$ indicates a neural-network based surrogate with $\theta$ denoting the network parameters. If $\mathcal{G}_{\theta}$ is a `good approximation of $\mathcal{G}$, then, using the surrogate model gives rise to a surrogate likelihood, $\mathit{\Phi}_{\theta}(q; y)$. Substituting the surrogate likelihood in lieu of the true likelihood gives rise to a surrogate posterior $\mu_{\theta}^y$, which satisfies (see also Theorem \ref{th:1.1}):
\begin{align} \label{eq:2.7}
\frac{d\mu_{\theta}^y}{d\mu_0}(q)=\frac{1}{Z^{\theta}_y}\exp(-\mathit{\Phi}_{\theta}(q; y))\,.
\end{align}
where $Z^{\theta}_{y}:=\int_{X}\exp(-\mathit{\Phi}_{\theta}(q; y)) d\mu_{0}$. In such cases, we would further like to understand the `closeness' of the surrogate posterior to the true posterior. One way to quantify this is to use the notion of Hellinger distance between the two posterior measures, which is defined by:
\begin{align}\label{eq:2.3}
\lvert \mu^y-\mu^{y}_{\theta}\rvert^2_{\text{Hell}}&={\frac{1}{2}}\int_{X} \bigg(\sqrt{\frac{\exp(-\mathit{\Phi}(q; y))}{Z(y)}}-\sqrt{\frac{\exp(-\mathit{\Phi}_{\theta}(q; y))}{Z^{\theta}(y)}} \bigg)^2 d\mu_0\,. 
\end{align}
In subsequent sections, we will make the following notion precise: {\textit{If $\mathcal{G}(q)\approxeq \mathcal{G}_{\theta}(q)$} in some appropriate sense, then $\lvert \mu^y-\mu^{y}_{\theta}\rvert^2_{\text{Hell}}$ is also small. This provides a principled rationale for replacing the FEM solver implementation of the forward model with a fast surrogate neural network implementation of the forward model instead.}}}

\section{The three imaging inverse problems.}\label{imagingprobs}
We will describe here the mathematical formulation of the three imaging problems we propose to study in this work, namely Electrical Impedance Tomography, Diffuse Optical Tomography, and Quantitative Photoacoustic Tomography.
\subsection{Electrical Impedance Tomography}
\label{subsection_eitmodel}
Electrical Impedance Tomography (EIT) is a non-invasive imaging technique that estimates the electrical conductivity distribution of an object by injecting currents and measuring the resulting boundary voltages via electrodes. EIT is widely used in medical imaging, geophysical exploration, and industrial monitoring applications. Usually, the forward model for EIT is described by a Complete Electrode Model (CEM) \cite{SoChIs_92}, which we briefly describe below. Let $\Omega\in \mathbb{R}^d, d\geq 2$ be the region to be imaged with its boundary denoted by $\partial \Omega$. Assume that there are $L$ electrodes, $\{e_l \}_{l=1}^L$, placed along $\partial \Omega$ where the boundary measurements of current and electric potential (voltage) can be made. Let $\sigma(x)$ be the conductivity coefficient, $u$ denotes the electric potential, $(I_l)_{l=1}^L$ be the applied current simulation pattern on the $L$ electrodes, and $(U_l)_{l=1}^L$ be the corresponding voltages across the electrodes. Let $z_l$ be the contact impedance of the \(l\)-th electrode. We assume that the parameter $\sigma(x)$ is a real valued, positive, bounded function in $L^{\infty}(\bar{\Omega})$. The CEM model corresponds to the following mathematical formulation, 
\begin{align} \label{eq:1}
-\nabla\cdot (\sigma(x)\nabla u(x)) &=0, \quad x\in \Omega\\
\int _{e_l} \sigma\dfrac{\partial u}{\partial n}dS &= I_l, \quad l\in \{1,\dots, L\}\\
\sigma(x)\frac{\partial u}{\partial n}(x)&=0, \quad \text{on } \partial \Omega\setminus \bigcup_{l=1}^{L}e_l\label{eq:2.3_bc}\\
u(x)+z_l \sigma(x)\frac{\partial u}{\partial n} &= U_l, \quad x\in e_l, \quad l\in \{1,\dots, L\}\label{eq:2.2}.
\end{align}
We aim to solve the inverse problem of recovering the conductivity field using boundary measurements of voltages, \((U_l)_{l=1}^L\), corresponding to a series of applied boundary currents, \((I_l)_{l=1}^L\). Assume that $J$ linearly independent current patterns are applied. For each $I^{(j)}$, $j\in \{1,\dots,J\}$, let us represent the non-noisy voltage across $L$ electrodes by $U^{(j)} (\in \Rb^L)$. Then, we can formally write,
\begin{align}
U^{(j)}={G}_j (\sigma),
\end{align}
where \(G_{(j)}(\sigma) = R(\sigma)I^{(j)}\), and \(R(\sigma) \in \mathbb{R}^{L \times L}\) is the resistivity matrix. Consider the case when the measured data is corrupted by a Gaussian noise.
Let \( y_j \) represent the noisy voltage measurements taken on the boundary:
\begin{equation}\label{eq:2.6}
y_j={G}_{j}(\sigma)+\eta _j, \quad j\in \{1,\dots,J\}\quad \text{and } \,\, \eta_j\sim N(0,\Gamma_0) \ i.i.d.
\end{equation}
Here $N(0,\Gamma_0)$ is used to denote a Gaussian random variable with mean ${0}$ and variance $\Gamma_0$. Concatenating all the vectors $y_j\in \mathbb{R}^{L}$ we can write:
\begin{align}\label{eitfdm}
y={G}(\sigma)+\eta
\end{align}
where $y\in \mathbb{R}^{JL}$ and $\eta \sim N(\textbf{0},\Gamma)$ where $\rm{\Gamma}=diag(\Gamma_0,\dots,\Gamma_0)$. In typical experiments, $y\in \Rb^{16\times 16}$ which corresponds to measurements taken across $16$ electrodes. The statistical inverse problem can now be formulated as a recovery of the parameter, $\sigma(x)$, from observed (noisy) data $y$.

In this work, we will work with a particular level-set formulation for the EIT inverse problem that was considered in \cite{stuart_16}. Thus, we will assume that the conductivity $\sigma$ can be expressed with the help of level-set functions. In particular, we assume that $\sigma(x)$ is piecewise constant and can be expressed by:
\begin{align}
    \sigma(x)=\sum_{i=1}^M \sigma_{i}\mathbb{I}(\Omega_{i})
\end{align} for some $M\in \Nb$. Here, $\mathbb{I}(S)$ denotes the characteristic function of some set $S\subset \Omega$. Also, for $i\neq k$, $\Omega_i \cap \Omega_k=\emptyset$. Besides, $\cup_{i=1}^M \Omega_i=\Omega$. The constants $\sigma_i $ are known, bounded, strictly positive numbers. Choose numbers $\{c_i\}_{i=1}^M$ with $c_0<c_1\dots<c_M$ and a continuous function, called a {\textit{level-set function}} $w:\Omega \to \Rb$ such that:
\begin{align*}
\Omega_i&=\{x\in \Omega: c_{i-1}\leq w(x)<c_i\}.
\end{align*}
Now let $H$ be an operator such that $H:H^s(\bar{\Omega})\to \mathcal{A}_{\sigma}(\Omega)$:
\begin{align*}
H(w)=\sum_{i=1}^M \sigma_{i}\mathbb{I}(\Omega_{i})=\sigma(x) 
\end{align*}
Here, for $s>1$, $H^s(\Omega)$ represents a space of Sobolev smooth functions, and $\mathcal{A}_{\sigma}(\Omega)$ denotes the space of admissible conductivities.
In terms of the level set function, equation \eqref{eitfdm} can be rewritten as:
\begin{align}
 y={G}(H(w))+\eta=:\mathcal{G}(w)+\eta 
\end{align}
where $\mathcal{G}=G\circ H.$
For more details about EIT and its level-set formulation in the Bayesian setting, see \cite{SoChIs_92,stuart_16,cheney1999electrical,borcea2002electrical,stuart16a} and the references therein.
Our goal will be to replace the true forward operator $\mathcal{G}$ (or, rather, its FEM implementation) with a neural network surrogate $\mathcal{G}_{\theta}$. We will describe this in {section \ref{OperatorLearning}}.

\subsection{Diffuse Optical Tomography}
Diffuse Optical Tomography (DOT) is an imaging technique that uses low-energy visible or near-infrared light to probe highly scattering media, including biological tissue. The forward problem in DOT is governed by the diffusion approximation of the radiative transfer equation, which models light transport in highly scattering media like biological tissue. It is often described by the following PDE:
\begin{equation}
\left.\begin{aligned}
-\nabla \cdot (\rho(x) \nabla u(x)) + \mu(x) u(x) &= h(x) \quad \text{in} \,\,\, \Omega\,, \\
u(x) + 2\rho(x) \frac{\partial u(x)}{\partial \nu} &= f(x) \quad \text{on} \,\,\, \partial \Omega\,,
\end{aligned}\right\}
\end{equation}
where \(u\) represents the photon density, \(\rho\) is the diffusion coefficient related to the scattering properties of the tissue, \(\mu\) is the absorption coefficient, and \(h\) represents the source term corresponding to the injected light. The boundary condition reflects the relationship between the photon density and its flux at the boundary \(\partial \Omega\), where \(f\) denotes the outgoing light on the surface. In DOT, the inverse problem involves reconstructing the tissue’s optical properties (absorption and diffusion coefficients) from photon density measurements on the boundary. DOT is widely used in medical imaging for applications like brain function monitoring, breast cancer detection, and tissue oxygenation assessment. In the present work, we consider a simplified one-parameter DOT model. The governing equation for this model is given by:
\begin{equation}
-\nabla \cdot (\rho(x) \nabla u(x)) + \mu (x) u (x) = 0 \quad \text{in } \Omega.
\end{equation}
For this simplified model, we focus solely on reconstructing the absorption coefficient \(\mu\) from boundary measurements of the photon density \(u\) while assuming that the diffusion coefficient \(\rho\) is known and fixed.  For more details about DOT and its discrete measurement model, see \cite{abhishek2022optimal, natt_wub, harrach2009uniqueness} and the references therein. We note that the inverse problem of EIT is a close cousin of the inverse problem of DOT, especially for the one-parameter DOT inverse problem considered here. For the sake of brevity, we do not present the formulation of the discrete problem for DOT here, except to mention that in our simulations, we do not use a level-set formulation for the priors but instead the formulation prescribed in \cite{abhishek2022optimal}. This is so that we can compare our reconstruction results against the ones shown in \cite{abhishek2022optimal}.

\subsection{Quantitative photoacoustic tomography}
We consider the inverse problem in Quantitative Photoacoustic Tomography (QPAT) as a third example. This imaging technique leverages the photoacoustic effect to assess the optical properties of target tissues. Light-absorbing molecules (chromophores) generate an acoustic pressure wave due to heat expansion by illuminating tissue with a short pulse of near-infrared light. The inverse problem reconstructs optical parameter images from acoustic pressure waves measured at the boundary of the tissue. QPAT has many applications in breast and skin cancer detection, small animal imaging, and vascular imaging.

In QPAT, light propagation through a scattering medium is often modeled using an elliptic PDE:
\begin{equation}
\left.\begin{aligned}
-\nabla \cdot (\rho(x) \nabla u(x)) + \gamma(x) u(x) &= 0 \quad \,\,\,\,\,\,\,\, \text{in } \Omega\,, \\
u(x)  &= g(x) \quad \text{on } \partial \Omega\,,
\end{aligned}\right\}
\end{equation}
where $\rho$ and $\gamma$ are the diffusion and absorption coefficients, respectively. The boundary condition \( u = g \) defines the incoming radiation source. In the present study, the statistical inverse problem in QPAT aims to reconstruct the parameter \( \gamma \), assuming \( \rho \in L_{\Lambda_\rho}^2(\Omega) \) is fixed, from the observed (noisy) data \( Y \) as described in \eqref{qpat_fdm1}.
For this problem, we consider an observation model as
\begin{align}\label{qpat_fdm1}
 Y = G(\gamma) + \varepsilon_n \xi\,,
\end{align}
where the continuous forward map \( G \) is considered as
\begin{align}\label{true_fdm_qpat}
G : \gamma \mapsto H := \gamma u, \quad G : L^2_\Lambda(\Omega) \to L^2(\Omega). \end{align}
Here the space of parameters, \(
L^2_\Lambda(\Omega) = \left\{ f \in L^2(\Omega) : \Lambda^{-1} \leq f \leq \Lambda \text{ a.e.} \right\}
\) and $\Lambda>0$ is a constant, $\mathcal{Y}$ denotes a real separable Hilbert space with an orthonormal basis \( \{ e_k \}_{k=1}^\infty \).  \( \xi \) denotes `white noise' in \( \mathcal{Y} \), which can be defined as 
\[
\xi := \sum_{k=1}^\infty \xi_k e_k, \quad \xi_k \overset{\text{i.i.d.}}{\sim} N(0, 1)\,.
\]
In the above, \( \varepsilon_n = \frac{\gamma}{\sqrt{n}} \) is used to denote the noise level, where \( \gamma > 0 \) and \( n \in \mathbb{N} \). The term \(Y\) is understood to depend on \(n\) and \(\gamma\). A typical discrete observation model (e.g., see \cite{afkham2024bayesian}) to the continuous observation model \eqref{qpat_fdm1} for the numerical experiments is as follows:
\begin{align}\label{qpat_apprmodel}
   Y_k = \langle G(\gamma), e_k \rangle_{L^2(D)} + \varepsilon \xi_k, \quad k = 1, \ldots, N_d\,,    
\end{align}
where \( \{ e_k \}_{k=1}^\infty \) represents an orthonormal basis of \( L^2(\Omega) \), consisting of the eigenfunctions of the Dirichlet Laplacian on \( \Omega \), and \( N_d \in \mathbb{N} \) is a suitable number. The observation \( \mathbf{Y} = \{ Y_k \}_{k=1}^{N_d} \) is
the collection of coefficients of the projection of \( Y \) from \eqref{qpat_fdm1} to the span of \( \{ e_k \}_{k=1}^{N_d} \). As \( N_d \to \infty \), observing \( \mathbf{Y} \) is equivalent to observing \( Y \); see, for example, [\cite{NiAb19}, Theorem 26]. For further details on QPAT and the discrete observation model \ref{qpat_apprmodel}, we refer to \cite{afkham2024bayesian, suhonen2024single, tarvainen2012reconstructing} and references therein.

For this problem, we consider star-shaped prior for parametrization of the inclusions, that is, inclusions parametrized by their center and a radial function. In numerical simulations, it was shown in \cite{afkham2024bayesian} that compared to level-set parametrization, star-shaped parametrization produced better reconstructions of the absorption parameter $\gamma$ in QPAT, and we present that formulation here. To begin, consider star-shaped sets in the plane. Let \( \psi \) be a continuously differentiable \( 2\pi \)-periodic function. We consider \( \psi: \mathbb{T} \to \mathbb{R} \) as a function defined on the one-dimensional torus \( \mathbb{T} := \mathbb{R} / 2\pi \mathbb{Z} \).
First, we introduce the parametrization for a single inclusion. The boundary $\partial A$ of the star-shaped set is a deformed unit circle: for a point \( x \in \Omega \), it can be defined as
\begin{align}\label{qpat_starpara_partial_A}
\partial A(\psi) = x + \{\exp(\psi(\vartheta)) \nu(\vartheta) \,, 0 \leq \vartheta \leq 2\pi \},
\end{align}
where \( \nu(\vartheta) := (\cos \vartheta, \sin \vartheta) \). The interior of the set is then given by
\begin{align}\label{qpat_starpara_A}
A(\psi) = x + \{s \exp(\psi(\vartheta)) \nu(\vartheta) \mid 0 \leq s \leq 1, \, 0 \leq \vartheta \leq 2\pi \}.    
\end{align}
Now, for positive constants \( \kappa_1, \kappa_2 > 0 \) we can define the map \(\Phi : \Omega \times H^\beta(\mathbb{T})  \to \mathcal{A}_{\gamma}(\Omega)\) by
\begin{align}\label{qpat_starpara_phi}
\Phi( x, \psi) := \kappa_1 \mathbb{I}_{A(x, \psi)} + \kappa_2\,.    
\end{align}
In the above, $A$ is the Lebesgue measurable subset of $\Omega$, and $\mathcal{A}_{\gamma}(\Omega)$ denotes the space of admissible absorption coefficients. For simplicity, we fix the center $x \in \Omega$ then \( \Phi(\cdot,\psi) \) satisfies a Hölder continuity-type result; see \cite[Lemma 4.1]{afkham2024bayesian}. Additionally, $\Phi(\cdot,\psi)$ represents a star-shaped parametrization for a single inclusion. We can generalize this for multiple inclusions. Consider \(N \geq 1\) disjoint inclusions, each parametrized by its respective center \( x^{(i)} \in \Omega \) and radial function \( \psi^{(i)} \). The boundary of the \( i \)-th inclusion is given as:
\[
\partial A^{(i)}(\psi^{(i)}) = x^{(i)} + \left\{\exp\big(\psi^{(i)}(\vartheta)\big) \nu(\vartheta) \mid 0 \leq \vartheta \leq 2\pi \right\},
\]
where \( \nu(\vartheta) = (\cos \vartheta, \sin \vartheta) \).
The interior of the \( i \)-th inclusion is:
\[
A^{(i)}(\psi^{(i)}) = x^{(i)} + \left\{s \exp\big(\psi^{(i)}(\vartheta)\big) \nu(\vartheta) \mid 0 \leq s \leq 1, \, 0 \leq \vartheta \leq 2\pi \right\}\,,
\]
where \( A_i(\psi) \) are some Lebesgue measurable subsets of \( \Omega \).
To define the composite conductivity map for all \( N \) inclusions, we consider the map \( \Phi: \big(\Omega \times H^\beta(\mathbb{T})\big)^N \to \mathcal{A}_{\gamma}(\Omega) \) as:
\begin{align}
\Phi\big(\{x^{(i)}\}_{i=1}^N, \{\psi^{(i)}\}_{i=1}^N\big) = \sum_{i=1}^N \kappa_i \mathbb{I}_{A^{(i)}(\psi^{(i)})},
\end{align}
where \( \kappa_i > 0 \) denotes the conductivity value in the \( i \)-th inclusion. If the background conductivity \( \kappa_{N+1} \) is included, the domain \( \Omega \) is partitioned such that:
\[
A^{(N+1)} = \Omega \setminus \bigcup_{i=1}^N A^{(i)}(\psi^{(i)}),
\]
and the map \( \Phi \) can be written as:
\begin{align}
\Phi\big(\{x^{(i)}\}_{i=1}^N, \{\psi^{(i)}\}_{i=1}^N\big) = \sum_{i=1}^{N+1} \kappa_i \mathbb{I}_{A^{(i)}(\psi^{(i)})}.
\end{align}
In this case, \( \Phi(\psi) \) also satisfies a Hölder continuity-type result; see \cite[Lemma 4.3]{afkham2024bayesian}, and \( \Phi(\psi) \) represents a star-shaped parametrization for multiple inclusions.

If we \textit{fix} the center of the inclusions in terms of the star-shaped set, the observation model
\eqref{qpat_fdm1} can be rewritten as:
\begin{align}
&Y = G(\Phi(\psi)) + \varepsilon_n \xi:=\mathcal{G}(\psi) + \varepsilon_n \xi\,.
\end{align}
Similarly, as the EIT, our goal is to replace the true forward operator $\mathcal{G}$ with a neural network surrogate $\mathcal{G}_{\theta}$. We will describe this in the next section.

\section{Deep Learning for Operator Approximation}\label{OperatorLearning}
Deep Learning for operator approximation is a cutting-edge approach that leverages the power of neural networks to approximate complex mathematical operators. These operators often arise in various scientific and engineering problems, such as solving PDEs, modeling dynamical systems, and simulating physical processes. 

For many decades, well-established numerical techniques like finite differences, finite elements, finite volumes, and spectral methods \cite{quarteroni2008numerical} have been successfully used to approximate PDE solutions. Nevertheless, these methods are more computationally expensive, particularly for many query issues, including uncertainty quantification (UQ), inverse problems, PDE-constrained control, and optimization, as well as in high-dimensional settings. Hence, reducing the computational cost demands fast, robust, and accurate surrogate models. Consequently, data-driven machine-learning algorithms have become a popular method for solving PDEs \cite{karniadakis2021physics}.

A very selected list of architectures for operator learning includes operator networks \cite{chen1995universal}, DeepONets \cite{lu2021learning} and its variants \cite{cai2021deepm,mao2021deepm}, PCA-net \cite{bhattacharya2021model}, neural operators \cite{kovachki2023neural} such as graph neural operator \cite{li2020neural}, multipole neural operator \cite{li2020multipole} and the very popular Fourier Neural Operator \cite{li2020fourier} and its variants \cite{li2024physics,pathak2022fourcastnet}, VIDON \cite{prasthofer2022variable}, the spectral neural operator \cite{fanaskov2023spectral}, LOCA \cite{kissas2022learning}, NOMAD \cite{seidman2022nomad}, De Hoop et al.\cite{de2022deep,de2023convergence}, Furuya et al. \cite{furuya2024globally}, and transformer-based operator learning architectures \cite{cao2021choose}. See also similar approaches considered for the advection-dispersion-reaction mechanisms considered in \cite{Uslu_24,Uslu_24a}.

Deep learning, mainly through architectures like neural operator networks, offers a scalable and efficient alternative by learning the underlying mappings from input functions to output solutions. This approach not only accelerates computations but also enhances accuracy and generalization, making it a promising tool for tackling challenging problems in computational science and engineering. Neural Operator Networks, such as FNOs \cite{li2020fourier} and DeepONets \cite{lu2021learning}, have shown significant promise in this field. 

The advantage of using deep learning for operator approximation lies in its ability to generalize from training data to unseen scenarios, enabling rapid predictions without requiring extensive re-computation. This is particularly beneficial in real-time applications and scenarios requiring repeated evaluations of operators, such as in MCMC-based algorithms for posterior exploration.

\subsection{Forward operator learning for EIT and QPAT}
With respect to the algorithmic implementation as described in section \ref{AlgorithmicImplementation}, we note that our proposed method for accelerating the MCMC computations is to replace the accurate but slow FEM implementation of the forward operator with a learned Neural-Operator approximation of the corresponding forward operator.
\begin{figure}[h]
\begin{center} 
\begin{tikzpicture}
\tikzset{boldarrow/.style={->, line width=1.5pt}}
\draw[boldarrow] (0.7,3.5) to node[above] {$ \mathcal{G} $} (7.2,3.5);
\draw node at (0.1,3.45) {$X$};
\draw node at (8.0,3.45) {$Y$};
\draw node at (4.0,-0.6) {$\mathbb{R}^{P}$};

\draw[boldarrow] (0.4,3.0) to node[left] {$\mathcal{E}_P$} (3.4,-0.4);
\draw[boldarrow] (3.6,-0.2) to node[right] {$\mathcal{D}_P$} (0.6,3.2);

\draw[boldarrow] (4.4,-0.2) to node[right] {$\mathcal{A}$} (7.3,3.2);
\end{tikzpicture}
\end{center}
\caption{The true map \( \mathcal{G} \) is approximated by a composition of two maps, encoder $\mathcal{E}$ and approximator $\mathcal{A}$.}\label{approxmapnew}
\end{figure}

 Note that in the Bayesian inverse problem formulated in section \ref{imagingprobs} for EIT, the space $H^s(\bar{\Omega}):=X$ is a separable Hilbert space, and $H: X\to \Gamma$ is a level set map from the space $X$ to the space of piecewise constant admissible conductivities. In this case, we denote the forward problem as:
\begin{align}
    y=\Gc_{\text{EIT}}(w)+\eta
\end{align}
where $\Gc_{\text{EIT}}:X\to \Rb^{L\times L}$ and $w$ is the level set representation of the unknown parameter. The goal is to approximate the map $\Gc_{\text{EIT}}$ by a neural operator $\Gc_{\theta_{\text{EIT}}}:X\to \Rb^{L\times L}$. Similarly, for QPAT, when using star-shaped parametrizations, we have a map $\Phi:  X\to L^2_{\Lambda}$ as described in \cite[section 4.1.1]{afkham2024bayesian}. As a result, we can again denote the forward problem for QPAT as:
\begin{align}
    y=\Gc_{\text{QPAT}}(w)+\eta
\end{align} where we now have, $y\in \mathbb{R}^{N_d}$, $\Gc_{\text{QPAT}}:X\to \mathbb{R}^{N_d}$, and $w$ is the star set representation of the unknown parameter. Again, the goal will; be to approximate the map $\Gc_{\text{QPAT}}$ by a neural operator $\Gc_{\theta_{\text{QPAT}}}:X\to \Rb^{N_d}$. In the context of the present study, a neural operator can be understood as a parametric mapping (depending upon the parameter $\theta$ which constitutes the various design choices of the network in terms of weights and biases) that connects an input space (such as the separable Hilbert space $X$ above) to an output space (like the space of matrices of size $\Rb^{L\times L}$ corresponding to the N-t-D data). Depending upon the context, i.e., whether the inverse problem pertains to EIT or QPAT, let us denote the true map between the input space $X$ and output space $Y$ by $\Gc$, then we will denote a neural-operator approximation of the true map by $\Gc_{\theta}$. The mapping $\Gc_{\theta}$ can be viewed as a composition of two distinct maps: an encoder and an approximator, shown in Figure \ref{approxmapnew}. In our notation, $\Gc_{\theta}=\mathcal{A}\circ\mathcal{E}_{P}$. Consequently,  the overall upper bound on the error arising from the approximation of the true operator by a neural operator can be decomposed into distinct upper bounds corresponding to the errors of the encoder and the approximator.

We now focus on the DeepONet framework introduced in \cite{lu2021learning,lanthaler2022error} and will adopt a similar concept for our present research. In this context, we define the following operators, maintaining the same terminology used in \cite{lanthaler2022error}.

\subsection*{Encoder}
We define the encoder \(\mathcal{E}_P\) as an operator that maps the function $w\in X$ to a set of discrete values \(\{w(\mathbf{x}_i)\}_{i=1}^{P}\) in \(\mathbb{R}^P\) at a set of regular grid points \(\{\mathbf{x}_i\}_{i=1}^{P} \in \Omega\). The encoder can be described as:
\[
\mathcal{E}_P : X \to \mathbb{R}^P, \quad \mathcal{E}_P ({w}) = (w_1, w_2, \dots, w_P):={\bf{\bar{w}}},
\]
where \(w_i = w(\mathbf{x}_i), \, i = 1, \dots, P\). When it is clear from the context, we will denote $\Ec_{P}:=\Ec$ for brevity.

\subsection*{Decoder}
While the encoder is a projection of an infinite-dimensional object (i.e., \(w \in X\)) into a finite-dimensional space \(\mathbb{R}^P\), the decoder lifts a finite-dimensional object back into the infinite-dimensional space \(X\). In this work, $\mathcal{D}_P(\bf{\bar{w}})$ is the unique trigonometric polynomial of degree $P$ such that, $\mathcal{D}_P({\bf{\bar{w}}})(x_j)=w(x_j)$. When it is clear from the context, we will denote $\Dc_P$ by $\Dc$ for brevity.
\subsection*{Approximator}
For EIT, the approximator \(\mathcal{A}_{\text{EIT}}\) maps the \(P\)-dimensional encoded vector in \(\mathbb{R}^P\) to an \(L \times L\) matrix in \(\mathbb{R}^{L \times L}\). This transformation can be described as:
\[
\mathcal{A}_{\text{EIT}} : \mathbb{R}^P \to \mathbb{R}^{L \times L}, \quad \mathcal{A}_{\text{EIT}}({\bf{\bar{w}}}) = \Lambda_P,
\]
where the vector ${\bf{\bar{w}}}=(w_1,w_2,\cdots,w_P)$ is the encoded representation of $w$, and \(\Lambda_P \in \mathbb{R}^{L \times L}\) is the output of the approximator map.
Using the above definitions and referring to \ref{approxmapnew}, note that our neural operator can be written in the form \(\mathcal{G}_{\theta_{\text{EIT}}} \equiv  \mathcal{A}_{\text{EIT}} \circ \mathcal{E}\), i.e., \(\mathcal{G}_{\theta} (w) =\mathcal{A}_{\text{EIT}} \circ \mathcal{E} (w)\). For QPAT, the approximator \(\mathcal{A}_{\text{QPAT}}\) maps the \(P\)-dimensional encoded vector in \(\mathbb{R}^P\) to a vector in $\mathbb{R}^{N_d}$. This transformation can be described as:
\[
\mathcal{A}_{\text{QPAT}} : \mathbb{R}^P \to \mathbb{R}^{N_d}, \quad \mathcal{A}_{\text{EIT}}({\bf{\bar{w}}}) = {\bf{Y}}=\{Y_k\}_{k=1}^{N_d},
\]
where we again denote a finite-dimensional encoding of $w$ by a vector ${\bf{\bar{w}}}=(w_1,w_2,\cdots,w_P)$. We note that this finite dimensional representation of $w$ and, in particular, the dimension $P$ may be different for QPAT and EIT, but we will denote both finite dimensional representations by $\mathbb{R}^P$, and it will be clear from the context whether we speak of EIT or QPAT.
Using the above definitions and referring to \ref{approxmapnew}, note that our neural operator for QPAT can be written in the form \(\mathcal{G}_{\theta_{\text{QPAT}}} \equiv  \mathcal{A}_{\text{QPAT}} \circ \mathcal{E}\), i.e., \(\mathcal{G}_{\theta_{\text{QPAT}}} (w) =\mathcal{A}_{\text{QPAT}} \circ \mathcal{E} (w)\).


\begin{theorem}\label{Th:1} Here and below, we will denote the true map between the input space of parameters $X$ and the respective finite-dimensional output space of observations $Y$ by $\Gc$, i.e. $\Gc\equiv \Gc_{\text{EIT}}$ for the EIT case and $\Gc\equiv \Gc_{\text{QPAT}}$ for the QPAT case. Furthermore, we will denote the corresponding neural operator approximation of the true map $\Gc$ by $\Gc_{\theta}$. Similarly, $Y=\mathbb{R}^{L \times L}$ for the EIT inverse problem and $Y=\Rb^{N_d}$ for QPAT inverse problem.
Let $\mu_0$ be a (prior) probability measure on $X$. We will show that for every \( \epsilon > 0 \), there exists a finite-dimensional spaces, $\mathbb{R}^{P}$ and  continuous maps
\[
\mathcal{E} : X\to \mathbb{R}^P, \quad \text{and} \quad \mathcal{A} : \mathbb{R}^P \to Y, 
\]
such that
\[
\| \mathcal{G}(w) - \mathcal{G}_{\theta} (w) \|_{L^2{(\mu)}} =\bigg(\int_{X}\|\mathcal{G}(w) - \mathcal{G}_{\theta} (w) \|^2_{*}d\mu(w)\bigg)^{1/2}\leq \epsilon.
\]
where $\lVert \cdot\rVert_{*}$ represents the corresponding Frobenius norm in $\mathbb{R}^{L\times L}$for EIT inverse problem and it denotes the $l^2$ norm in $\mathbb{R}^{N_{d}}$  for the QPAT inverse problem.
\end{theorem}

\noindent \textit{Proof:} 
Let $\mu_0$ be a Gaussian measure on $X$.
If the neural operator architecture is given by \(\mathcal{G}_{\theta} \equiv \mathcal{A} \circ \mathcal{E}\), then the error \scalebox{1.5}{\(\hat{\bm{\varepsilon}}\)} in representing the true opertaor $\Gc(w)$ by its approximation $\Gc_{\theta}(w)$ measured in the $L^2(\mu)$-norm is given by:
\begin{align}
\scalebox{1.5}{\(\hat{\bm{\varepsilon}}\)} &= \left(\int_{X}   \| \mathcal{G}(w)-\mathcal{G}_{\theta} (w)\|^{2}_{*} d \mu (w) \right)^{1/2}
\end{align}
We can split the term $\| \mathcal{G}(w)  - \mathcal{G}_{\theta}(w) \|_{*}$
in the following way:

\begin{align}\label{cnoneterror}
\| \mathcal{G}(w)  -  \mathcal{A} \circ \mathcal{E}_P (w) \|_{*}
&\quad\leq \left.  \underbrace{ \| \mathcal{G}(w) - \mathcal{G} \circ \mathcal{D}_P \circ \mathcal{E}_P(w) \|_{*} }_{T_1} \right. \notag\\
&\quad \big. + \underbrace{ \| \mathcal{G}\circ \mathcal{D}_P \circ \mathcal{E}_P (w) - \mathcal{A} \circ \mathcal{E}_P (w) \|_{*} }_{T_2} \big. \Big)
\end{align}
where $T_1$ and $T_2$ will contribute respectively to encoding error and approximation error. Let us first consider the term $T_1$
\begin{align}
 T_1=\| \mathcal{G}(w) - \mathcal{G} \circ \mathcal{D}_P \circ \mathcal{E}_P(w) \|_{*}
\end{align}
In the Bayesian formulations of the EIT/QPAT inverse problem, whether we use level-set or star-shaped priors, the forward map $\mathcal{G}$ is a continuous and bounded map, see e.g. \cite{stuart_16,stuart16a,afkham2024bayesian}. More precisely, $\lVert\Gc(w)\rVert_{*}\leq M$, ($M$ is a fixed constant) see e.g. proof of \cite[Proposition 3.10]{stuart_16} or \cite[section 3]{afkham2024bayesian}. Furthermore, as $X$ is a separable Banach space, by Lusin's theorem, there exists a compact set $K$ such that $\mu_0(X\setminus K)<\frac{\epsilon}{8 M^2}$, see \cite[Lemma C.1]{lanthaler2022error}. Clearly, the map $\Gc$ is uniformly continuous on $K$. We have defined the maps $\mathcal{E}_P$ and $\mathcal{D}_P$ as the ones used in trigonometric interpolation, i.e., $ \mathcal{D}_P \circ \mathcal{E}_P:=\mathcal{I}_P$ is the pseudo-spectral projection (trigonometric interpolation) as described in [\cite{kovachki2021universal}, page 32]. From Theorems \textbf{39} and \textbf{40} in \cite{kovachki2021universal}, for any $w \in X$, we have 
\begin{align*}
    \| (Id - \mathcal{I}_P) w\|_{L^{\infty}} \lesssim P^{-\xi(s)} \quad \text{for some}\  \xi(s) >0\,.
\end{align*}
Now, consider the set $Z=K\cup_{P=1}^{\infty} \mathcal{I}_P(K)$. By \cite{kovachki2022machine}[Lemma 21], $Z$ is a compact set. Hence, there exists a modulus of continuity $\omega$ such that, $\lVert \Gc(w_1)-\Gc(w_2)\rVert \leq \omega\big(\lVert(w_1-w_2)\rVert_{X}\big)$ for $w_1,w_2 \in Z$.
\noindent
Here and below, the symbol $\lesssim$ will indicate that the inequality holds up to some constant.
Thus,
\begin{align}
\lVert T_1\rVert_{L^2(\mu)}=&\int_{K}\|  \mathcal{G} (Id - \mathcal{I}_P)w\|^2_{*} d\mu(w) +   \int_{X\setminus K}\|\mathcal{G} (Id - \mathcal{I}_P)w\|^2_{*}d\mu(w)\nonumber\\
&\leq \int_{Z}\|  \mathcal{G} (Id - \mathcal{I}_P)w\|^2_{*} d\mu(w) +   \frac{\epsilon}{8M^2}4M^2\nonumber\\
&\leq \omega(P^{-\xi(s)}) +   \frac{\epsilon}{2} \lesssim \epsilon
\end{align} 
for $P$ large enough.
Now, we look at
\begin{align}\label{eq:23}
    \| T_2 \|_{L^2(\mu)} &\leq \int_{\mathcal{E}(X)} \| \mathcal{A}(\hat{w}) - \mathcal{G} \circ \mathcal{D}_P(\hat{w}) \|^{2}_{*} d (\mathcal{E}_{\#} \mu) (\hat{w}) \\
    & =\| \mathcal{A}(\hat{w} )- \mathcal{G} \circ \mathcal{D}_P(\hat{w}) \|^{2}_{L^2(\mathcal{E}_{\#} \mu)}\,.
\end{align}

Recall that $\Gc$ is a continuous map and $\Dc_P$ is Lipschitz continuous. As a result, the composition $\Gc\circ\Dc_P$ is also continuous. For any given $P$ (suitably chosen so that eq. \eqref{eq:23} is satisfied), any continuous map can be well approximated by some ReLu DNN (see, e.g., \cite[Theorem 2]{yar18a}) or a deep CNN \cite[Theorem 1]{ZHOU2020} $\Ac$. As a result, we can make $T_2\lesssim \epsilon$. This concludes the proof.
\begin{theorem}
    Let us assume that the noise, $\eta$, is i.i.d. Gaussian. Let the true posterior measure, when using the true forward model, $\Gc$, be $\mu$ and the approximate posterior measure when using the Neural network surrogate, $\Gc_{\theta}$, be $\mu_{\theta}$.
    Then $\lvert\mu^y-\mu^y_{\theta}\rvert_{\text{Hell}} \to 0$ as $\epsilon \to 0$ in Theorem \ref{Th:1}, where $\lvert \cdot \rvert_{Hell}$ denotes the Hellinger distance between the posterior densities.
\end{theorem}
\noindent \textit{Proof:} Define the potential (likelihood) functions, $\mathit{\Phi}(w;y)$ and $\mathit{\Phi}_{\theta}(w;y)$ for the true and approximate forward models respectively as:
\begin{align}
    \mathit{\Phi}(w;y)&=\frac{1}{2}\lVert y- G(w)\rVert_{*}^2 \\
    \mathit{\Phi}_{\theta}(w;y)&=\frac{1}{2}\lVert y- G_{\theta}(w)\rVert_{*}^2.
    \end{align}
Then following similar calculations as in \cite[Lemma 4.2]{marzouk_09} and \cite[Appendix 1]{stuart16a}, we get
\begin{align}\label{eq:28}
    \lvert Z^{\theta}_y- Z_y\rvert\lesssim \lVert \Gc(w)-\Gc_{\theta}(w)\rVert_{L^2(\mu)}.
\end{align}
where $Z^{\theta}_y$ and $Z_y$ are the corresponding evidence terms as per \eqref{eq:2.7}.
Thus for the Hellinger distance between the two posterior measures $\mu^y$ and $\mu^y_{\theta}$, {see e.g. Theorem \ref{th:1.1}} or, \cite{stuart16a}, we have
\begin{align}\label{eq:5.13}
2\lvert \mu^y-\mu^y_{\theta}\rvert_{\text{Hell}}^2&=\int_{X} \left(\sqrt{\frac{\exp(-\mathit{\Phi}(w;y))}{Z_y}}-\sqrt{\frac{\exp(-\mathit{\Phi}_{\theta}(w;y))}{Z^{\theta}_y}} \right)^2 d\mu_0 \nonumber\\
&\leq I_1+ I_2
\end{align}
where \begin{align*}
I_1&=\frac{2}{Z_y}\int_{X} \bigg(\sqrt{\exp(-\mathit{\Phi}(w;y))}-\sqrt{\exp(-\mathit{\Phi}_{\theta}(w;y))} \bigg)^2 d\mu_0 \quad \quad \text{and}\\
I_2&={2}{\lvert Z_y^{-1/2}-{Z_y^{\theta}}^{-1/2}\rvert^2}\int_{X}\exp(-\mathit{\Phi}_{\theta}(w;y))d\mu_0
\end{align*}

We note that $Z_y$ is bounded below, and thus
\begin{align}
I_1 &= \frac{2}{Z_y}\int_{X}\bigg({\exp(-\frac{1}{2}\mathit{\Phi}(w;y))}-{\exp(-\frac{1}{2}\mathit{\Phi}_{\theta}(w;y))} \bigg)^2 d\mu_0 \nonumber\\
&\lesssim \lVert \Gc(w)-\Gc_{\theta}(w)\rVert_{L^2(\mu)} \quad \textnormal{(similar to calculations in \eqref{eq:28} and \cite[Lemma~4.2]{marzouk_09})}\,. 
\end{align}

In addition, following \cite[Appendix 1]{stuart16a}, 
\begin{align}
I_2&\lesssim {\lvert {Z_y}^{-\frac{1}{2}}-{{Z^{\theta}_y}^{-\frac{1}{2}}}\rvert^2}\int_{X}\exp(-\mathit{\Phi}_{\theta}(w;y))d\mu_0 \quad \nonumber \\
&\lesssim \lVert \Gc(w)-\Gc_{\theta}(w)\rVert^2_{L^2(\mu)}\label{eq:5.14}
\end{align}
where we use the fact that 
\begin{align}
    \mathit{\Phi}_{\theta}(w;y))\leq \lVert \Gc_{\theta}(w)\rVert^2+\lVert y\rVert^2\leq \lVert \Gc_{\theta}(w)-\Gc(w)\rVert^2+\lVert \Gc(w)\rVert^2+\lVert y\rVert^2.
\end{align}
is bounded if $y$ is bounded.
Thus as $\lVert \Gc(w)-\Gc_{\theta}(w)\rVert_{L^2(\mu)} \to 0$, then both $I_{1}\to0$ and $I_2\to 0$ and this concludes the proof.

\section{Numerical Experiments}
\label{NumericalExperiments}
In this section, we discuss the computational method and inversion results to demonstrate the effectiveness and accuracy of the MCMC-Net for accelerating the MCMC method on slow likelihood functions. We start with a detailed overview of the algorithm, model training, and the neural network architecture, specifically the CNN architecture, designed for various inverse problems. In our approach, we replace the traditional forward solver, such as FEM, with a trained CNN and refer to this framework as MCMC-Net. Finally, we compare the inversion results obtained using MCMC-Net with those from the classical MCMC-FEM approach, which does not involve any neural network components.


\subsection{Algorithmic Implementation}
\label{AlgorithmicImplementation}
Here, we describe the MCMC algorithm used in the numerical simulations. The central concept is to draw samples from the posterior distribution to compute the Bayesian estimate for the unknown parameter $q$ introduced in {section \ref{BayesianInverseProblems}}. This estimate can be approximated by the Monte Carlo average  \( E_\Pi(q^* \mid \Lambda^*_{\text{true-noisy}}) \approx \frac{1}{n} \sum_{i=1}^n q^*_i \), where \(q^*\) represents the discretized random variable of \(q\), \(q^*_i\) are the individual discrete samples of \(q^*\), and \(\Lambda^*\) corresponds to a finite set of measurements. We utilize the \textit{pCN} algorithm with a Gaussian process (GP) prior to sample from the posterior distribution. To begin, we define the log-likelihood function for the inverse problem for each sample \(q^*_i\), where \(i \in \{1, 2, \dots, N\}\):
\[
L(q^*_i) := -\frac{1}{2\sigma_{\text{noise}}^2} \|\Lambda^*_{\text{true-noisy}} - \Lambda^*_{q^*_i}\|_F^2,
\]
where \(\Lambda^*_{q^*_i}\) represents the solution to a single forward problem, and \(\|\cdot\|_F^2\) denotes the squared Frobenius norm.

Let \(\Pi\) denote a Gaussian prior with mean zero and covariance matrix \(\mathcal{C}\) for \(q^*\). Set the initial value \(q_0^*\) to the typical background value of \(q^*\). Then, repeat the following steps until the required number of samples is obtained:
\begin{enumerate}
\item Sample \(\Upsilon \sim \Pi\) and compute the proposal \(q^*_{\text{Prop.}} := \sqrt{1 - 2\Delta} q_i^* + \sqrt{2\Delta} \Upsilon\), where \(\Delta > 0\).
\item Update \(q^*_{i+1}\) as follows:
\begin{align*}
   q^*_{i+1} = \begin{cases}
   q^*_{\text{Prop.}} & \text{with probability } \min\left(1, L(q^*_{\text{Prop.}}) - L(q^*_i)\right), \\
   q_i^* & \text{otherwise}.
   \end{cases}    
\end{align*}
\end{enumerate}

It is important to note that a substantial number of burn-in iterations is required before reaching the high-probability regions of the posterior distribution. Only then can the samples be used to compute the Monte Carlo average. To ensure faster convergence, it is recommended to select a parameter \(\Delta > 0\), which is used to scale the covariance matrix of the proposal distribution to ensure that approximately 25\% of the proposed samples are accepted after burn-in. This is done to keep it in line with generally accepted optimal acceptance ratio of 0.234, see \cite{roberts98}. In practice, this problem is not trivial. Therefore, we adjust \(\Delta\) during the burn-in phase to maintain the acceptance rate near 0.25. This approach follows the method outlined in \cite{NiAb19} or a similar strategy proposed in \cite{ahmad2019comparison,khan_strauss}. Recall that the parameter of interest \(q\) is assumed to lie in an infinite-dimensional Sobolev space. In a computational setting, however, we approximate the parameter space by a finite, albeit very high-dimensional space, where functions are piecewise constant over each triangle in the finite element discretization. Traditional MCMC methods applied to such high-dimensional spaces suffer from the ``curse of dimensionality," which can lead to a large number of iterations being required for convergence. As discussed in \cite[Section 4]{NiAb19}, it can be shown that \(q^*_{i+1}\) forms a Markov chain. For more information on how the \textit{pCN} method enhances the algorithm and leads to significant computational speed-ups in Bayesian estimation, we direct readers to \cite{cotter2013mcmc}.

To achieve properly regularized reconstructions, we select the covariance matrix for \(\Pi\) to be defined by a Matérn kernel, with the parameters \(\upsilon \geq 3\) and \(\ell\) chosen heuristically. This Gaussian process prior imposes a regularization effect similar to a Sobolev-norm penalty, as discussed in \cite{NiAb19}. The Matérn kernel is given by
\[
k_{\upsilon, \ell}(d) := \frac{2^{1-\upsilon}}{\Gamma(\upsilon)} \left( \frac{d \sqrt{2 \upsilon}}{\ell} \right)^\upsilon K_\upsilon \left( \frac{d \sqrt{2\upsilon}}{\ell} \right),
\]
where \(d:= \|x_i - x_j\|_2\) is the Euclidean distance between two PDE mesh centroids, and \(K_\upsilon\) represents the modified Bessel function of the second kind. Reconstructions were done in this section by running the MCMC method with 50,000 samples for burn-in and another 50,000 samples that we considered from the true posterior distribution samples. Note that for high quality reconstructions and uncertainty quantification in such ill-posed inverse problems, one needs far more number of MCMC iterations so that one can reliably draw from the stationary posterior distribution. In our set-up with multiple experiments, this would cause a huge computational burden, especially when working with the traditional MCMC approach where the forward problem is realized by solving a PDE with an FEM solver. This is why we limit the MCMC iterations to 100,000 and compare the reconstruction quality of MCMC-Net with MCMC-FEM by running the sampler for an equal number of iterations. The reconstructed images show, that even with these many iterations, qualitatively the reconstructions are good. At the same time, we are able to show the relative advantage of the proposed method in terms of computational speed vis-a-vis the traditional method, while also outperforming the traditional method in various error metrics, including the $L^{\infty}$-loss, mean absolute error (MAE), and mean square error (MSE). This gives credence to the idea of using neural-net surrogates in Bayesian inversion, especially in light of the theoretical results presented in section \ref{OperatorLearning}. We also present an example reconstruction for EIT inversion in the appendix, where the MCMC-net was allowed to run for 200,000 iterations, see Figure \ref{EITReconstructionAN1Loss}. Convergence plots to illustrate convergence of the MCMC chain and good mixing properties have been shown in Figures \ref{EITan1Loss} and \ref{ConvergencediagnosticsEIT}

\subsection{Hardware}

We trained all the different neural networks used in this study on a high-performance computing (HPC) cluster utilizing the Slurm workload manager. Each training job was submitted to the GPU partition and allocated one node equipped with Dual 8-core Intel Xeon CPUs @ 3.20GHz and 1 NVIDIA Tesla V100S GPU with 32 GB HBM2 RAM. MATLAB R2024b was used for all training processes, with the necessary scripts executed in non-interactive mode. The same computing power was used during the EIT and DOT inversion (reconstruction) processes. For the QPAT inversion, all the computations were performed on a desktop running Windows 11 and MATLAB version R2024b. The hardware configuration includes an Intel Core i7 processor with 14 cores clocked at 2.30 GHz, 16 GB of RAM, and a dedicated NVIDIA GeForce RTX 3050 Ti GPU with 4 GB of VRAM.

\subsection{Data Generation and Training}\label{sec:5.3}
Since the inverse problems considered—EIT, DOT, and QPAT—differ in nature, we employ three distinct CNN architectures to replace the forward model evaluations in each respective problem. For training these CNNs, we generate datasets comprising input-output pairs: inputs represent problem-specific parameters (e.g., conductivity for EIT), and outputs correspond to measurements (e.g., boundary voltage measurements for EIT). Detailed descriptions of the data generation process, CNN architectures, and training methodologies for each inverse problem are provided in the subsequent sections. We systematically experimented with different CNN networks and training parameters to find the best setup (in terms of error and computation speed) for each problem. {Finally, we would also like to mention that the MCMC-FEM inversion code for QPAT is adapted from the github repository \cite{afkham_git}} and is based on the article \cite{afkham2024bayesian}.

\textbf{\textit{EIT:}} In the case of EIT, we generate a set of data pairs for conductivity distribution $\sigma$  and boundary measurements $\Phi_b$  for training the CNN. The dataset includes conductivity fields $\sigma$ with varying radii and locations of anomalies within the finite element mesh. Also, we generate $\sigma$ using the level set prior as described in section \ref{imagingprobs}. For the prior knowledge, we use the probability distribution with mean zero and Matérn kernel ($\upsilon=3$, $l=0.4$) as a covariance matrix. In total, we compute 7200 conductivity ($\sigma$) values. Then based on all these $\sigma$ values, we solve the forward EIT model using FEM to compute the associated boundary measurements $\Phi_b$. These paired data $\{\sigma,  \Phi_b\}$ are then used to train the CNN.
\begin{table}[ht]
\centering
\caption{EIT Forward Network Architecture}
\label{EITNet}
\scalebox{0.95}{
\begin{tabular}[t]{|l|c|c|c|c|}
\hline
\textbf{Network Block} & \textbf{Input Size} & \textbf{Output Size} & \textbf{Kernel Size} & \textbf{Activation Function} \\
\hline
Fully Connected & \(1345 \times 1\) & \(16 \times 16\) & – & ReLU \\
Convolutional   & \([16, 16]\)    & \([16, 16]\)    & \([3, 3]\) & ReLU \\
Convolutional   & \([16, 16]\)    & \([16, 16]\)    & \([3, 3]\) & ReLU \\
Convolutional   & \([16, 16]\)    & \([16, 16]\)    & \([3, 3]\) & ReLU \\
Convolutional   & \([16, 16]\)    & \([16, 16]\)    & \([3, 3]\) & ReLU \\
\hline
\end{tabular}
}
\end{table}

In this problem, the CNN architecture consists of a fully connected layer followed by a sequence of convolutional layers, as detailed in Table~\ref{EITNet}. Padding is applied to preserve the spatial dimensions throughout the convolutional blocks. The network is trained for 2000 epochs using a mini-batch size of 128 and a learning rate of 0.001.

\textbf{\textit{DOT:}} Similar to the EIT case, for the DOT inverse problem, we generate 6400 data pairs consisting of absorption coefficient distributions $\mu$ and the corresponding boundary measurements $\Phi_b$. The dataset includes samples of $\mu$ with varying anomaly radii and locations within the finite element mesh. The samples are drawn from a zero-mean Gaussian distribution with a Matérn covariance kernel ($\upsilon=3$, $l=0.2$). Unlike the EIT case, the level set parametrization is not used here.
\begin{table}[ht]
\centering
\caption{DOT Forward Network Architecture}
\label{DOTNet}
\scalebox{0.95}{
\begin{tabular}[t]{|l|c|c|c|c|}
\hline
\textbf{Network Block} & \textbf{Input Size} & \textbf{Output Size} & \textbf{Kernel Size} & \textbf{Activation Function} \\
\hline
Fully Connected & \(549 \times 1\) & \(16 \times 16\) & – & ReLU \\
Convolutional   & \([16, 16]\)    & \([16, 16]\)    & \([3, 3]\) & ReLU \\
Convolutional   & \([16, 16]\)    & \([16, 16]\)    & \([3, 3]\) & ReLU \\
Convolutional   & \([16, 16]\)    & \([16, 16]\)    & \([3, 3]\) & ReLU \\
Convolutional   & \([16, 16]\)    & \([16, 16]\)    & \([3, 3]\) & ReLU \\
\hline
\end{tabular}
}
\end{table}
In this problem, the CNN architecture consists of a fully connected layer followed by a sequence of convolutional layers, as shown in Table~\ref{DOTNet}. Padding is applied in all convolutional layers to preserve the spatial dimensions. The network is trained for 100 epochs using a mini-batch size of 8. The initial learning rate is set to 0.001 and is reduced by a factor of 0.1 every 50 epochs to facilitate better convergence.

\textbf{\textit{QPAT:}} For the QPAT inverse problem, we generate a dataset of 4800 paired samples consisting of absorption coefficients $\gamma$ and the corresponding measurements $\Phi_q$. The dataset includes circular-shaped $\gamma$ profiles with varying radii, as well as shapes resembling the simulated ground truth absorption coefficient shown in Figure~\ref{qpat_gamma_true}. To introduce diversity, the two inclusions in the ground truth configuration are independently rotated by random angles within the range $[-\pi, \pi]$. Additionally, star-shaped priors are used to generate further variations of $\gamma$, which are also included in the dataset. For each $\gamma$ sample, the forward QPAT model is solved using FEM to compute the corresponding measurement $\Phi_q$. These paired data $\{\gamma,  \Phi_q\}$ are then used to train the CNN.
\begin{table}[ht]
\centering
\caption{QPAT Forward Network Architecture}
\label{QPATNet}
\scalebox{0.95}{
\begin{tabular}[t]{|l|c|c|c|c|}
\hline
\textbf{Network Block} & \textbf{Input Size} & \textbf{Output Size} & \textbf{Kernel Size} & \textbf{Activation Function} \\
\hline
Fully Connected & \(2904 \times 1\) & \(16 \times 16\) & – & ReLU \\
Convolutional   & \([16, 16]\)    & \([16, 16]\)    & \([3, 3]\) & ReLU \\
Convolutional   & \([16, 16]\)    & \([16, 16]\)    & \([3, 3]\) & ReLU \\
Convolutional   & \([16, 16]\)    & \([16, 16]\)    & \([3, 3]\) & ReLU \\
Convolutional   & \([16, 16]\)    & \([16, 16]\)    & \([3, 3]\) & ReLU \\
\hline
\end{tabular}
}
\end{table}
In this problem, the CNN architecture consists of a fully connected layer followed by a sequence of convolutional layers, as detailed in Table~\ref{QPATNet}. Padding is applied in all convolutional layers to preserve spatial dimensions. The network is trained for 100 epochs using a mini-batch size of 8. The initial learning rate is set to 0.001, and a reduction factor of 0.1 is applied every 20 epochs to refine the training process.

In addition to the methodologies described above, we train all three CNN architectures using the ADAM optimizer. The results for EIT, DOT, and QPAT are presented in the following three subsections. We also emphasize the critical role of training data in determining the quality of the reconstructions across all cases. A significant finding from our numerical experiments is the generalization capability of the forward neural-operator surrogate—it was able to adequately recover geometric shapes that were not included in the training set, as illustrated in Figure~\ref{DOTConductivityReconstructionsUnseenData}.

\subsection{Electrical Impedance Tomography}
\begin{figure}
\centering

\begin{minipage}{.235\textwidth}
\begin{subfigure}{\textwidth}
\includegraphics[scale=0.32]{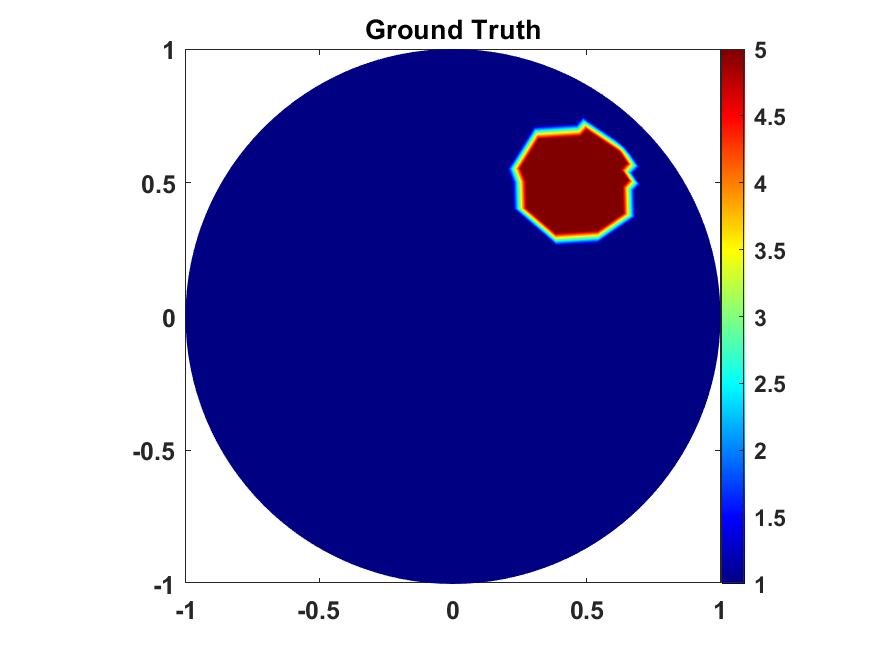}
\subcaption{True}
\label{eit_sigma_an1_true}
\end{subfigure}
\end{minipage}
\hfill
\begin{minipage}{.235\textwidth}
\begin{subfigure}{\textwidth}
\includegraphics[scale=0.32]{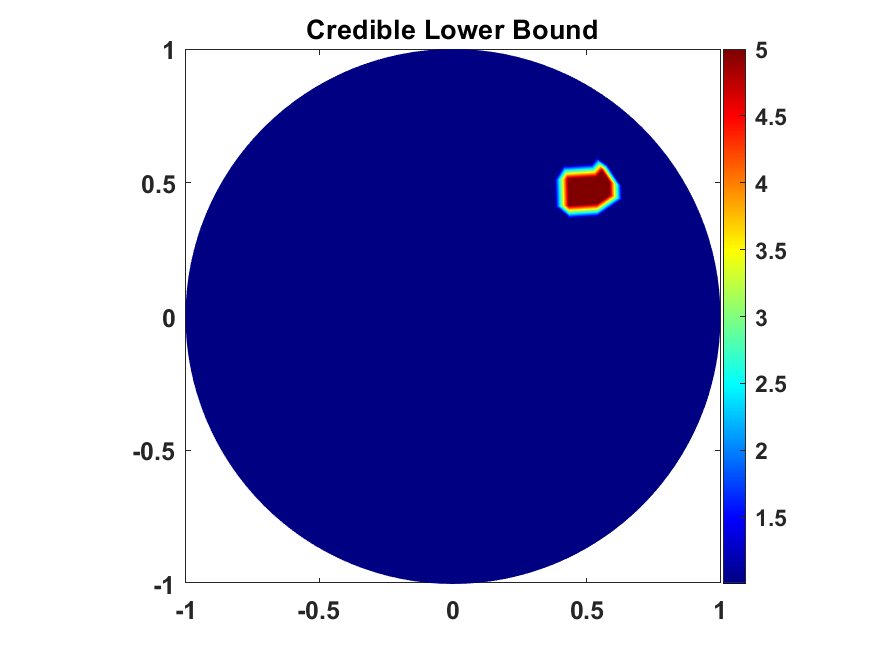}
\subcaption{FEM}
\end{subfigure}
\begin{subfigure}{\textwidth}
\includegraphics[scale=0.32]{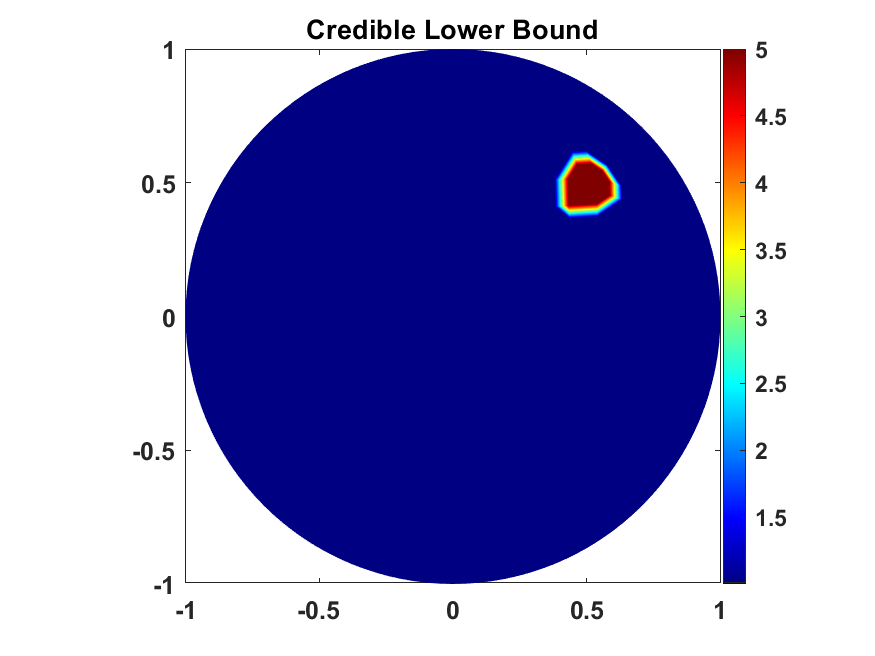}
\subcaption{Net}
\end{subfigure}
\end{minipage}
\hfill
\begin{minipage}{.235\textwidth}
\begin{subfigure}{\textwidth}
\includegraphics[scale=0.32]{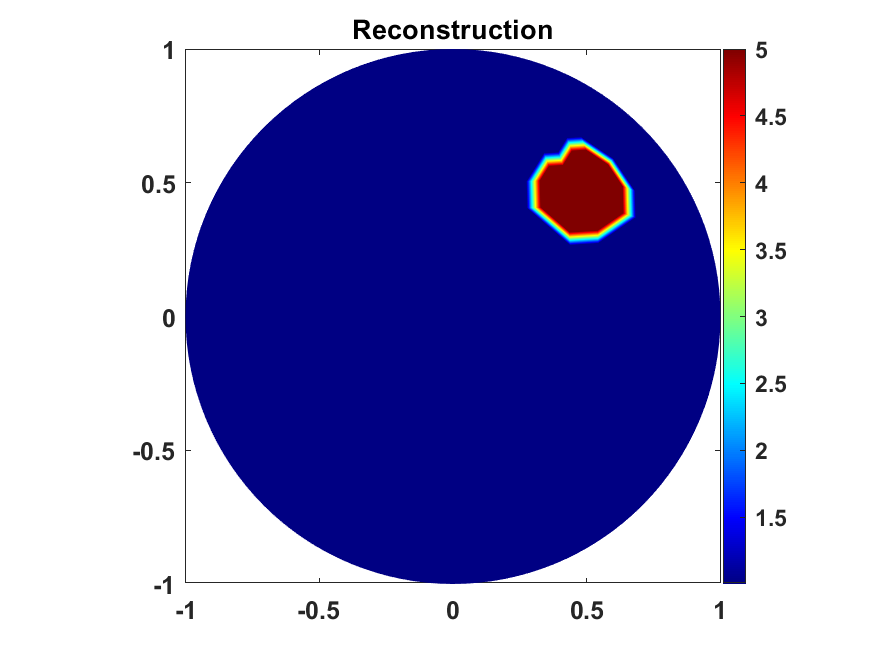}
\subcaption{FEM}
\end{subfigure}
\begin{subfigure}{\textwidth}
\includegraphics[scale=0.32]{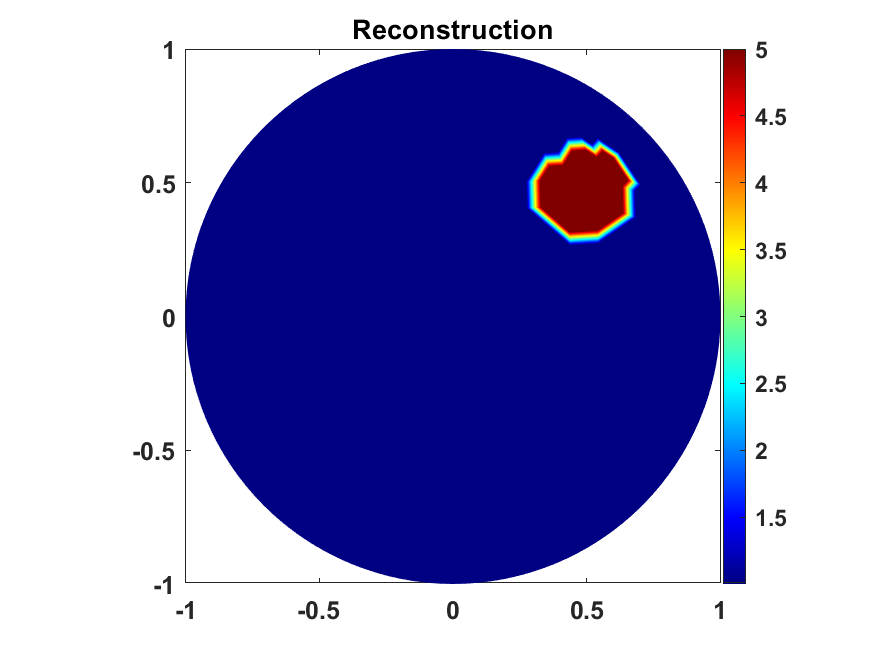}
\subcaption{Net}
\end{subfigure}
\end{minipage}
\hfill
\begin{minipage}{.235\textwidth}
\begin{subfigure}{\textwidth}
\includegraphics[scale=0.32]{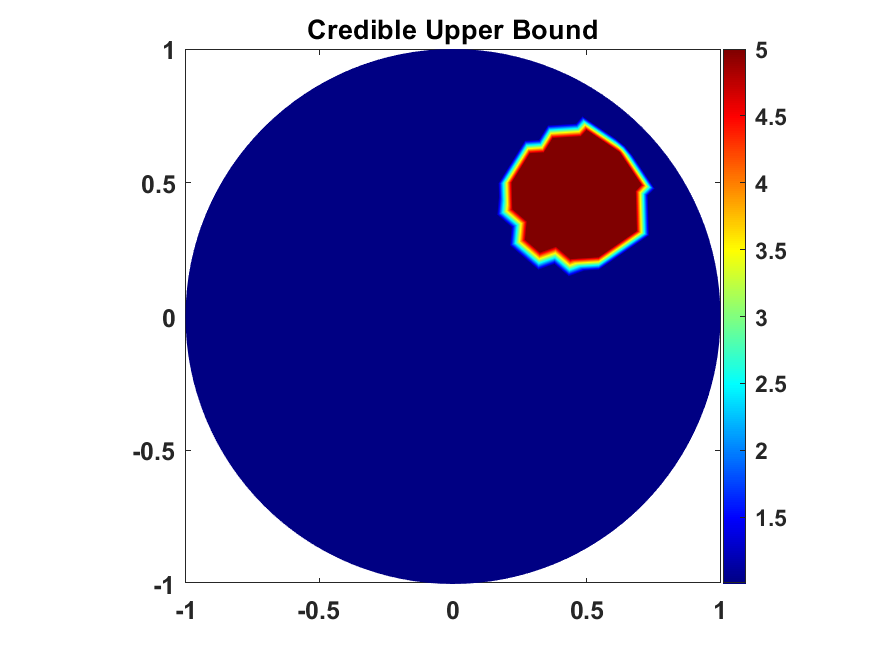}
\subcaption{FEM}
\end{subfigure}
\begin{subfigure}{\textwidth}
\includegraphics[scale=0.32]{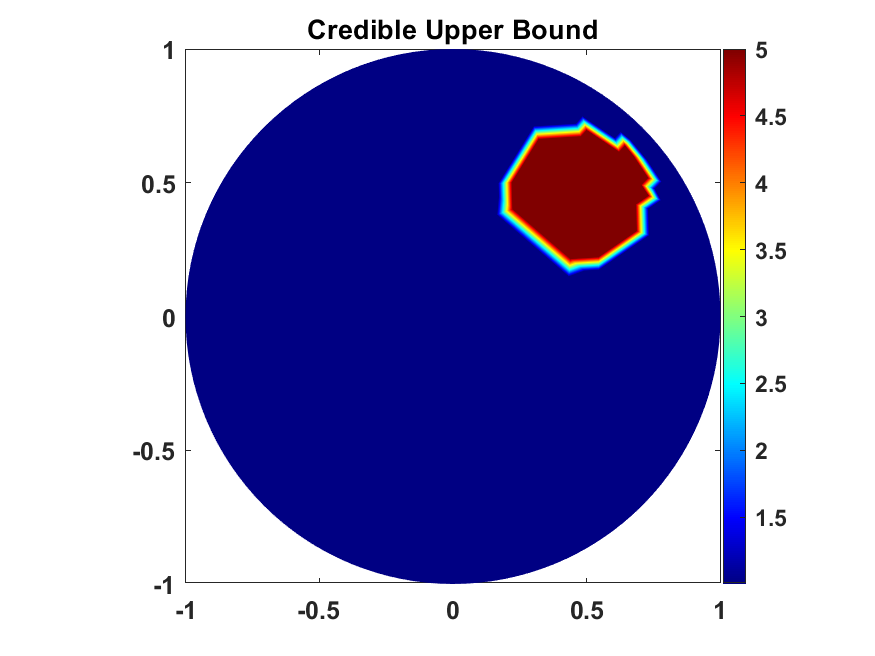}
\subcaption{Net}
\label{eit_sigma_an1_cnn_ucr}
\end{subfigure}
\end{minipage}

\begin{minipage}{.235\textwidth}
\begin{subfigure}{\textwidth}
\includegraphics[scale=0.32]{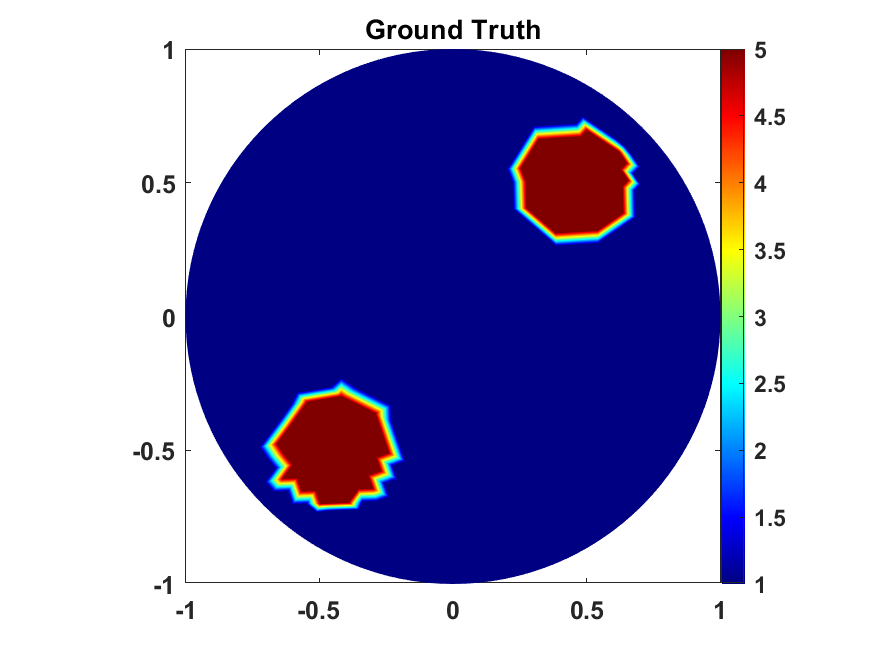}
\subcaption{True}
\label{eit_sigma_an2_true}
\end{subfigure}
\end{minipage}
\hfill
\begin{minipage}{.235\textwidth}
\begin{subfigure}{\textwidth}
\includegraphics[scale=0.32]{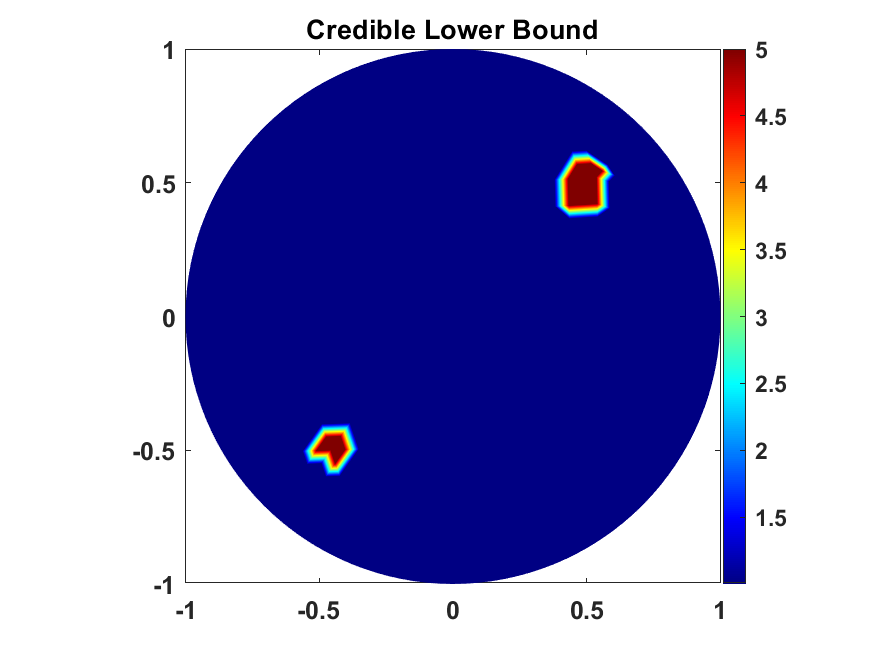}
\subcaption{FEM}
\end{subfigure}
\begin{subfigure}{\textwidth}
\includegraphics[scale=0.32]{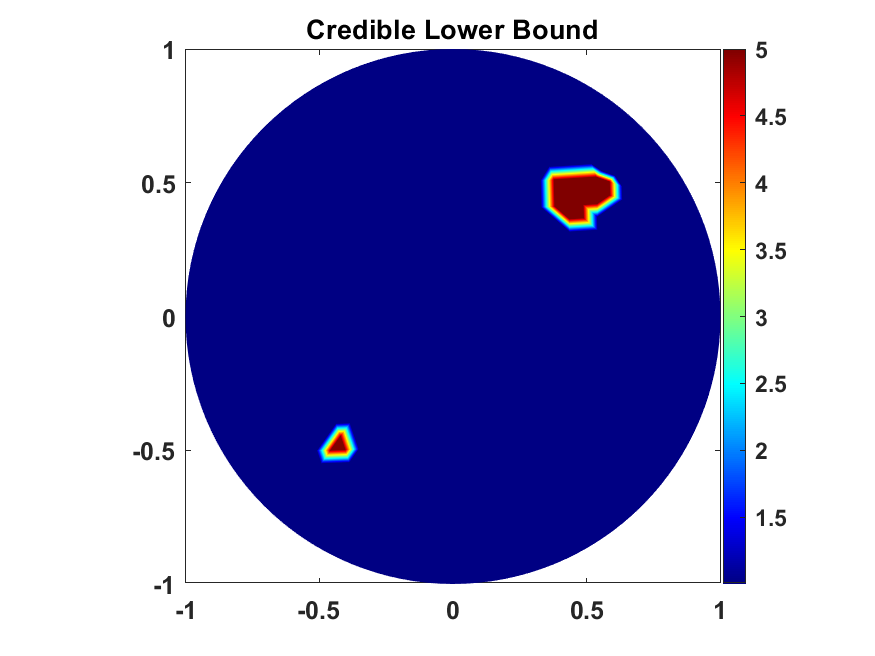}
\subcaption{Net}
\end{subfigure}
\end{minipage}
\hfill
\begin{minipage}{.235\textwidth}
\begin{subfigure}{\textwidth}
\includegraphics[scale=0.32]{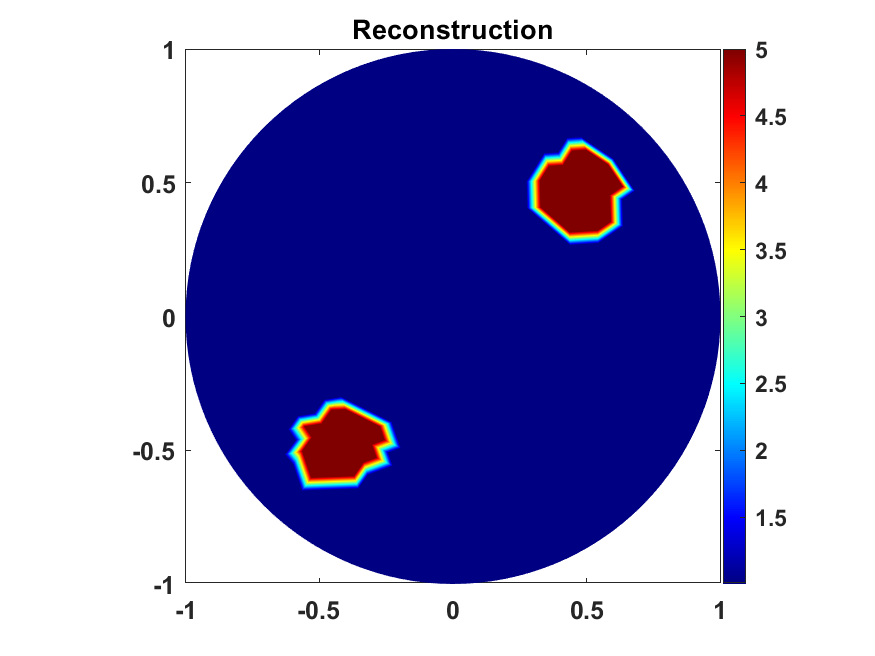}
\subcaption{FEM}
\end{subfigure}
\begin{subfigure}{\textwidth}
\includegraphics[scale=0.32]{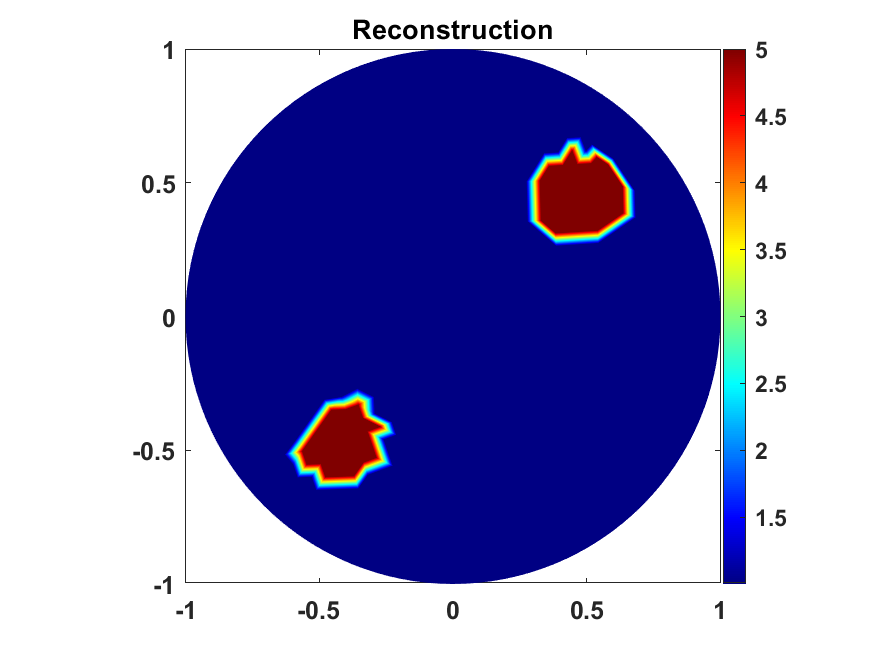}
\subcaption{Net}
\end{subfigure}
\end{minipage}
\hfill
\begin{minipage}{.235\textwidth}
\begin{subfigure}{\textwidth}
\includegraphics[scale=0.32]{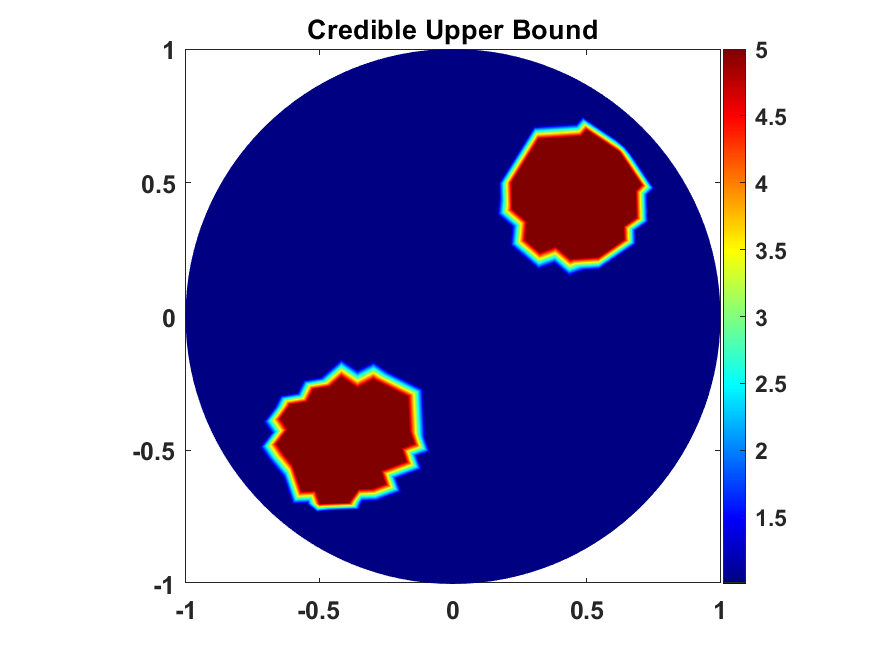}
\subcaption{FEM}
\end{subfigure}
\begin{subfigure}{\textwidth}
\includegraphics[scale=0.32]{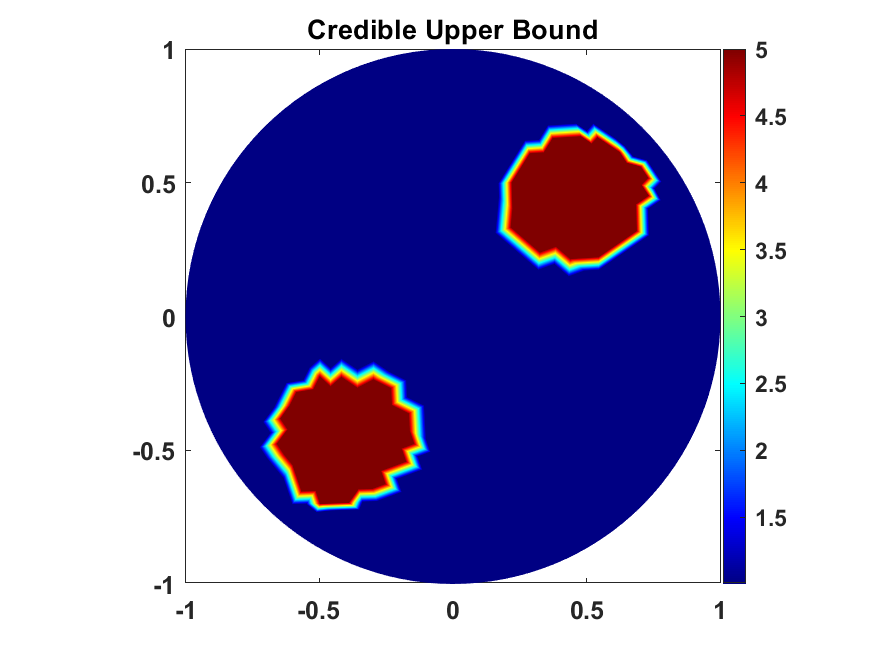}
\subcaption{Net}
\label{eit_sigma_an2_cnn_ucr}
\end{subfigure}
\end{minipage}

\caption{Reconstruction of electrical conductivity using MCMC-FEM and MCMC-Net for EIT (data obtained 1$\%$ relative noise). From left to right: First column: Ground-truth. Second column: Lower Bayesian credible bound (20 $\%$ percentile). Third column: Bayesian estimate in coarse reconstruction mesh. Fourth column: Upper Bayesian credible bound (80 $\%$ percentile). Figure \ref{eit_sigma_an1_true}--\ref{eit_sigma_an1_cnn_ucr} for one anomaly. Figure \ref{eit_sigma_an2_true}--\ref{eit_sigma_an2_cnn_ucr} for two anomalies.}\label{EITan1Emu}
\end{figure}

\begin{figure}
        \centering 
        \begin{subfigure}[b]{0.4\textwidth}
                \includegraphics[scale=0.45]{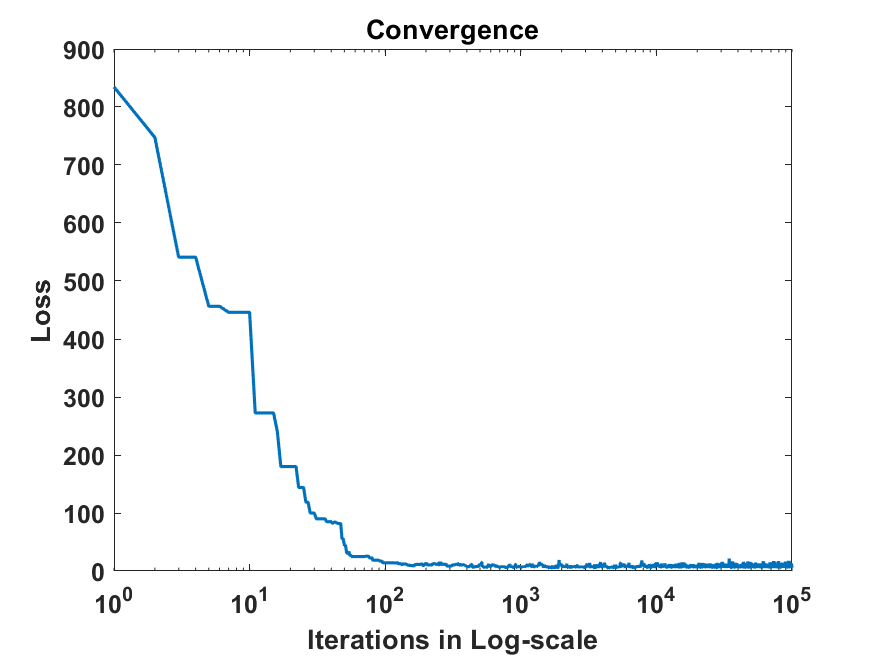}
                \caption{FEM}
                \label{eit_loss_an1_fem}
         \end{subfigure}
        \begin{subfigure}[b]{0.4\textwidth}
                \includegraphics[scale=0.45]{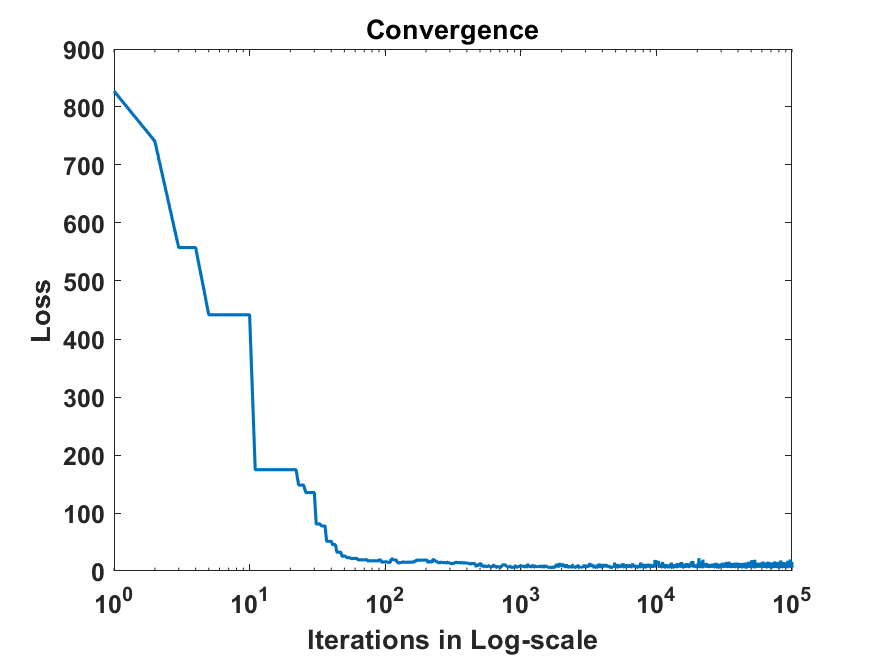}
                \caption{Net}
                \label{eit_loss_an1_cnn}
         \end{subfigure}
\caption{EIT: Convergence of the MCMC method for FEM and CNN. Data obtained at 1 $\%$ relative noise.}\label{EITan1Loss}
\end{figure}

\begin{figure}[!h]
    \centering 

    \begin{subfigure}[b]{0.32\textwidth}
      \includegraphics[height=0.15\textheight]{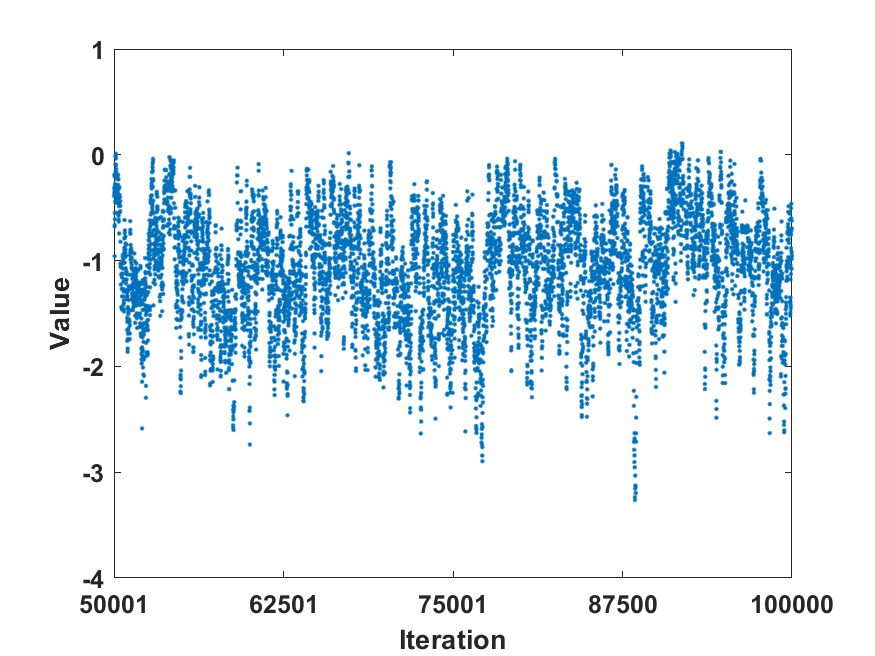}
        \label{TracePlot_an1_600_TS100KBRN50K}
    \end{subfigure}
    \begin{subfigure}[b]{0.32\textwidth}
    \includegraphics[height=0.15\textheight]{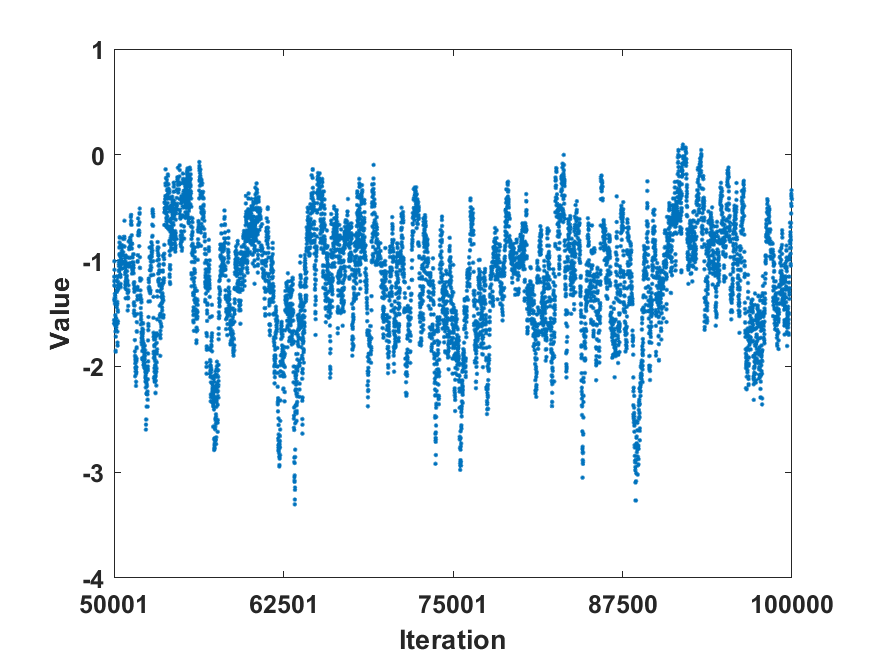}
        \label{TracePlot_an2_600_TS100KBRN50K}
    \end{subfigure}
    \begin{subfigure}[b]{0.32\textwidth}
    \includegraphics[height=0.15\textheight]{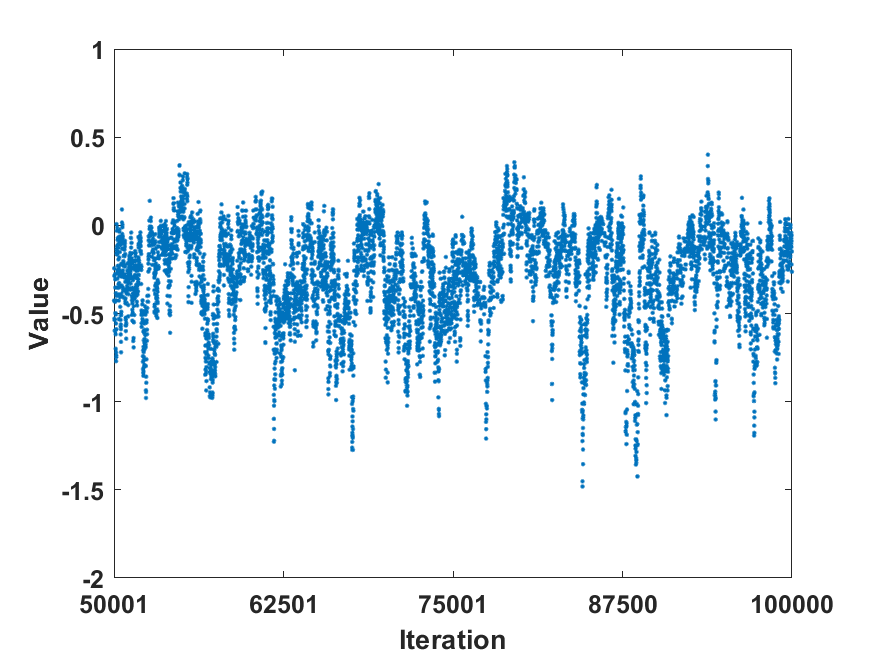}
    \label{TracePlot_an3_600_TS100KBRN50K}
    \end{subfigure}

    \begin{subfigure}[b]{0.32\textwidth}
      \includegraphics[height=0.15\textheight]{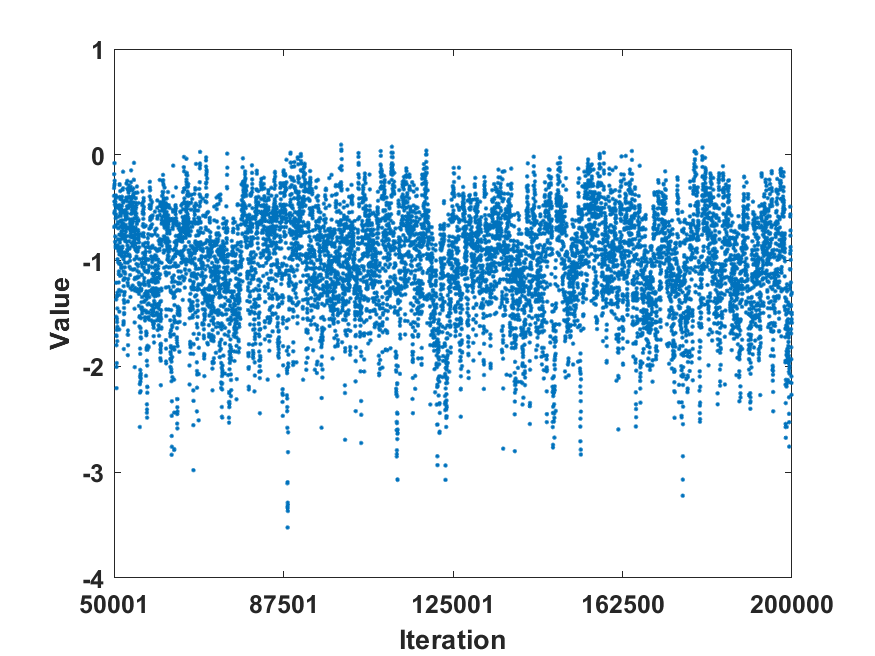}
        \label{TracePlot_an1_600_TS200KBRN50K}
    \end{subfigure}
    \begin{subfigure}[b]{0.32\textwidth}
    \includegraphics[height=0.15\textheight]{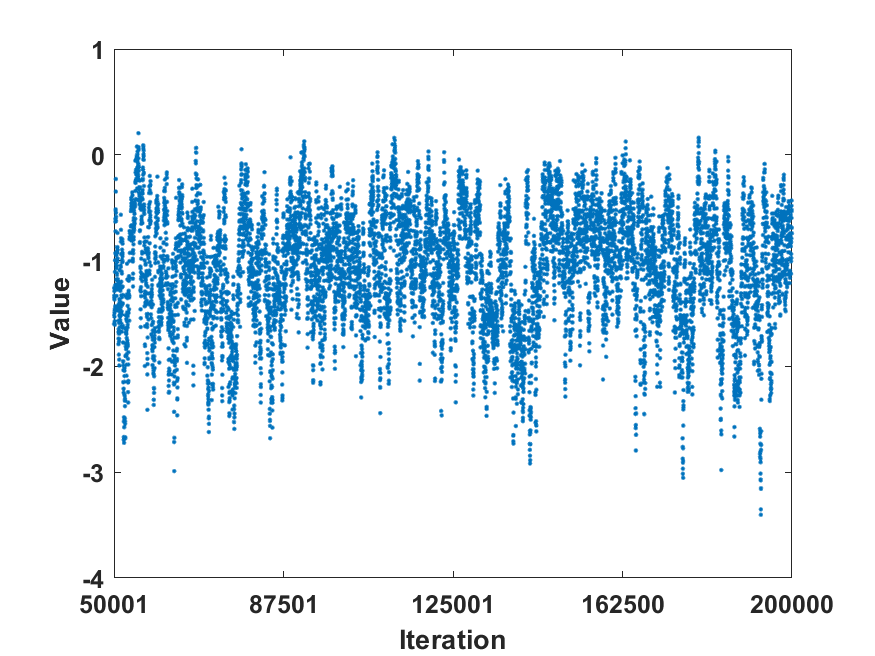}
        \label{TracePlot_an2_600_TS200KBRN50K}
    \end{subfigure}
    \begin{subfigure}[b]{0.32\textwidth}
    \includegraphics[height=0.15\textheight]{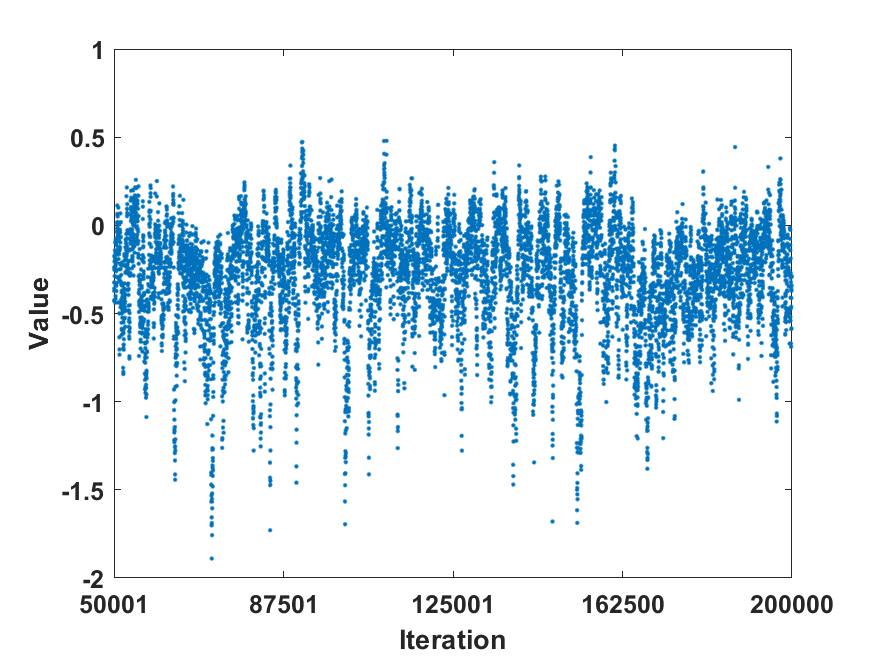}
    \label{TracePlot_an3_600_TS200KBRN50K}
    \end{subfigure}

\caption{
Trace plots of the posterior samples of level-set functions using the proposed approach for EIT inverse problem at a randomly chosen grid point. The top row corresponds to plotted values at this grid point with $100{,}000$ total samples and a burn-in of $50{,}000$, while the bottom row uses $200{,}000$ total samples with the same burn-in. The three columns represent different anomaly configurations: one anomaly (left), two anomalies, and three anomalies (right). These plots illustrate the convergence behavior and mixing of the MCMC chains at the selected grid location.
}
\label{ConvergencediagnosticsEIT}
\end{figure}

In Figure \ref{EITan1Emu}, we compare the electrical conductivity reconstruction results for the EIT inverse problem using MCMC-FEM and MCMC-Net, considering one and two anomalies. To visualize the convergence of the MCMC method for both approaches, see Figure \ref{EITan1Loss}. Figure \ref{ConvergencediagnosticsEIT} shows the corresponding mixing plots for EIT. It is observed that the MCMC chain shows good mixing after discarding the first 50,000 samples as burn-in.Table \ref{Reconstruction_errors_eit} presents the different errors calculated between the ground truth and the reconstructions for both MCMC-FEM and MCMC-Net. The last column of Table \ref{Reconstruction_errors_eit} describes the corresponding computation times. In addition to the Bayesian estimates, the associated credible regions are also shown in Figure~\ref{EITan1Emu}.

In the one-dimensional case, a Bayesian credible interval of size \(1-\alpha\) is an interval \((a, b)\) such that \(P(a \leq \tau \leq b | \{X_i\}_{i=1}^n) = 1-\alpha\), where \(\tau\) is a random variable, and \(\{X_i\}_{i=1}^n\) are the given samples of \(\tau\). In this paper, a \(1-\alpha\) credible region refers to a two-dimensional extension of the Bayesian credible interval. Specifically, a region of \(1-\alpha\) credibility is defined over a set of non-overlapping triangles in the reconstruction mesh \(\Omega_{\text{discrete}}\) covering \(\Omega\). For each triangle \(t \in \Omega_{\text{discrete}}\), we have an interval \((a(t), b(t))\) such that \(P(a(t) \leq \tau(t) \leq b(t) | \{X_i(t)\}_{i=1}^n) = 1-\alpha\), where \(\tau(t)\) is the random variable within the corresponding triangle. 

By analyzing all the results described above, we can conclude that by replacing the FEM forward solver with a CNN forward solver, we can reconstruct the internal conductivity using MCMC-Net with the same accuracy as the one using MCMC-FEM but at about thirty fold lower computing cost.


\begin{table}
\caption{EIT: Comparison of reconstruction errors and inversion time (Inv. Time) in minutes between MCMC-FEM (FEM) and MCMC-Net (Net).}
\label{Reconstruction_errors_eit}
\centering
\renewcommand{\arraystretch}{1.6} 
\scalebox{0.95}{
\begin{tabular}[t]{|c|cc|cc|cc|}
\hline
\multirow{2}{*}{\shortstack{Anomaly type \\ (Circular)}}   & \multicolumn{2}{c|}{MAE}  & \multicolumn{2}{c|}{MSE} &
\multicolumn{2}{c|}{Inv. Time (m)}  \\
\cline{2-7}
& \multicolumn{1}{c}{FEM} & \multicolumn{1}{c|}{Net} & 
\multicolumn{1}{c}{FEM} & \multicolumn{1}{c|}{Net} & \multicolumn{1}{c}{FEM} & \multicolumn{1}{c|}{Net}   \\
\hline
One      & 0.071397 & 0.068423  & 0.285499 & 0.273603  & 124.85  & 4.09 \\
Two      & 0.154669 & 0.151695  & 0.618581 & 0.606685  & 124.86  & 4.12 \\
Three    & 0.282551 & 0.350951  & 1.130100 & 1.403706  & 125.93  & 4.11 \\
\hline    
\end{tabular}
}
\end{table}

\subsection{Diffuse Optical Tomography}
\begin{figure}
\centering
\begin{minipage}{.235\textwidth}
\begin{subfigure}{\textwidth}
\includegraphics[scale=0.32]{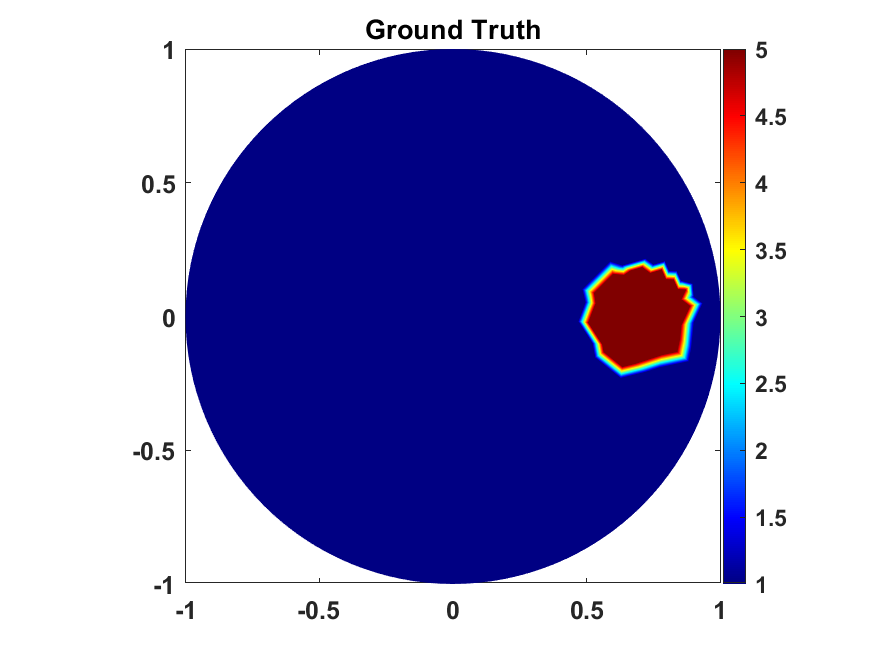}
\subcaption{True}
\label{dot_mu_an1_true}
\end{subfigure}
\end{minipage}
\hfill
\begin{minipage}{.235\textwidth}
\begin{subfigure}{\textwidth}
\includegraphics[scale=0.32]{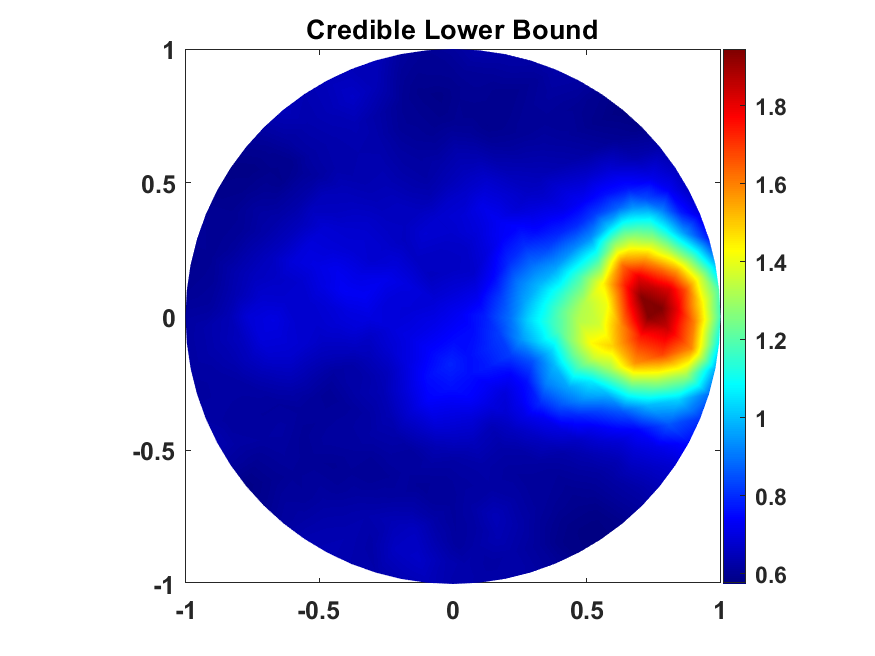}
\subcaption{FEM}
\end{subfigure}
\begin{subfigure}{\textwidth}
\includegraphics[scale=0.32]{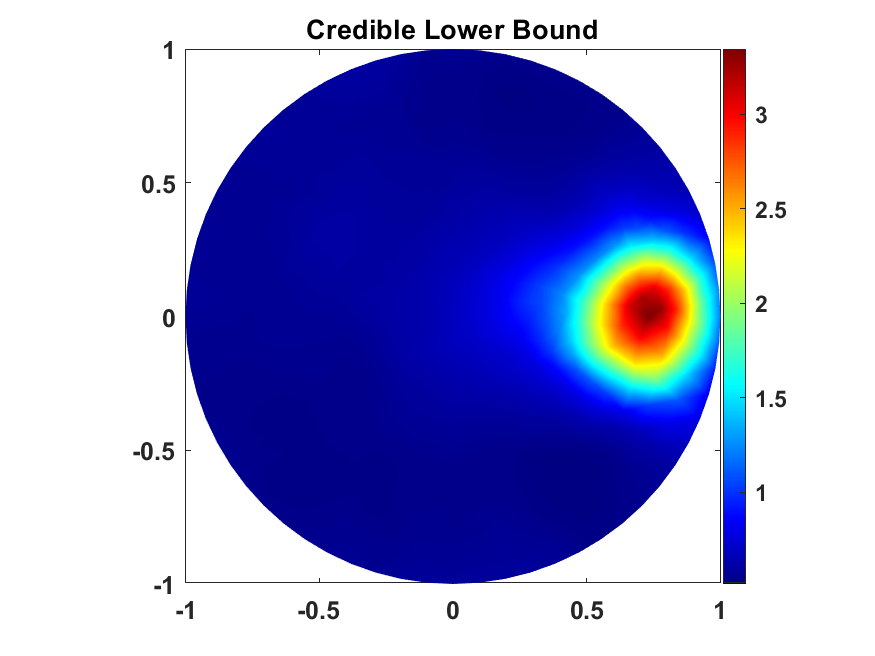}
\subcaption{Net}
\end{subfigure}
\end{minipage}
\hfill
\begin{minipage}{.235\textwidth}
\begin{subfigure}{\textwidth}
\includegraphics[scale=0.32]{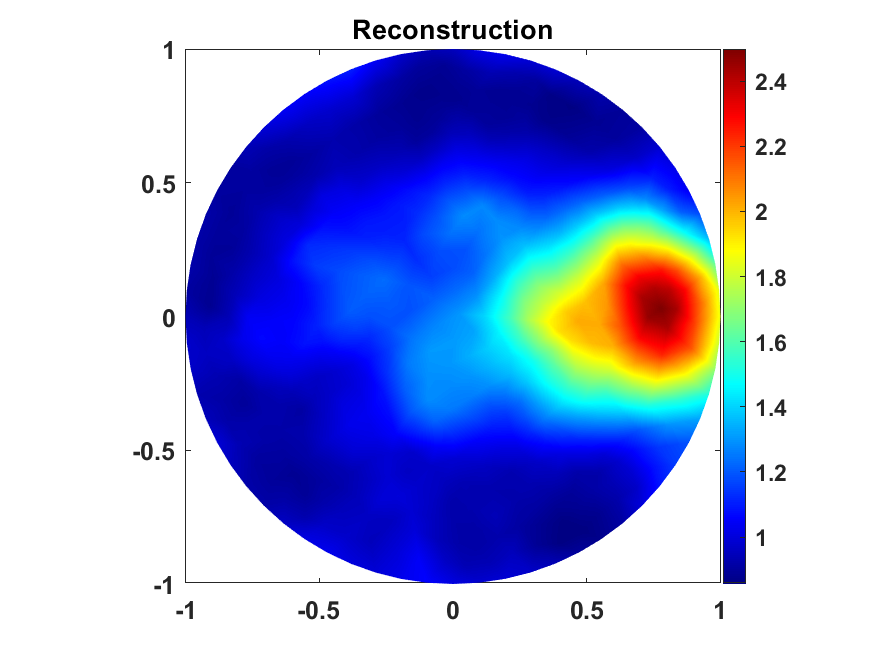}
\subcaption{FEM}
\end{subfigure}
\begin{subfigure}{\textwidth}
\includegraphics[scale=0.32]{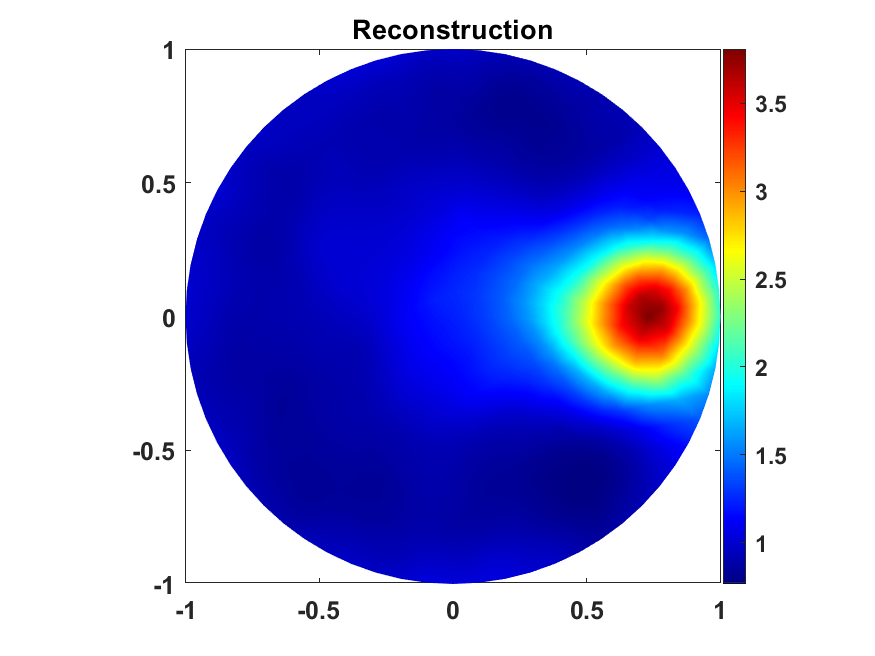}
\subcaption{Net}
\end{subfigure}
\end{minipage}
\hfill
\begin{minipage}{.235\textwidth}
\begin{subfigure}{\textwidth}
\includegraphics[scale=0.32]{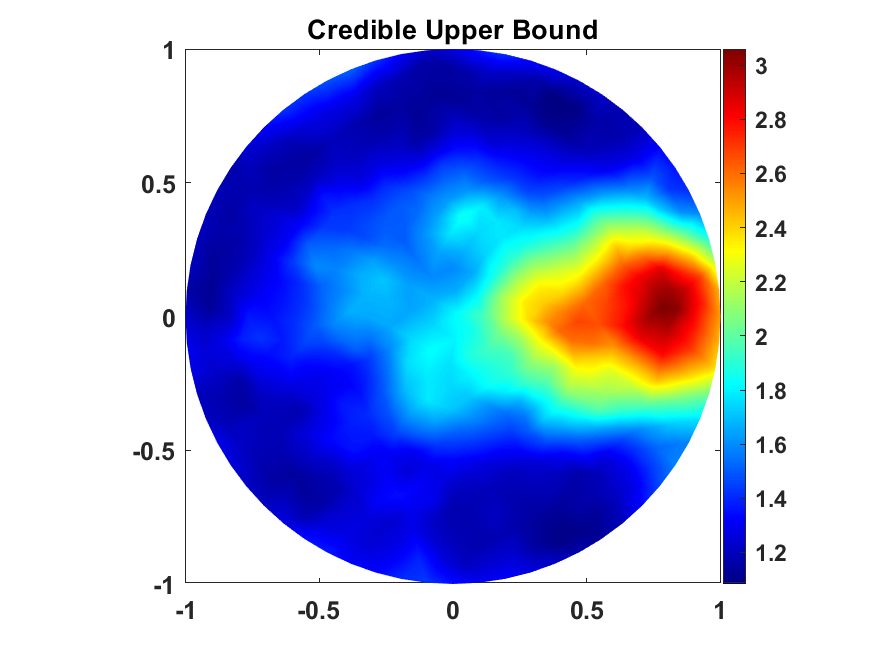}
\subcaption{FEM}
\end{subfigure}
\begin{subfigure}{\textwidth}
\includegraphics[scale=0.32]{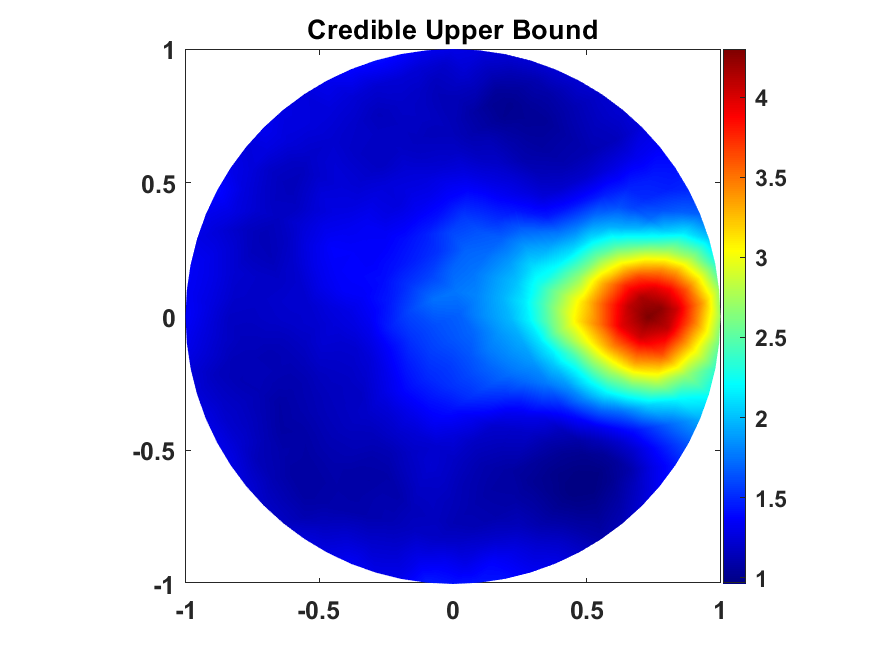}
\subcaption{Net}
\label{dot_mu_an1_cnn_ucr}
\end{subfigure}
\end{minipage}

\begin{minipage}{.235\textwidth}
\begin{subfigure}{\textwidth}
\includegraphics[scale=0.32]{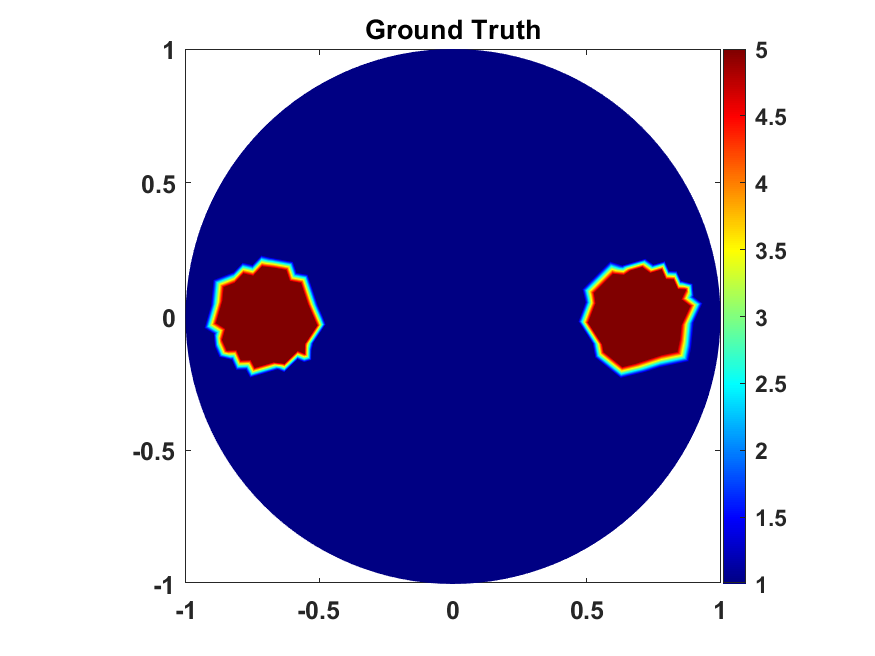}
\subcaption{True}
\label{dot_mu_an2_true}
\end{subfigure}
\end{minipage}
\hfill
\begin{minipage}{.235\textwidth}
\begin{subfigure}{\textwidth}
\includegraphics[scale=0.32]{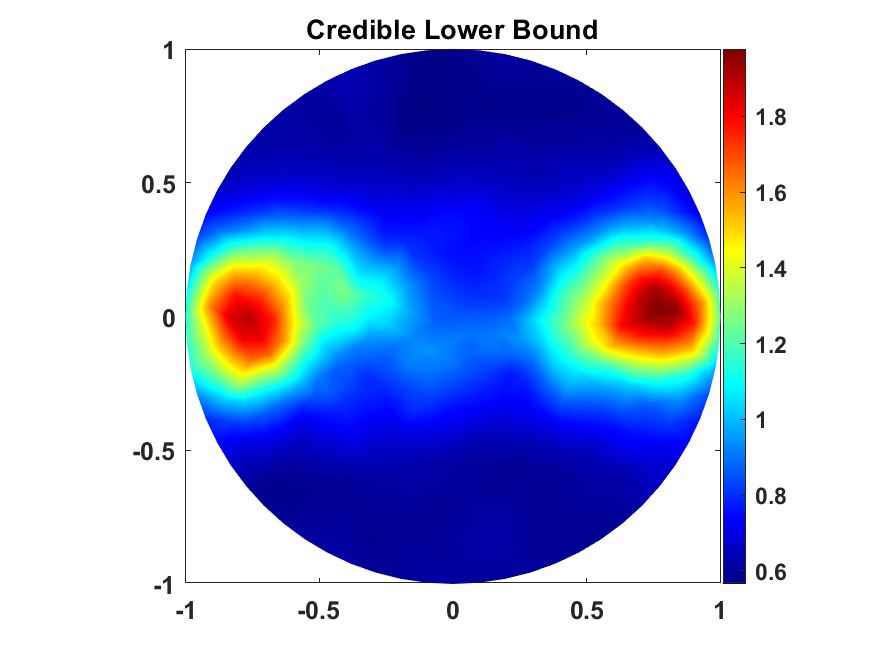}
\subcaption{FEM}
\end{subfigure}
\begin{subfigure}{\textwidth}
\includegraphics[scale=0.32]{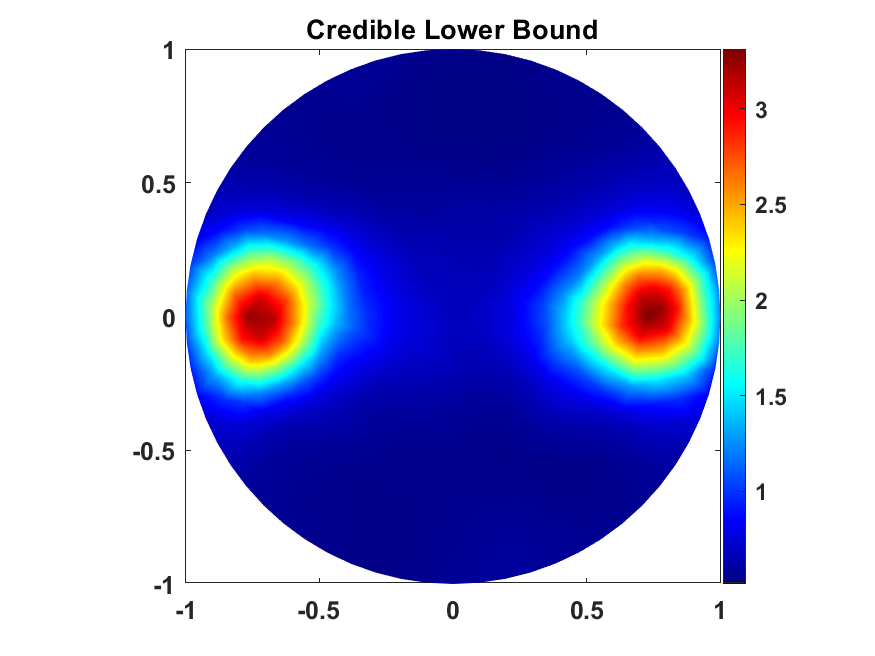}
\subcaption{Net}
\end{subfigure}
\end{minipage}
\hfill
\begin{minipage}{.235\textwidth}
\begin{subfigure}{\textwidth}
\includegraphics[scale=0.32]{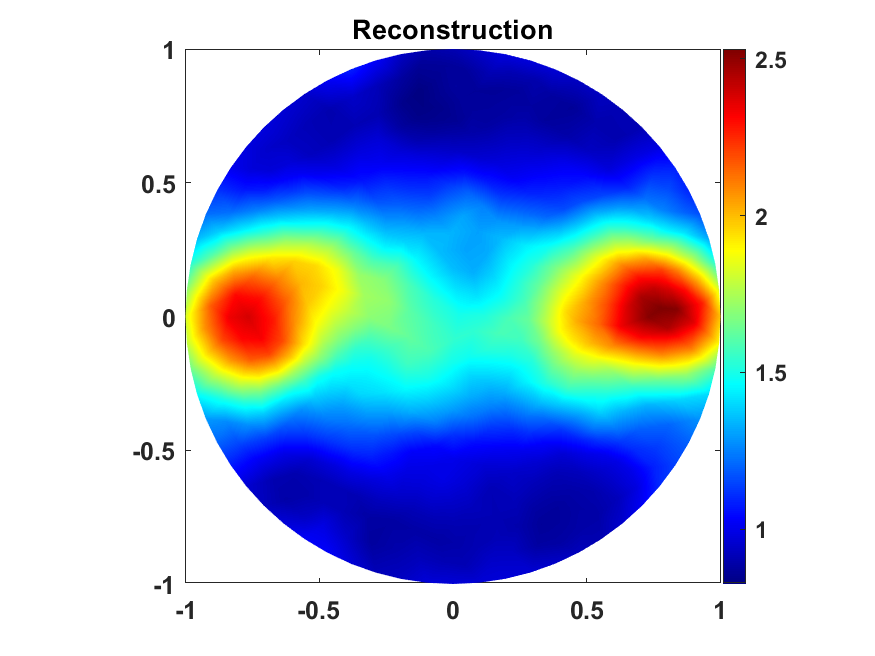}
\subcaption{FEM}
\end{subfigure}
\begin{subfigure}{\textwidth}
\includegraphics[scale=0.32]{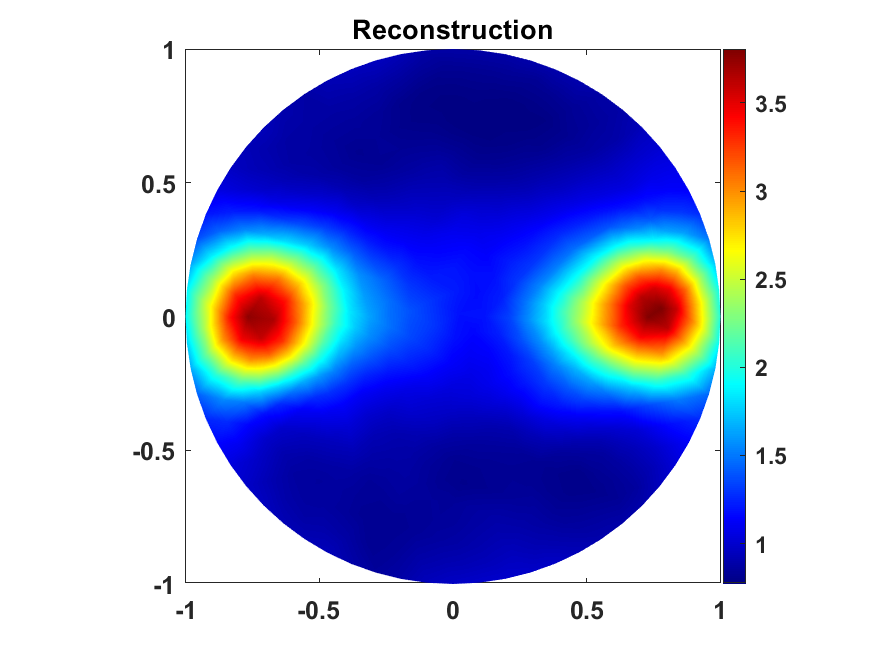}
\subcaption{Net}
\end{subfigure}
\end{minipage}
\hfill
\begin{minipage}{.235\textwidth}
\begin{subfigure}{\textwidth}
\includegraphics[scale=0.32]{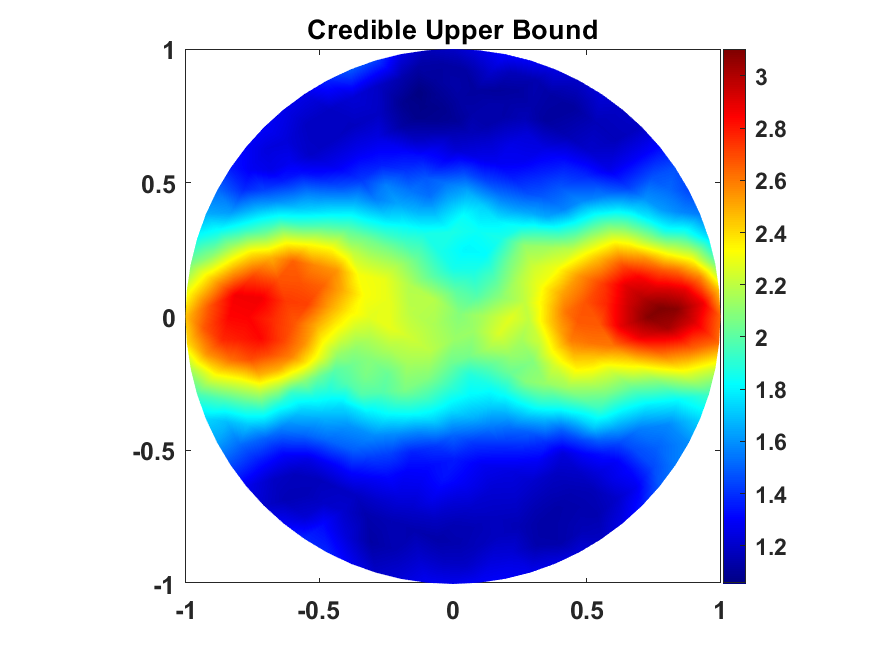}
\subcaption{FEM}
\end{subfigure}
\begin{subfigure}{\textwidth}
\includegraphics[scale=0.32]{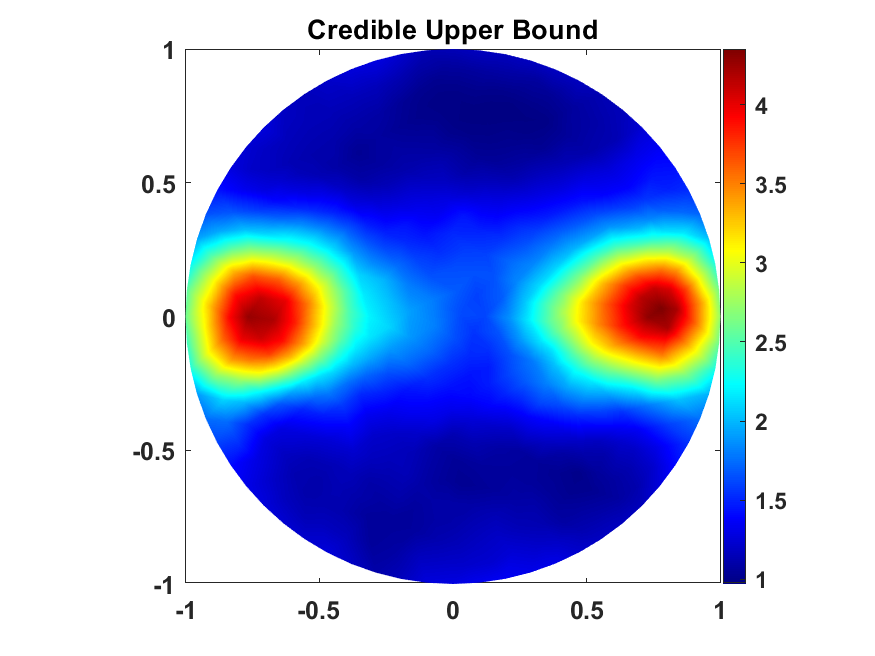}
\subcaption{Net}
\label{dot_mu_an2_cnn_ucr}
\end{subfigure}
\end{minipage}

\caption{Reconstruction of absorption coefficient ($\mu$) using MCMC-FEM and MCMC-Net for DOT (data obtained 1$\%$ relative noise). From left to right: First column: Ground-truth. Second column: Lower Bayesian credible bound (20 $\%$ percentile). Third column: Bayesian estimate in coarse reconstruction mesh. Fourth column: Upper Bayesian credible bound (80 $\%$ percentile). Figure \ref{dot_mu_an1_true}--\ref{dot_mu_an1_cnn_ucr} for one anomaly. Figure \ref{dot_mu_an2_true}--\ref{dot_mu_an2_cnn_ucr} for two anomalies.}\label{DOTan1Emu}
\end{figure}

In this section, we discuss the numerical simulation results of the inverse problem of DOT. This problem is solved using the MCMC-FEM approach and discussed in \cite{abhishek2022optimal}. In this study, we replace the FEM forward solver with CNN and compute the reconstruction results for the absorption coefficient and represent the results in Figure \ref{DOTan1Emu}. Similar to the previous case of EIT results, we can see that by replacing the FEM forward solver with CNN forward solver, we can reconstruct the absorption coefficient similarly well. 

Different errors computed between the ground truth and the reconstructed results are presented in Table \ref{Reconstruction_errors_dot}, further supporting our findings. The computational time for both approaches is shown in Column 4 of Table \ref{Reconstruction_errors_dot}. The proposed approach takes approximately nine times less computational time to achieve similar accuracy as the FEM forward solver.

\begin{table}
\caption{DOT: Comparison of reconstruction errors and inversion time in minutes between MCMC-FEM and MCMC-Net.}
\label{Reconstruction_errors_dot}
\centering
\renewcommand{\arraystretch}{1.5} 

\scalebox{0.9}{
\begin{tabular}[t]{|c|cc|cc|cc|cc|}
\hline
\multirow{2}{*}{\shortstack{Anomaly type \\ (Circular)}} & \multicolumn{2}{c|}{$L^{\infty}$ -loss}   & \multicolumn{2}{c|}{MAE}  & \multicolumn{2}{c|}{MSE}  &
\multicolumn{2}{c|}{Inv. Time (m)}  \\
\cline{2-9}		
& \multicolumn{1}{c}{FEM} & \multicolumn{1}{c|}{Net}  & \multicolumn{1}{c}{FEM} & \multicolumn{1}{c|}{Net} &
\multicolumn{1}{c}{FEM} & \multicolumn{1}{c|}{Net} & \multicolumn{1}{c}{FEM} & \multicolumn{1}{c|}{Net}   \\
\hline
One    & 2.9909 & 2.3438  & 0.2818 & 0.2589  & 0.3632 & 0.2416 &  34.25 & 3.60 \\
Two    & 3.1235 & 2.5087  & 0.4863 & 0.4322  & 0.7260 & 0.4850 &  34.29 & 3.62\\
\hline      
\end{tabular}
}
\end{table}
\subsection{Quantitative Photoacoustic Tomography}
\begin{figure}
\centering

\begin{minipage}{.235\textwidth}
\begin{subfigure}{\textwidth}
\includegraphics[scale=0.32]{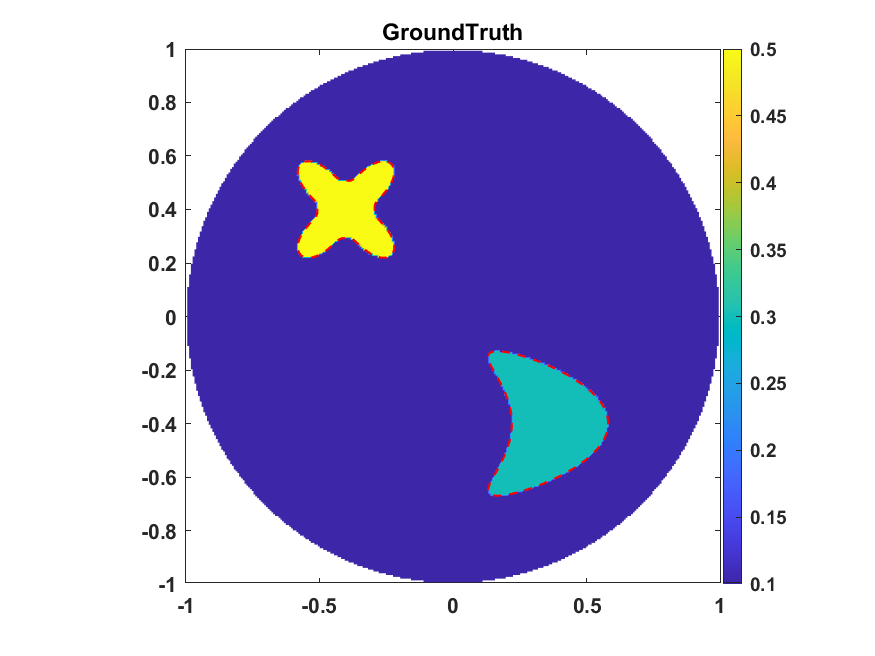}
\subcaption{True}
\label{qpat_gamma_true}
\end{subfigure}
\end{minipage}
\hfill
\begin{minipage}{.235\textwidth}
\begin{subfigure}{\textwidth}
\includegraphics[scale=0.32]{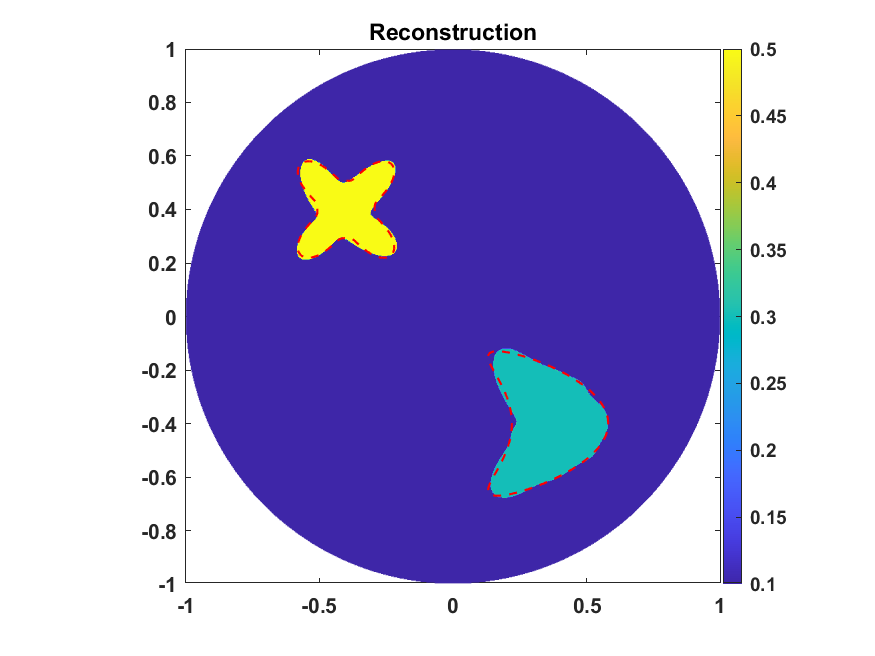}
\subcaption{FEM}
\end{subfigure}
\begin{subfigure}{\textwidth}
\includegraphics[scale=0.32]{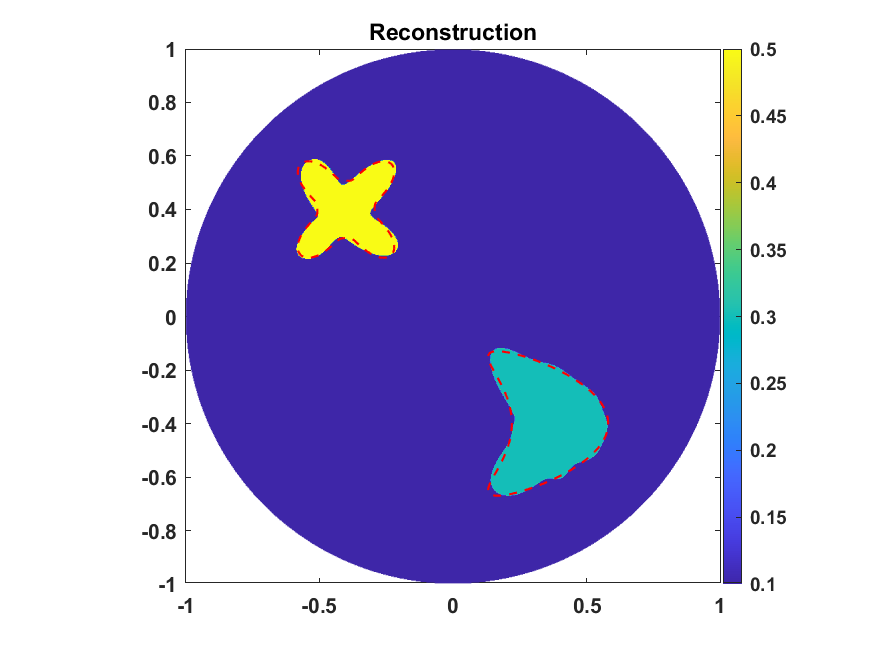}
\subcaption{Net}
\end{subfigure}
\end{minipage}
\hfill
\begin{minipage}{.235\textwidth}
\begin{subfigure}{\textwidth}
\includegraphics[scale=0.32]{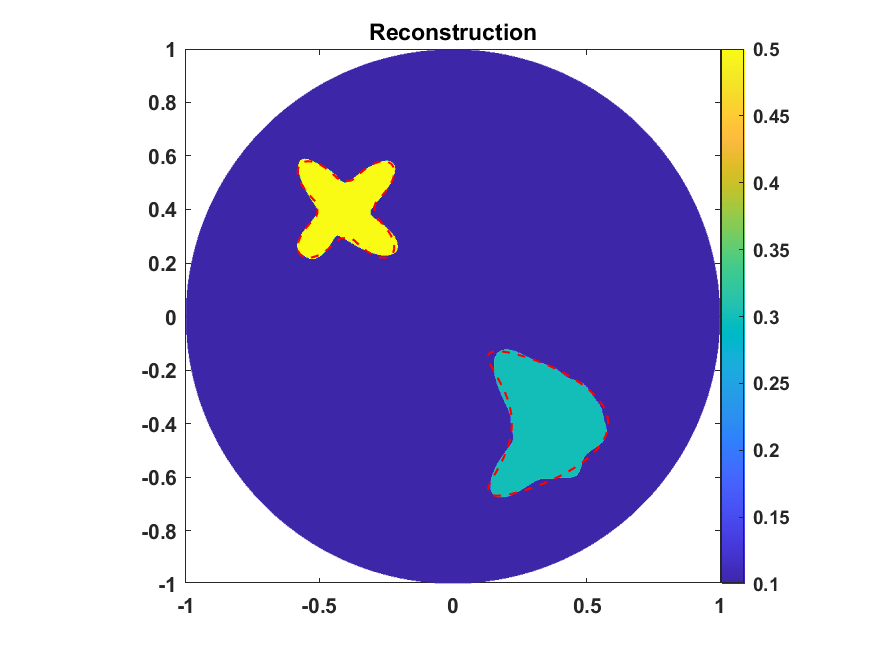}
\subcaption{FEM}
\end{subfigure}
\begin{subfigure}{\textwidth}
\includegraphics[scale=0.32]{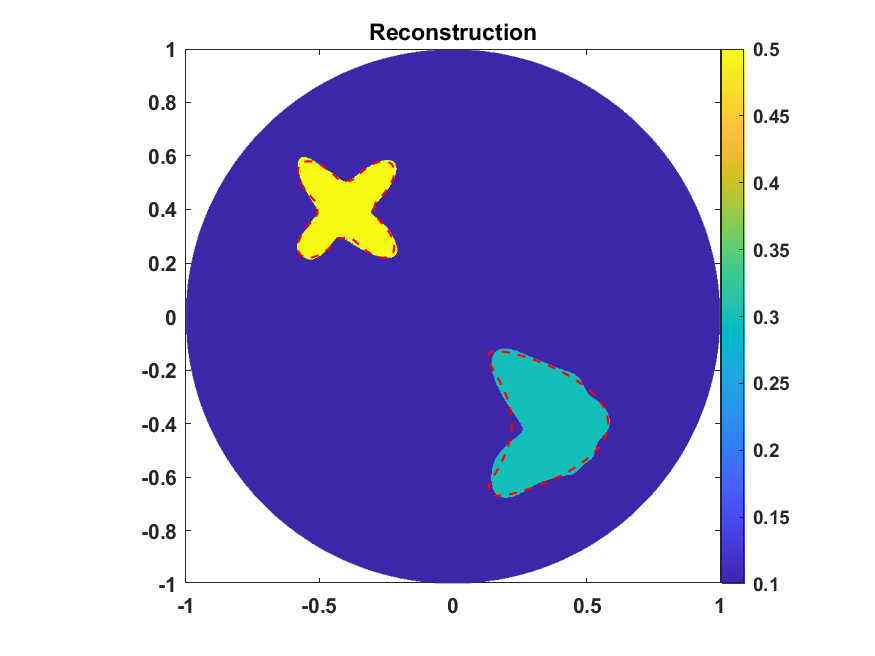}
\subcaption{Net}
\end{subfigure}
\end{minipage}
\hfill
\begin{minipage}{.235\textwidth}
\begin{subfigure}{\textwidth}
\includegraphics[scale=0.32]{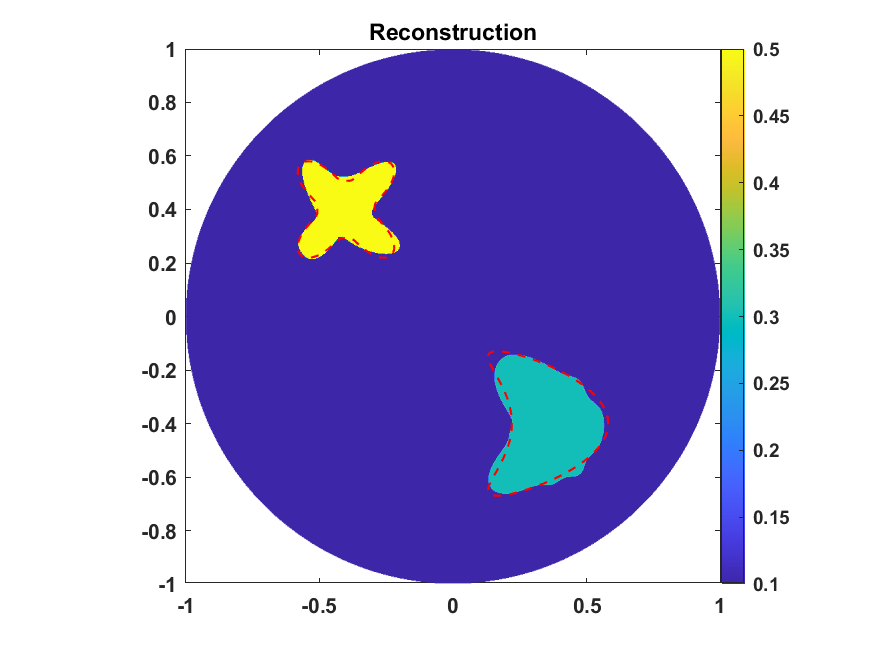}
\subcaption{FEM}
\end{subfigure}
\begin{subfigure}{\textwidth}
\includegraphics[scale=0.32]{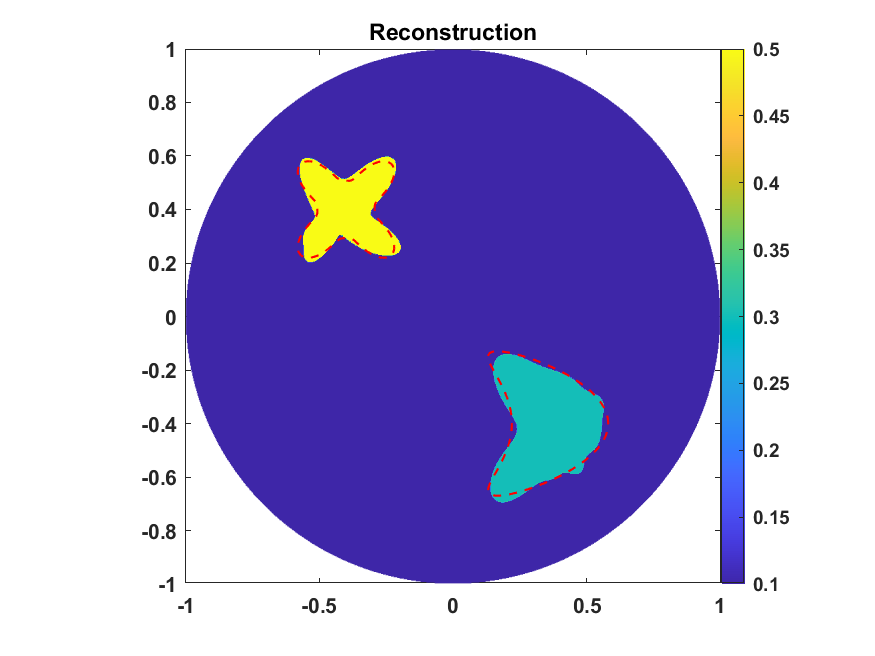}
\subcaption{Net}
\end{subfigure}
\end{minipage}

    \caption{Reconstruction of the absorption coefficient using MCMC-FEM and MCMC-Net for QPAT. From left to right: First column: Ground truth. Second column: Reconstruction (data obtained 2 $\%$ relative noise). Third column: Reconstruction (data obtained 4 $\%$ relative noise). Fourth column: Reconstruction (data obtained 8 $\%$ relative noise).}
    \label{qPATNoise2}
\end{figure}

In the final part of the numerical results, we describe the absorption coefficient reconstructions for the inverse problem of QPAT \cite{afkham2024bayesian}. For this problem, we reproduced the results for FEM forward solver using the code as in \cite{afkham2024bayesian}. Then we replaced the FEM forward solver with the proposed CNN solver. We display the reconstruction results in Figure \ref{qPATNoise2} for various noise levels. The computed errors and computation times are summarized in Table \ref{Reconstruction_errors_qpat}. After analyzing all the results, we infer that replacing the FEM forward solver with the CNN solver enables us to reconstruct an absorption coefficient at a similar quality to that obtained with FEM but at a lowered computational cost.

\subsection{Network Sensitivity Analysis and Robustness Assessment} 
\begin{figure}[!h]
       \centering 
         \begin{subfigure}[b]{0.24\textwidth}
          \includegraphics[scale=0.31]{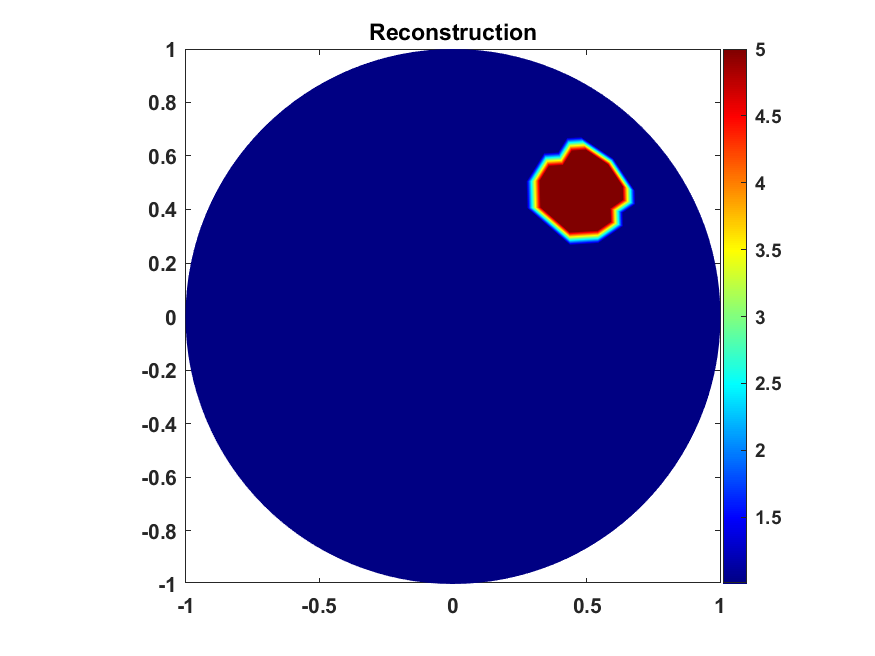}
                \caption{NL=4}
                \label{EIT_SIGAN1_L4_N16_MB64_EP2K}
         \end{subfigure}
         \begin{subfigure}[b]{0.24\textwidth}
        \includegraphics[scale=0.31]{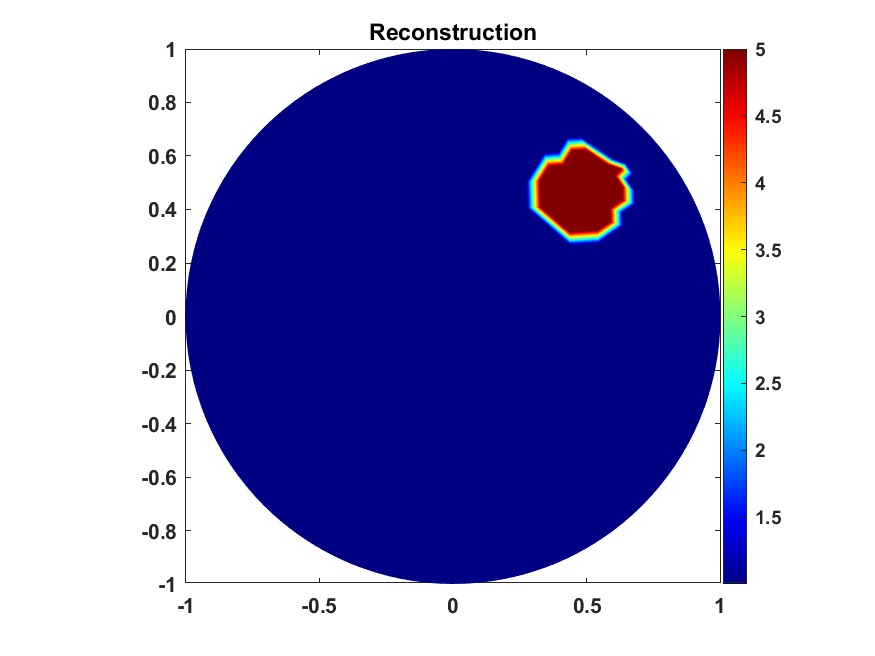}
                \caption{NL=6}
                \label{EIT_SIGAN1_L6_N16_MB64_EP2K}
         \end{subfigure}
          \begin{subfigure}[b]{0.24\textwidth}
         \includegraphics[scale=0.31]{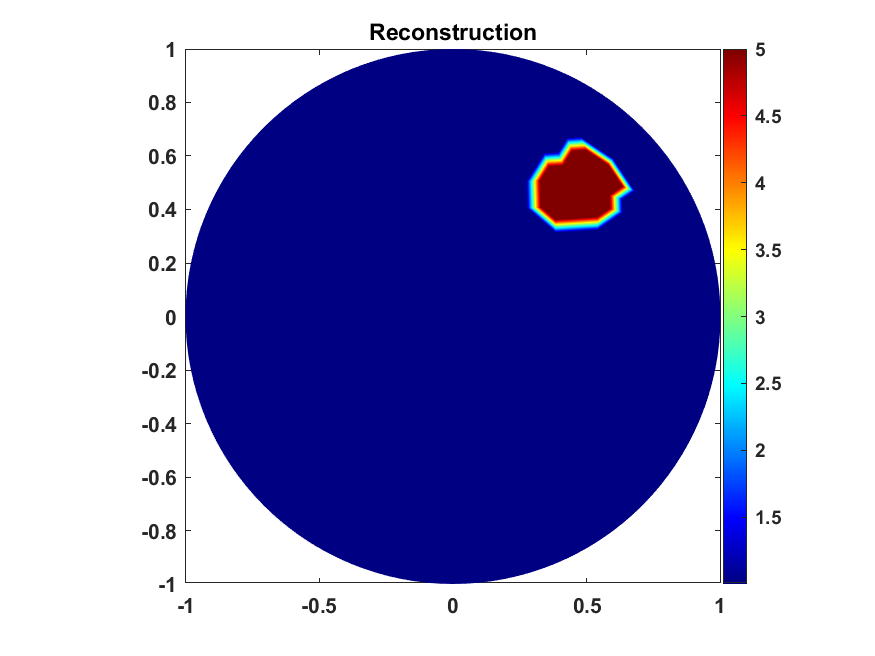}
               \caption{NL=8}
                \label{EIT_SIGAN1_L8_N16_MB64_EP2K}
         \end{subfigure}
        \begin{subfigure}[b]{0.24\textwidth}
          \includegraphics[scale=0.31]{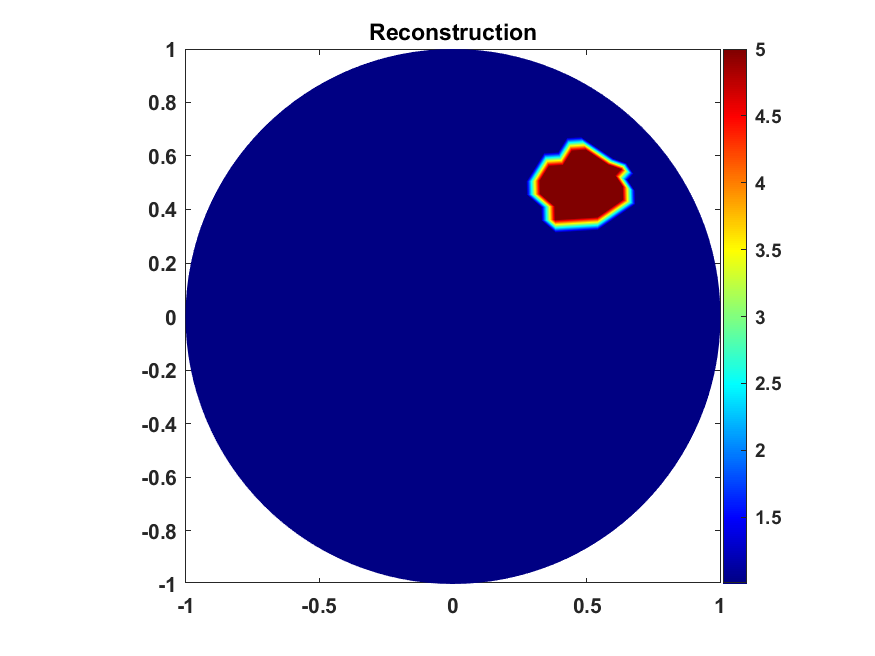}
                \caption{NL=4, EP=1000}
                \label{EIT_SIGAN1_L4_N16_MB128_EP1K}
         \end{subfigure}
         
            \begin{subfigure}[b]{0.24\textwidth}
          \includegraphics[scale=0.31]{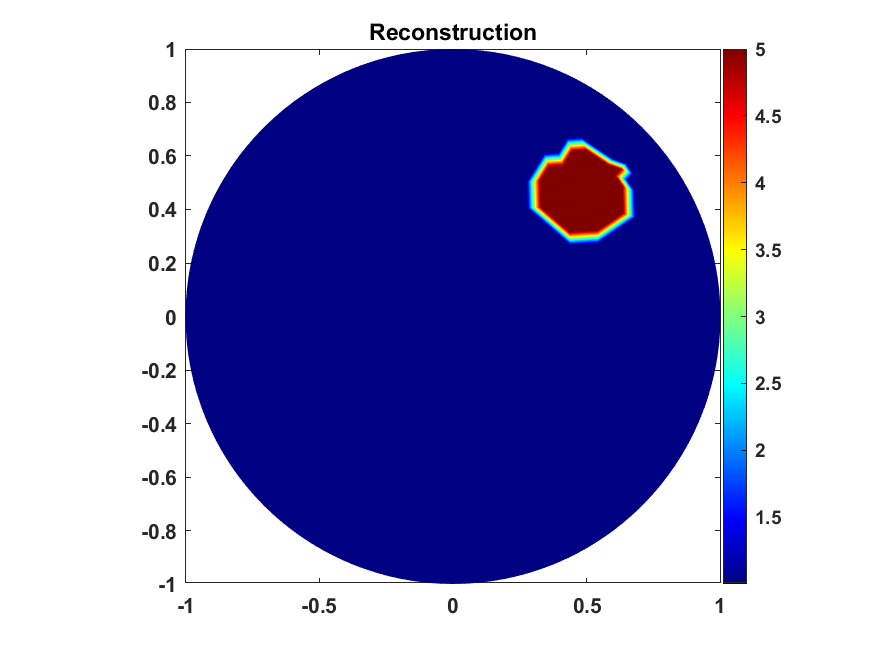}
                \caption{NL=4, EP=2000}
                \label{EIT_SIGAN1_L4_N16_MB128_EP2K}
         \end{subfigure}
            \begin{subfigure}[b]{0.24\textwidth}
          \includegraphics[scale=0.31]{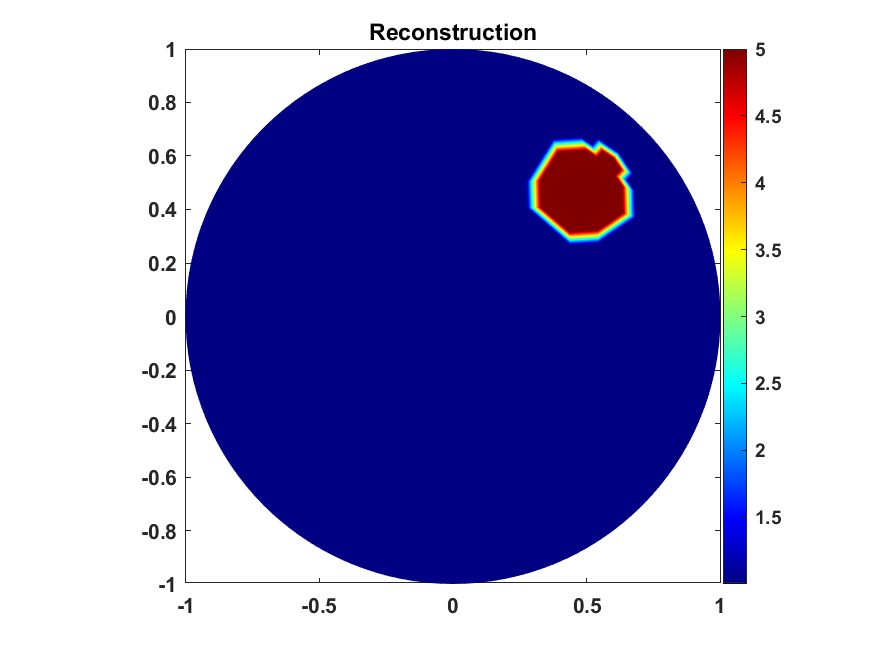}
                \caption{NL=4, EP=5000}
                \label{EIT_SIGAN1_L4_N16_MB128_EP5K}
         \end{subfigure}
         \begin{subfigure}[b]{0.24\textwidth}
                \includegraphics[scale=0.31]{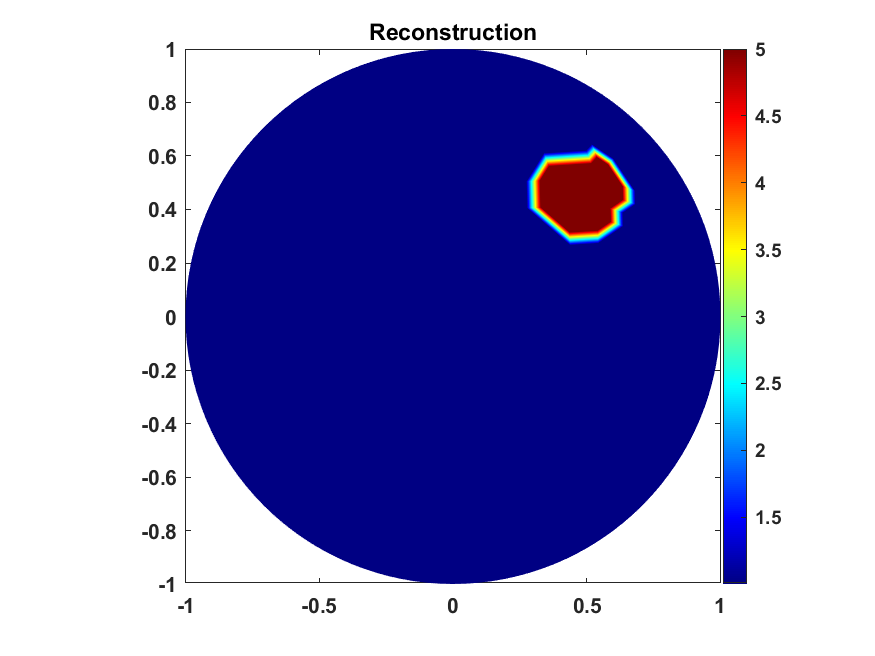}
                \caption{NL=6, EP=2000}
                \label{EIT_SIGAN1_L6_N16_MB128_EP2K}
         \end{subfigure}
          \begin{subfigure}[b]{0.24\textwidth}
                \includegraphics[scale=0.31]{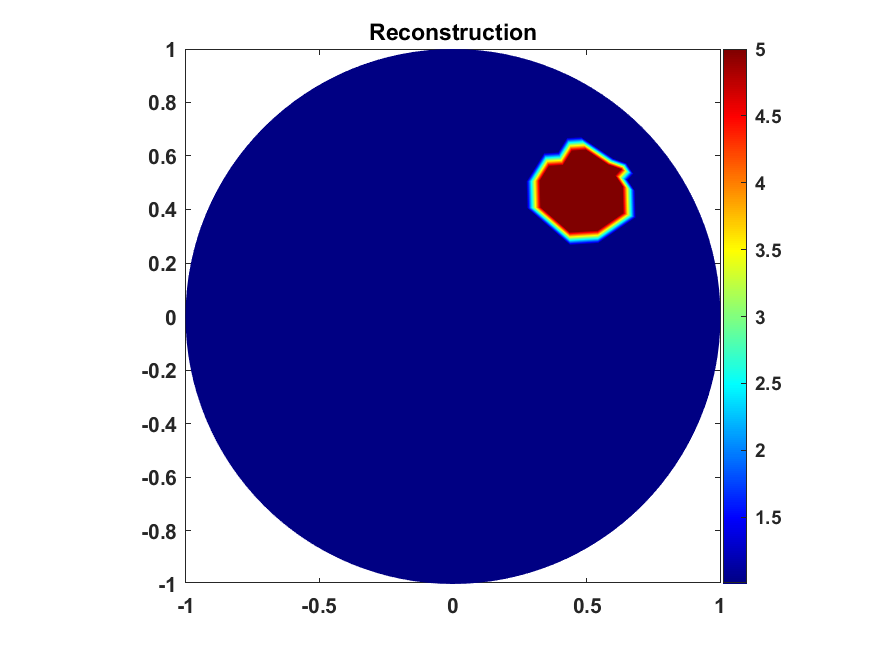}
               \caption{NL=8, EP=2000}
                \label{EIT_SIGAN1_L8_N16_MB128_EP2K}
         \end{subfigure}

       \begin{subfigure}[b]{0.24\textwidth}
          \includegraphics[scale=0.31]{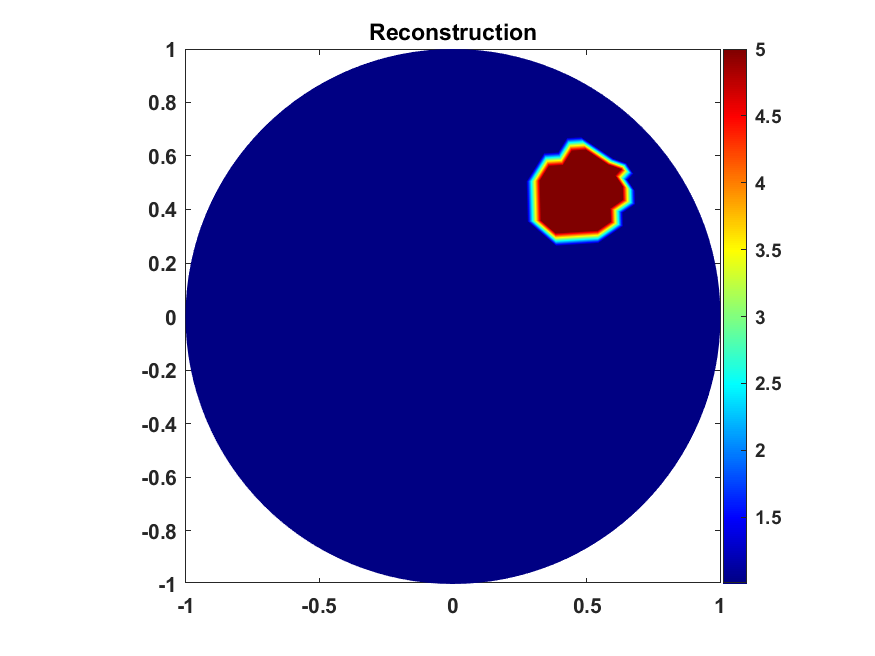}
                \caption{NL=4}
                \label{EIT_SIGAN1_L4_N32_MB128_EP2K}
         \end{subfigure}
         \begin{subfigure}[b]{0.24\textwidth}
        \includegraphics[scale=0.31]{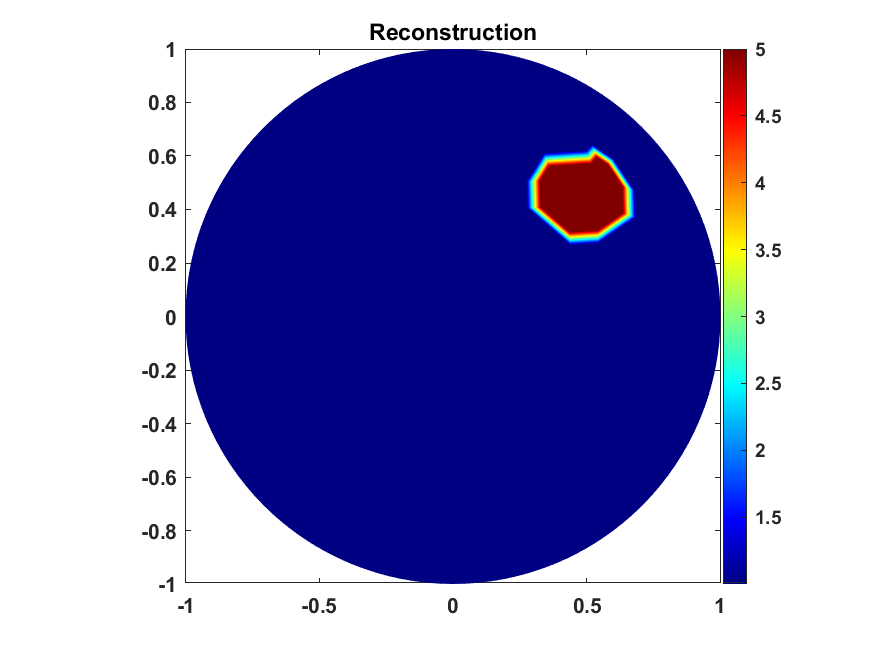}
                \caption{NL=6}
                \label{EIT_SIGAN1_L6_N32_MB128_EP2K}
         \end{subfigure}
          \begin{subfigure}[b]{0.24\textwidth}
        \includegraphics[scale=0.31]{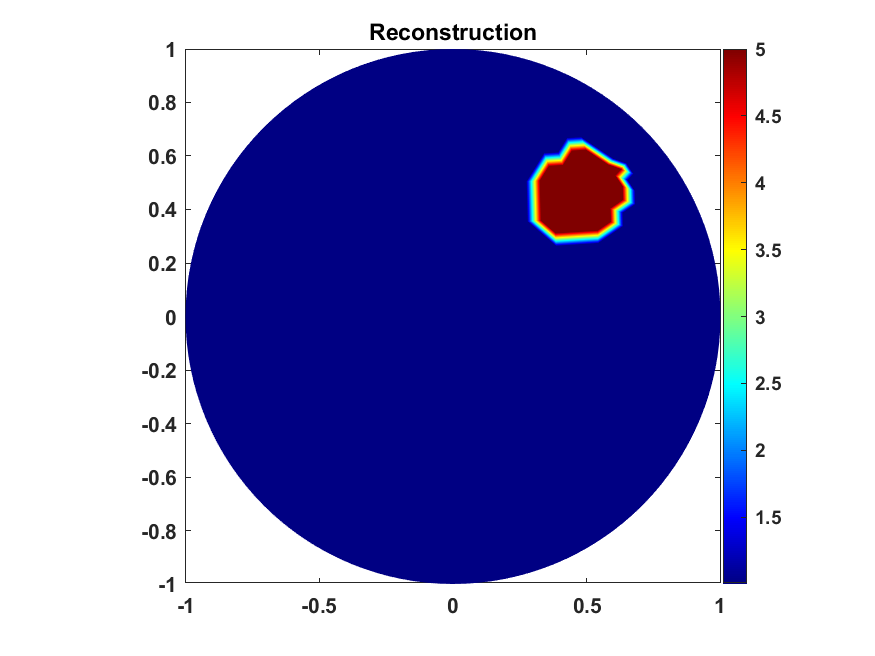}
               \caption{NL=8}
                \label{EIT_SIGAN1_L8_N32_MB128_EP2K}
         \end{subfigure}
        \begin{subfigure}[b]{0.24\textwidth}
          \includegraphics[scale=0.31]{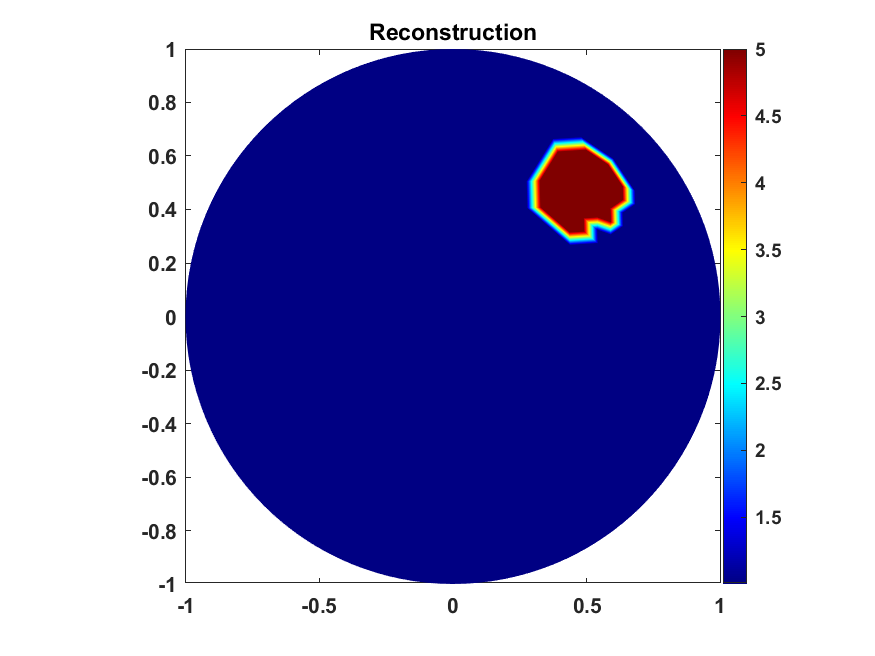}
                \caption{NL=4}
                \label{EIT_SIGAN1_L4_N64_MB128_EP2K}
         \end{subfigure}
         \begin{subfigure}[b]{0.24\textwidth}
                \includegraphics[scale=0.31]{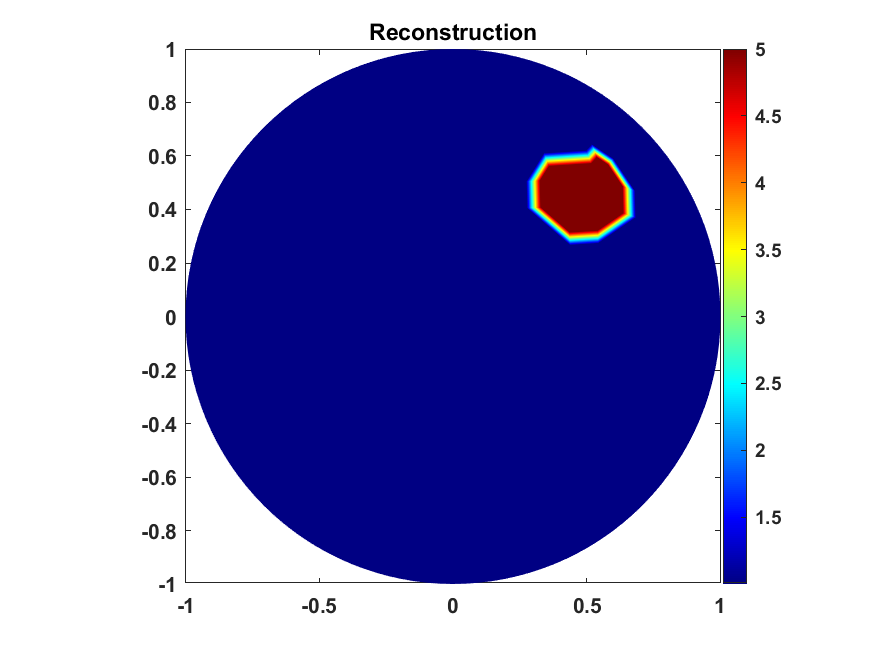}
                \caption{NL=6}
                \label{EIT_SIGAN1_L6_N64_MB128_EP2K}
         \end{subfigure}
          \begin{subfigure}[b]{0.24\textwidth}
                \includegraphics[scale=0.31]{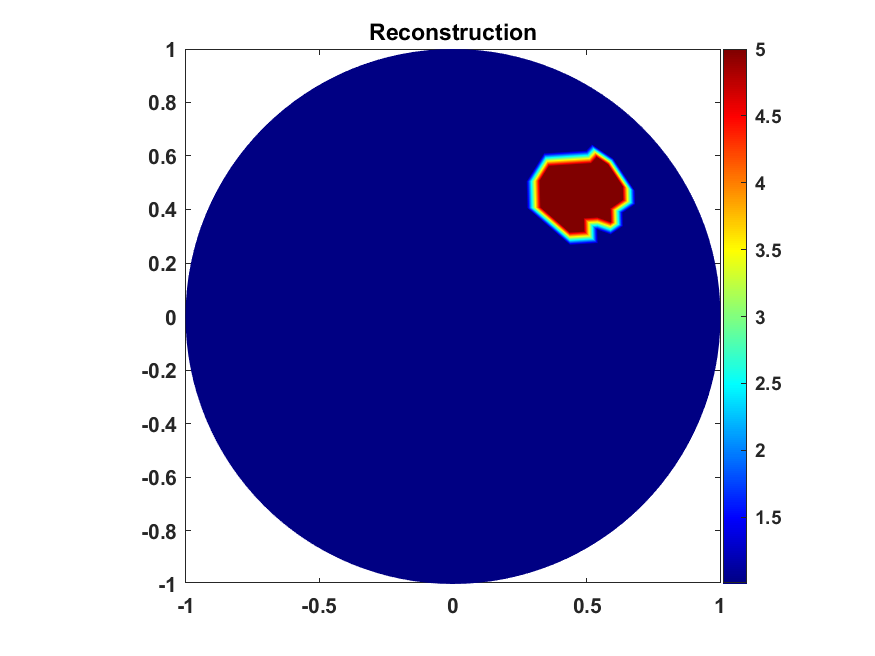}
               \caption{NL=8}
                \label{EIT_SIGAN1_L8_N64_MB128_EP2K}
         \end{subfigure}
        \begin{subfigure}[b]{0.24\textwidth}
          \includegraphics[scale=0.31]{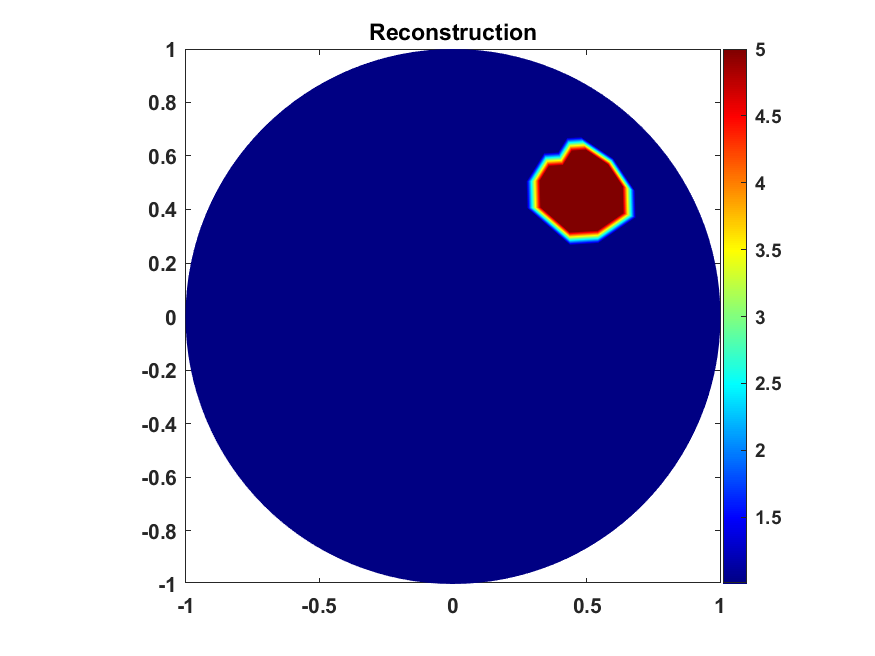}
                \caption{FEM}
                \label{EIT_SIGAN1_FEM_AN1_2024B}
         \end{subfigure}
             \begin{subfigure}[b]{0.24\textwidth}
          \includegraphics[scale=0.31]{sigma_eit_an1_true}
                \caption{True}
                \label{sigma_eit_an1_true1}
         \end{subfigure}

\caption{This figure shows qualitative results corresponding to the quantitative results in Table~\ref{tab:experiments_layers_neurons_mb_ep}, arranged top to bottom.
Figures~\ref{EIT_SIGAN1_L4_N16_MB64_EP2K}--\ref{EIT_SIGAN1_L8_N16_MB64_EP2K} show results for different numbers of layers (NL = 4, 6, 8) with the same number of neurons (NN = 16), minibatch size (MB = 64), and training epochs (EP = 2000). Figures~\ref{EIT_SIGAN1_L4_N16_MB128_EP1K}--\ref{EIT_SIGAN1_L8_N16_MB128_EP2K} present results for varying numbers of layers and epochs, with NN = 16 and MB = 128. Figures~\ref{EIT_SIGAN1_L4_N32_MB128_EP2K}--\ref{EIT_SIGAN1_L8_N32_MB128_EP2K} show results for different numbers of layers with NN = 32, MB = 128, and EP = 2000. Figures~\ref{EIT_SIGAN1_L4_N64_MB128_EP2K}--\ref{EIT_SIGAN1_L8_N64_MB128_EP2K} illustrate results for different numbers of layers with NN = 64, MB = 128, and EP = 2000. Figure~\ref{EIT_SIGAN1_FEM_AN1_2024B} illustrates the result using FEM-based MCMC inversion without using neural networks. Figure~\ref{sigma_eit_an1_true1} shows the simulated ground truth conductivity. Reconstructions are performed on a finite element mesh with 1345 grid points and 2560 triangles.}
\label{EITResultsDifferentHyperparameter}
\end{figure}

In our numerical experiments, we observed that the proposed method is highly robust to variations in the network architecture. We tabulate the result of one such experiment for EIT in Table \ref{tab:experiments_layers_neurons_mb_ep}. In conducting the sensitivity analysis for the network architecture and robustness assessment, we varied the network depth, the number of neurons in each layer, and the minibatch size and epochs during network training. We report the performance of the network using a variety of metrics, including training and inversion times, and the error incurred in the reconstruction using MAE and MSE. Similar tables for the DOT and QPAT set-ups are given in the Appendix; see Tables \ref{tab:dot_experiments_layers_neurons_mb_ep} and \ref{tab:qpat_experiments_layers_neurons_mb_ep}, respectively. We also include example reconstructions for this experiment in Figure \ref{EITResultsDifferentHyperparameter}. We highlight in bold the specific network configurations in Tables~\ref{tab:experiments_layers_neurons_mb_ep}, \ref{tab:dot_experiments_layers_neurons_mb_ep}, and \ref{tab:qpat_experiments_layers_neurons_mb_ep}, which were used as forward surrogates for the EIT, DOT, and QPAT inverse problems, respectively, in the main experiments presented in this paper.
\begin{table}
\caption{Quantitative results of conductivity reconstructions using MCMC-Net with varying numbers of layers (NL) and neurons (NN), trained under different minibatch sizes (MB) and epochs (EP) in the EIT set-up. Reported metrics include training time (Train. Time) and inversion time (Inv. Time), both in minutes, along with MAE and MSE. All reconstructions are based on MCMC sampling with 100,000 samples. The bottom row shows results from a traditional FEM-based MCMC inversion. The reconstructions correspond to a circular-shaped single-anomaly EIT inverse problem with simulated data corrupted by 1\% relative noise.}
\label{tab:experiments_layers_neurons_mb_ep}
\centering
\renewcommand{\arraystretch}{1.2} 

\scalebox{0.99}{
\begin{tabular}{|c|c|c|c|c|c|c|c|}
\hline
NL & NN & MB & EP & Train. Time (m) & Inv. Time (m) & MAE & MSE \\
\hline
4 & 16 & 64 & 2000 & 40.62 & 4.05 & 0.074371 & 0.297395 \\
6 & 16 & 64 & 2000 & 38.80 & 4.41 & 0.071397 & 0.285499 \\
8 & 16 & 64 & 2000 & 46.80 & 4.81 & 0.089240 & 0.356874 \\
\hline
4 & 16 & 128 & 1000 & 10.80 & 4.04 & 0.086267 & 0.344978 \\
\textbf{4} & \textbf{16} & \textbf{128} & \textbf{2000} & \textbf{17.58} & \textbf{4.09} & \textbf{0.068423} & \textbf{0.273603}\\
4 & 16 & 128 & 5000 & 58.87 & 4.08 & 0.059502 & 0.237916 \\
6 & 16 & 128 & 2000 & 20.97 & 4.45 & 0.080319 & 0.321186 \\
8 & 16 & 128 & 2000 & 23.02 & 4.72 & 0.068423 & 0.273603 \\
\hline
4 & 32 & 128 & 2000 & 19.53 & 4.06 & 0.065449 & 0.261707 \\
6 & 32 & 128 & 2000 & 24.83 & 4.45 & 0.077345 & 0.309291 \\
8 & 32 & 128 & 2000 & 28.36 & 4.61 & 0.065449 & 0.261707 \\
\hline
4 & 64 & 128 & 2000 & 68.98 & 4.61 & 0.074371 & 0.297395 \\
6 & 64 & 128 & 2000 & 94.18 & 5.12 & 0.077345 & 0.309291 \\
8 & 64 & 128 & 2000 & 119.27 & 5.01 & 0.083293 & 0.333082 \\
\hline
\multicolumn{5}{|c|}{\textbf{MCMC-FEM (no neural network)}} & 124.85 &  0.071397 & 0.285499 \\
\hline
\end{tabular}
}
\end{table}
\begin{table}[!h]
\caption{QPAT: Comparison of reconstruction errors and inversion time in minutes between MCMC-FEM and MCMC-Net at varying levels of noise in the data.}
\label{Reconstruction_errors_qpat}
\centering
\renewcommand{\arraystretch}{1.8} 

\scalebox{0.95}{
\begin{tabular}[t]{|l|cc|cc|cc|cc|}
\hline
\multirow{2}{*}{\shortstack{Anomaly type \\ Star $\&$ Kite } }   & \multicolumn{2}{c|}{MAE}  & \multicolumn{2}{c|}{MSE}  &
\multicolumn{2}{c|}{Inv. Time (m)}  \\
\cline{2-7}		
& \multicolumn{1}{c}{FEM} & \multicolumn{1}{c|}{Net} & 
\multicolumn{1}{c}{FEM} & \multicolumn{1}{c|}{Net} & \multicolumn{1}{c}{FEM} & \multicolumn{1}{c|}{Net} \\
\hline
 Noise Level 2 \%   & 0.001329 & 0.001708  & 0.000101 & 0.000156 & 10.62 & 6.61\\
 Noise Level 4 \%   & 0.001467 & 0.001908  & 0.000102 & 0.000149 & 10.76 & 6.48  \\
 Noise Level 8 \%   & 0.002121 & 0.002417  & 0.000195 & 0.000242 & 10.66 & 6.52 \\
\hline      
\end{tabular}
}
\end{table}

\subsection{Noise sensitivity analysis} Having established the robustness of the reconstruction process with respect to the network architecture and training hyperparameters in the preceding subsection, we now describe the impact of noise in the data on reconstruction quality. In Table \ref{Reconstruction_errors_qpat}, we show the results for reconstructions carried out for the inverse problem in QPAT using the network architecture as described in subsection \ref{sec:5.3} for varying levels of noise in the data. As expected, reconstruction errors increase with higher noise levels for both FEM-based and neural operator-based approaches. However, in terms of computational time, neural operator-based reconstructions maintain a clear advantage while outperforming the traditional method in terms of the errors incurred. Corresponding results for the EIT and DOT inverse problems are provided in the Appendix; see Tables~\ref{tab:noise_inversion_mae_mse_grouped} and~\ref{tab:dot_noise_inversion_mae_mse_grouped}, respectively.

\begin{figure}[!h]
    \centering 
\includegraphics[height=0.3\textheight]{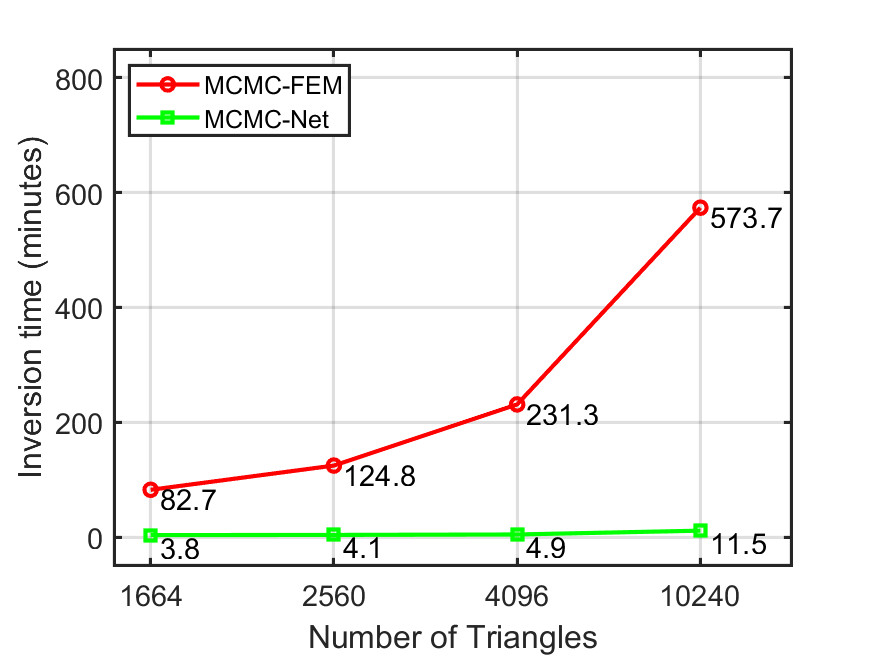}

\caption{
Comparison of inversion times between MCMC-FEM and MCMC-Net. For each method, 100{,}000 posterior samples were generated using four different mesh resolutions corresponding to 1,664; 2,560; 4,096; and 10,240 triangular elements. Note the exponential increase in inversion time for the MCMC-FEM solvers vs the  linear increase in inversion time for the MCMC-Net solver. The results are computed for the EIT inverse problem with a single circular-shaped anomaly.
}
\label{MCMCCNNFEM_Time_Comparison}
\end{figure}

\subsection{Speed-up due to MCMC vs dimension of the discretized parameter space}
One of the significant advantages of using the MCMC-Net over the traditional solver is the speed-up that is achieved when carrying out the inversion with a trained network. The speed-up is significant in high-dimensional inverse problems that we consider here. Figure~\ref{MCMCCNNFEM_Time_Comparison} illustrates the inversion times of MCMC-Net compared to the traditional FEM-based solver for the EIT problem vis-a-vis the discretization of the conductivity function.  We note that as the dimension of the discretized representation of the parameter space increases, the inversion time of the MCMC-FEM solver increases exponentially, while it increases only linearly for the MCMC-Net implementation. This shows that MCMC-Net is increasingly computationally efficient when the dimension of the discrete parameter space increases.



\section{Conclusion}
\label{conclusion}
In this work, we introduce a deep neural network based surrogate model that significantly reduces the computational cost of likelihood evaluation in MCMC for Bayesian inverse problems. We further analyze a universal approximation theorem-type result to show that the proposed network structure can approximate the respective forward maps in a suitable asymptotic sense. This in turns allows us to establish an asymptotic result that shows that the surrogate posterior derived by using the MCMC-Net approaches the posterior using the traditional FEM forward solver in pCN based sampler while at the same time cutting down computational cost by an order of magnitude which justifies using the proposed method. We validated this approach on three distinct inverse problems and demonstrated that the proposed surrogate technique achieves up to thirty times the speed of traditional methods when considering posterior samples up to 100,000.  Hence, the proposed technique is an important addition to the field of Bayesian inverse problems, offering both computational efficiency and reliability. We do however caution that the use of such neural network surrogates is only justified in presence of reliable training data. Finally, a more detailed analysis incorporating training and generalization error bounds is left for future work.
 

\clearpage




\appendix
\newpage

\section{Sensitivity analysis of Network architectures for DOT and QPAT inversions.}
\begin{figure}[!h]
    \centering 
    \begin{subfigure}[b]{0.24\textwidth}
        \includegraphics[scale=0.31]{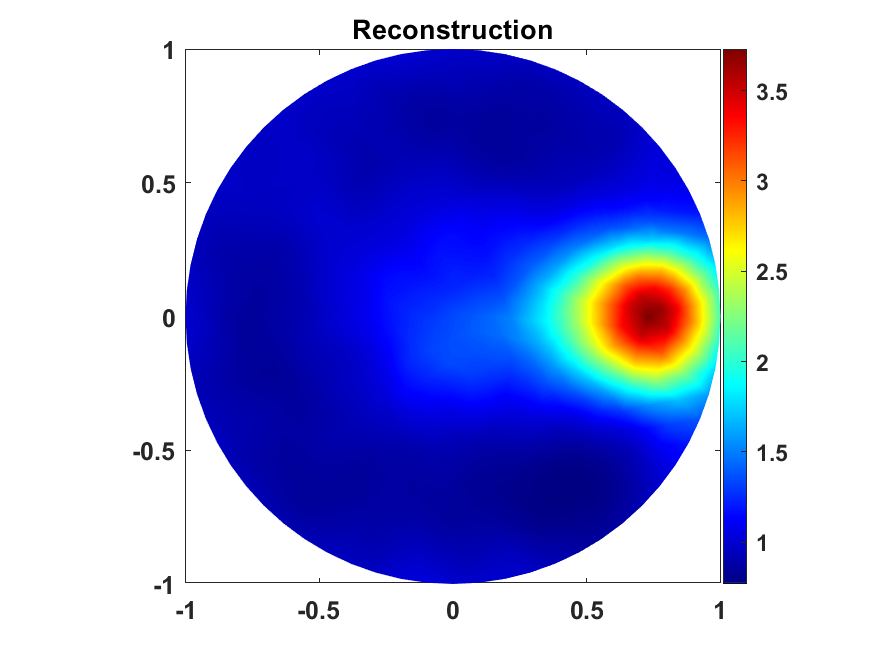}
        \caption{NL=4}
        \label{DOT_mu_CNN_AN1_nL1_L4_N16_MB8_EP100}
    \end{subfigure}
    \begin{subfigure}[b]{0.24\textwidth}
        \includegraphics[scale=0.31]{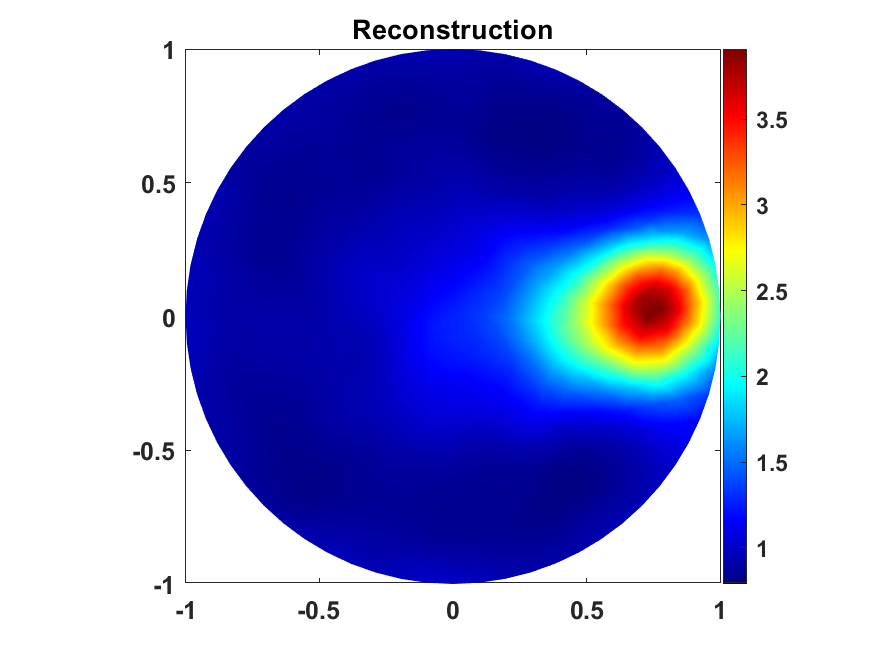}
        \caption{NL=4, MB=16}
        \label{DOT_mu_CNN_AN1_nL1_L4_N16_MB16_EP100}
    \end{subfigure}
        \begin{subfigure}[b]{0.24\textwidth}
        \includegraphics[scale=0.31]{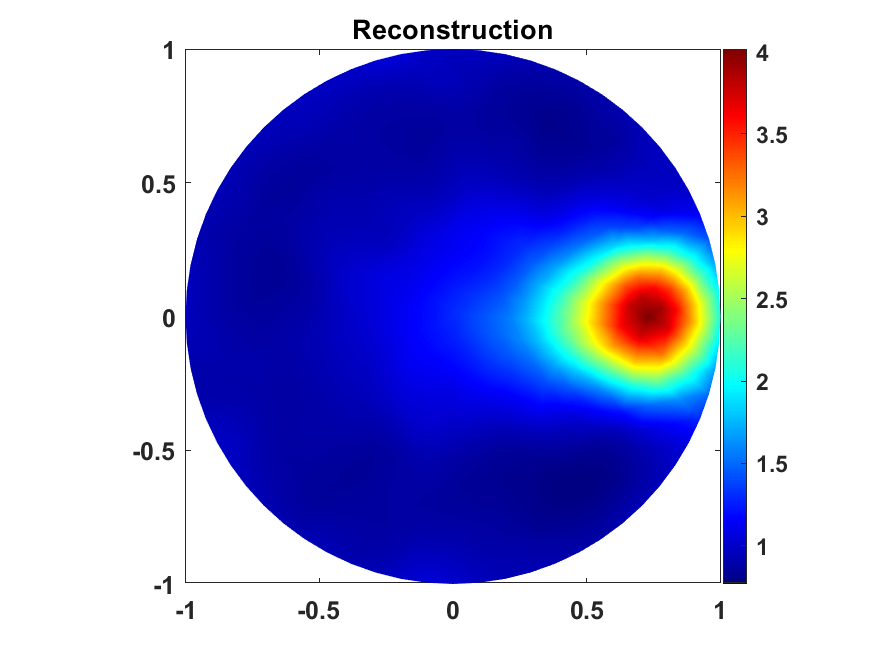}
        \caption{NL=6}
        \label{DOT_mu_CNN_AN1_nL1_L6_N16_MB8_EP100}
    \end{subfigure}
        \begin{subfigure}[b]{0.24\textwidth}
        \includegraphics[scale=0.31]{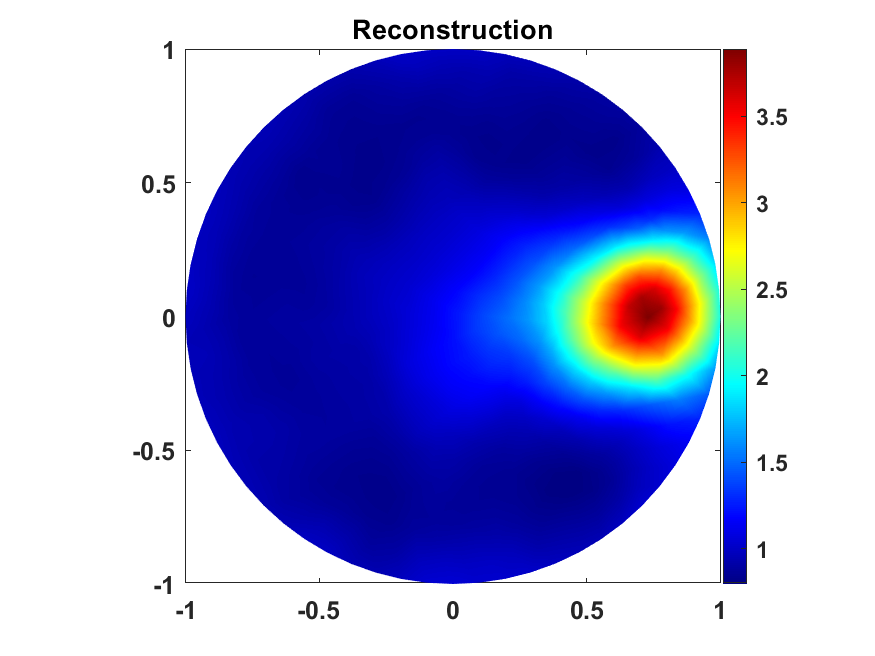}
        \caption{NL=8}
        \label{DOT_mu_CNN_AN1_nL1_L8_N16_MB8_EP100}
        \end{subfigure}

    \begin{subfigure}[b]{0.24\textwidth}
        \includegraphics[scale=0.31]{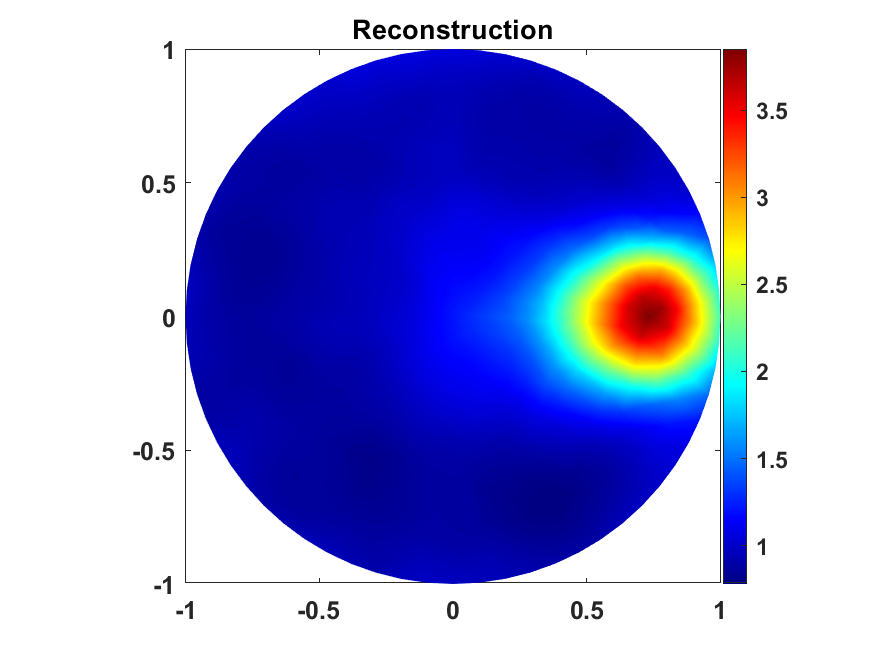}
        \caption{NL=4}
        \label{DOT_mu_CNN_AN1_nL1_L4_N16_MB8_EP200}
    \end{subfigure}
    \begin{subfigure}[b]{0.24\textwidth}
        \includegraphics[scale=0.31]{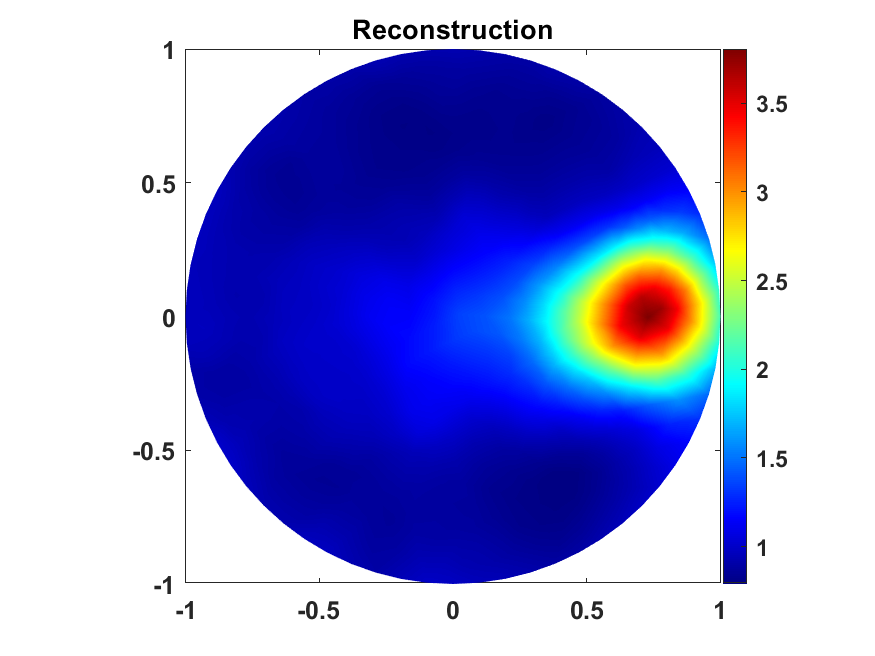}
        \caption{NL=4, MB=16}
        \label{DOT_mu_CNN_AN1_nL1_L4_N16_MB16_EP200}
    \end{subfigure}
        \begin{subfigure}[b]{0.24\textwidth}
        \includegraphics[scale=0.31]{mu_CNN_AN1_nL1_L6_N16_MB8_EP100}
        \caption{NL=6}
        \label{DOT_mu_CNN_AN1_nL1_L6_N16_MB8_EP200}
    \end{subfigure}
        \begin{subfigure}[b]{0.24\textwidth}
        \includegraphics[scale=0.31]{mu_CNN_AN1_nL1_L8_N16_MB8_EP100}
        \caption{NL=8}
        \label{DOT_mu_CNN_AN1_nL1_L8_N16_MB8_EP200}
        \end{subfigure}

       \begin{subfigure}[b]{0.24\textwidth}
        \includegraphics[scale=0.31]{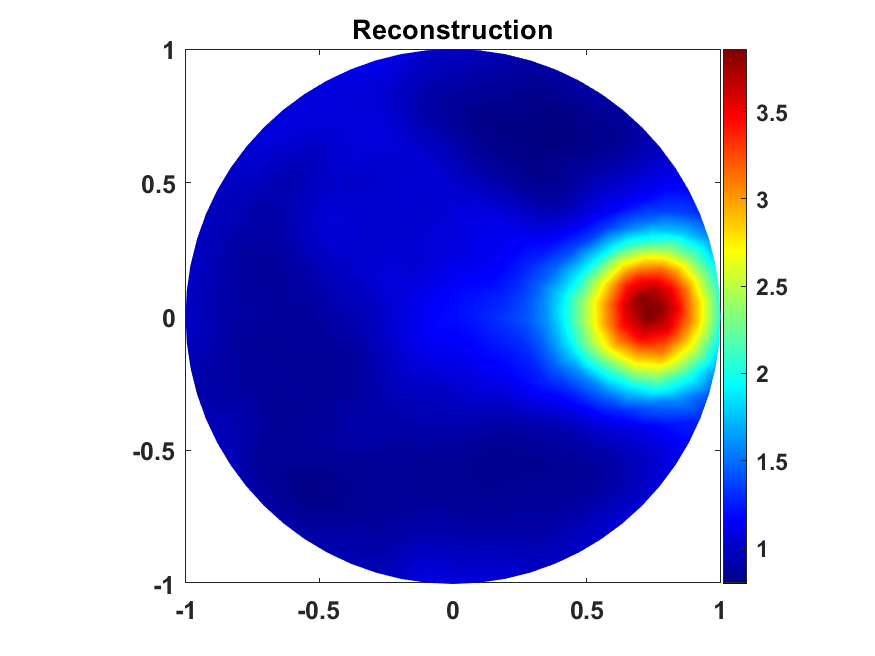}
        \caption{NL=4}
        \label{DOT_mu_CNN_AN1_nL1_L4_N32_MB8_EP100}
    \end{subfigure}
    \begin{subfigure}[b]{0.24\textwidth}
        \includegraphics[scale=0.31]{mu_CNN_AN1_nL1_L4_N32_MB8_EP100}
        \caption{NL=6}
        \label{DOT_mu_CNN_AN1_nL1_L6_N32_MB16_EP100}
    \end{subfigure}
        \begin{subfigure}[b]{0.24\textwidth}
        \includegraphics[scale=0.31]{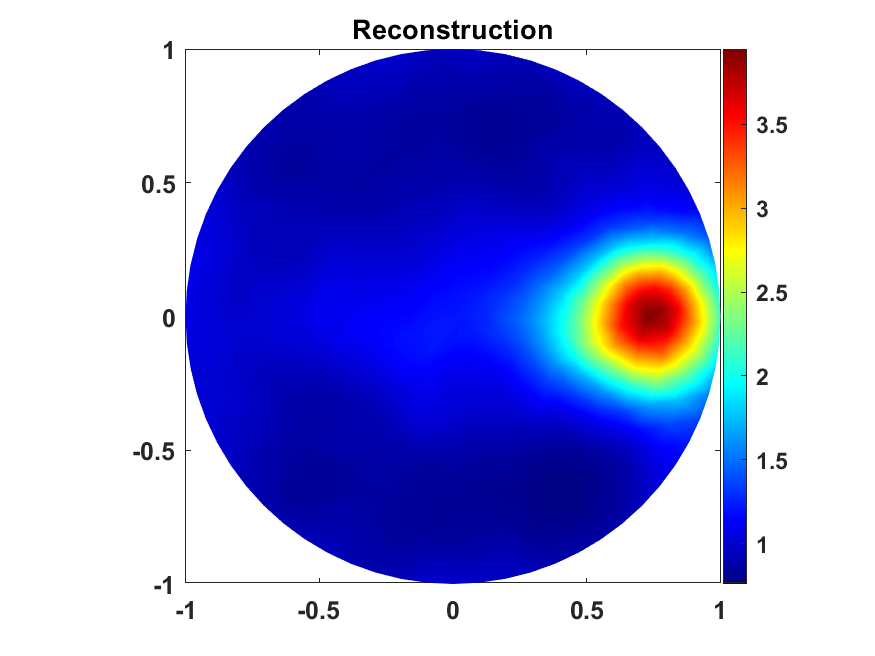}
        \caption{NL=8}
        \label{DOT_mu_CNN_AN1_nL1_L8_N32_MB8_EP100}
    \end{subfigure}
        \begin{subfigure}[b]{0.24\textwidth}
        \includegraphics[scale=0.31]{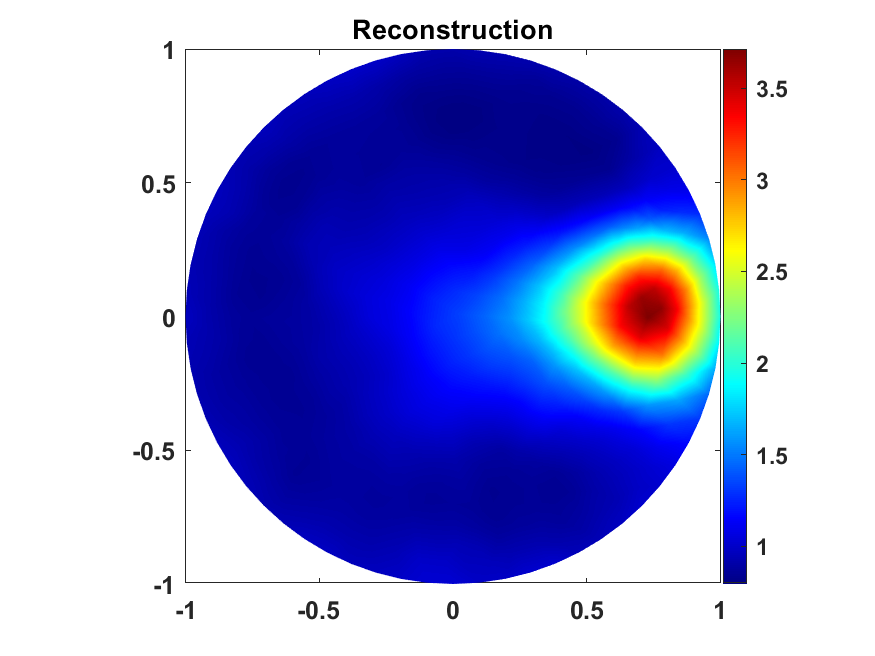}
        \caption{NL=4}
        \label{DOT_mu_CNN_AN1_nL1_L4_N64_MB8_EP100}
        \end{subfigure}

        \begin{subfigure}[b]{0.24\textwidth}
        \includegraphics[scale=0.31]{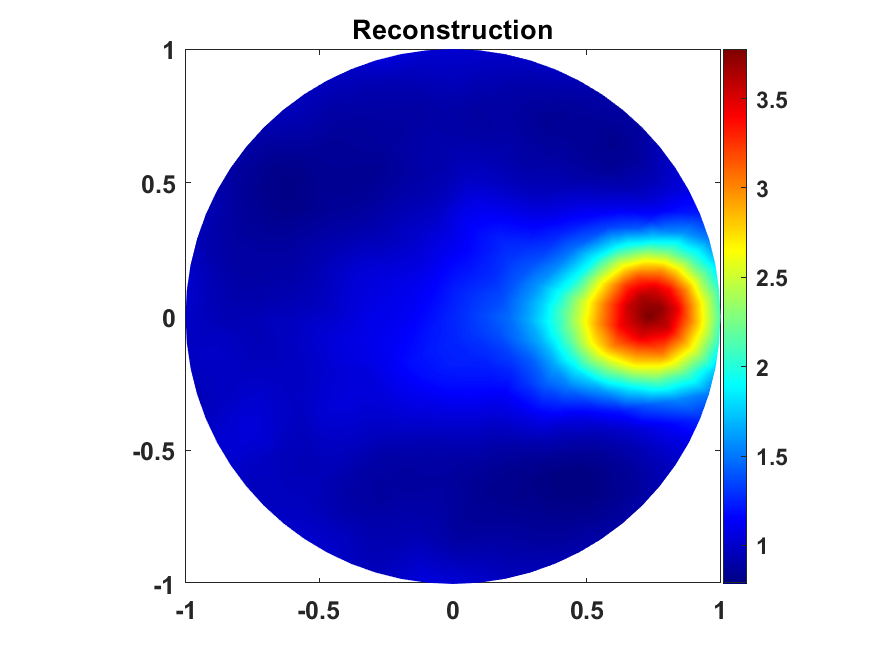}
        \caption{NL=6}
        \label{DOT_mu_CNN_AN1_nL1_L6_N64_MB8_EP100}
    \end{subfigure}
    \begin{subfigure}[b]{0.24\textwidth}
        \includegraphics[scale=0.31]{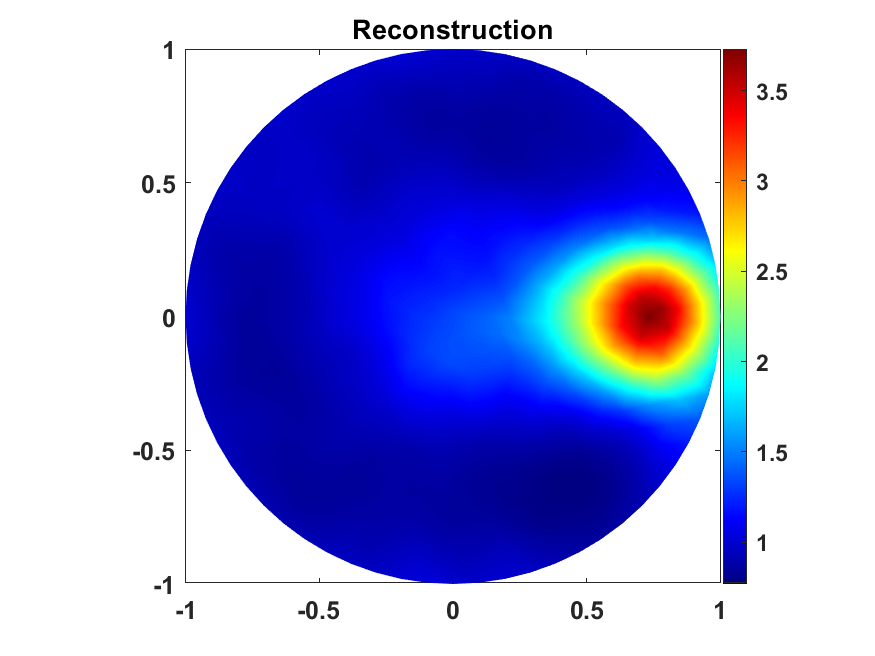}
        \caption{NL=8}
        \label{DOT_mu_CNN_AN1_nL1_L8_N64_MB16_EP100}
    \end{subfigure}
         \begin{subfigure}[b]{0.24\textwidth}
        \includegraphics[scale=0.31]{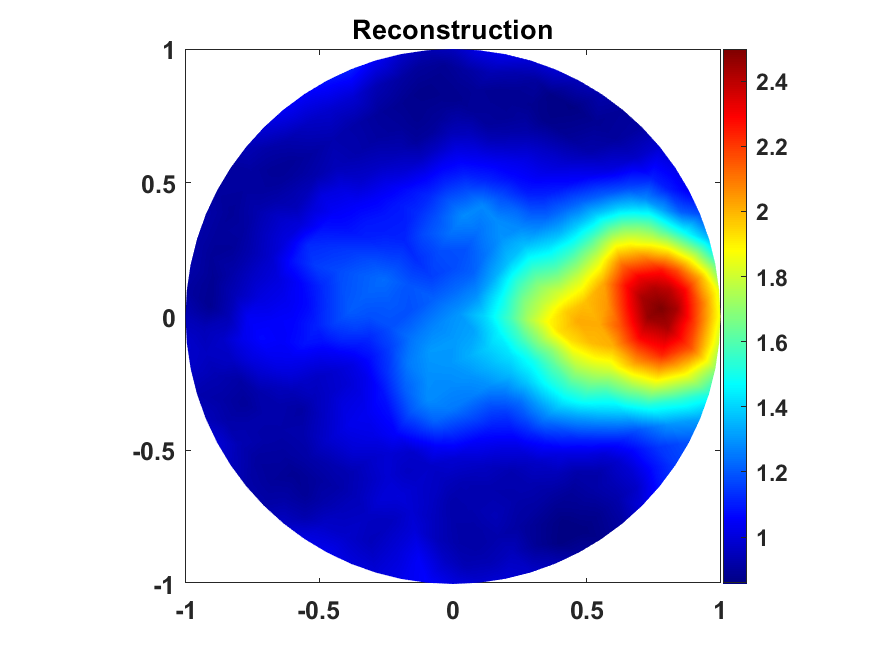}
        \caption{MCMC-FEM}
        \label{DOT_mu_FEM_AN1_nL1_2024b}
        \end{subfigure}
        \begin{subfigure}[b]{0.24\textwidth}
        \includegraphics[scale=0.31]{mu_dot_an1_truefine}
        \caption{True}
        \label{DOT_mu_an1_truefine1}
    \end{subfigure}
   
\caption{This figure shows qualitative results corresponding to the quantitative results in Table~\ref{tab:dot_experiments_layers_neurons_mb_ep}, arranged top to bottom.
Figures~\ref{DOT_mu_CNN_AN1_nL1_L4_N16_MB8_EP100}--\ref{DOT_mu_CNN_AN1_nL1_L8_N16_MB8_EP100} show results for different numbers of layers (NL = 4, 6, 8) with the same number of neurons (NN = 16), trained for 100 epochs. Figures~\ref{DOT_mu_CNN_AN1_nL1_L4_N16_MB8_EP100}--\ref{DOT_mu_CNN_AN1_nL1_L8_N16_MB8_EP100} show results for different numbers of layers (NL = 4, 6, 8) with the same number of neurons (NN = 16), trained for 200 epochs. Figures~\ref{DOT_mu_CNN_AN1_nL1_L4_N32_MB8_EP100}--\ref{DOT_mu_CNN_AN1_nL1_L8_N32_MB8_EP100} show results for different numbers of layers (NL = 4, 6, 8) with the same number of neurons (NN = 32), trained for 100 epochs. Figures~\ref{DOT_mu_CNN_AN1_nL1_L4_N64_MB8_EP100}--\ref{DOT_mu_CNN_AN1_nL1_L8_N64_MB16_EP100} show results for different numbers of layers (NL = 4, 6, 8) with the same number of neurons (NN = 64), trained for 100 epochs. Figure~\ref{DOT_mu_FEM_AN1_nL1_2024b} illustrates the result using traditional FEM-based MCMC inversion without using neural networks. Figure~\ref{DOT_mu_an1_truefine1} shows the simulated ground truth conductivity. Reconstructions are performed on a finite element mesh consisting of 549 grid points and 1032 triangles.}
\label{DOTResultsDifferentHyperparameter}
\end{figure}

\begin{figure}[h!]
    \centering 
    \begin{subfigure}[b]{0.28\textwidth}
        \includegraphics[scale=0.35]{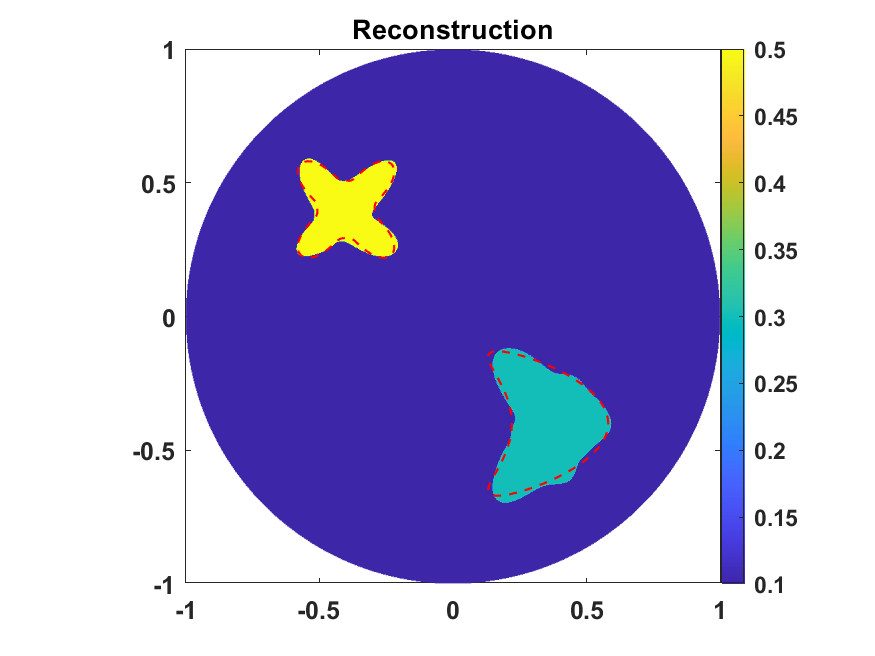}
        \caption{NL=4, NN=16}
        \label{gamma_QPAT_nL1_L4_N16_EP100_LD20_MB8}
    \end{subfigure}
    \begin{subfigure}[b]{0.28\textwidth}
        \includegraphics[scale=0.35]{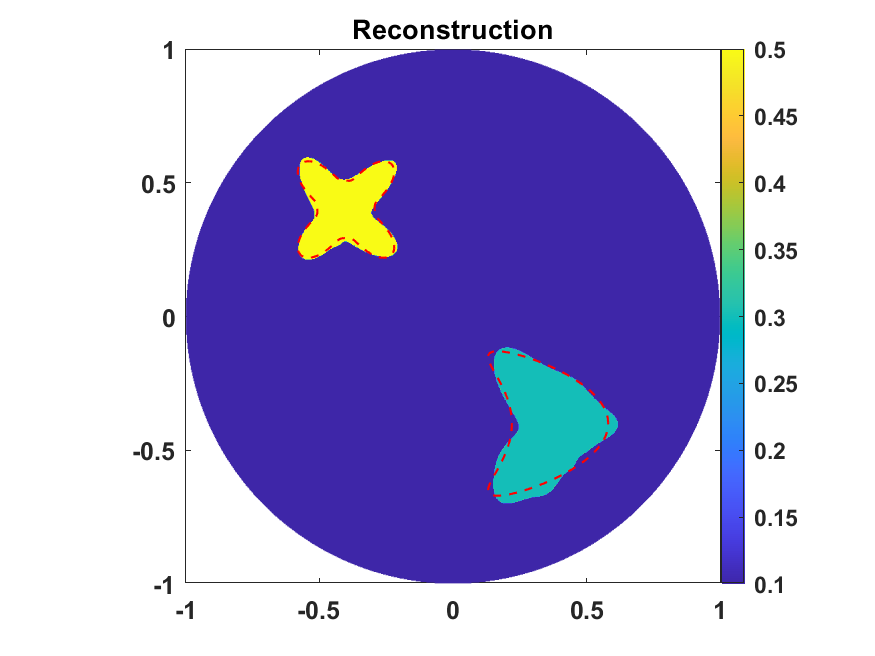}
        \caption{NL=4, NN=32}
        \label{gamma_QPAT_nL1_L4_N32_EP100_LD20_MB8}
    \end{subfigure}
        \begin{subfigure}[b]{0.28\textwidth}
        \includegraphics[scale=0.35]{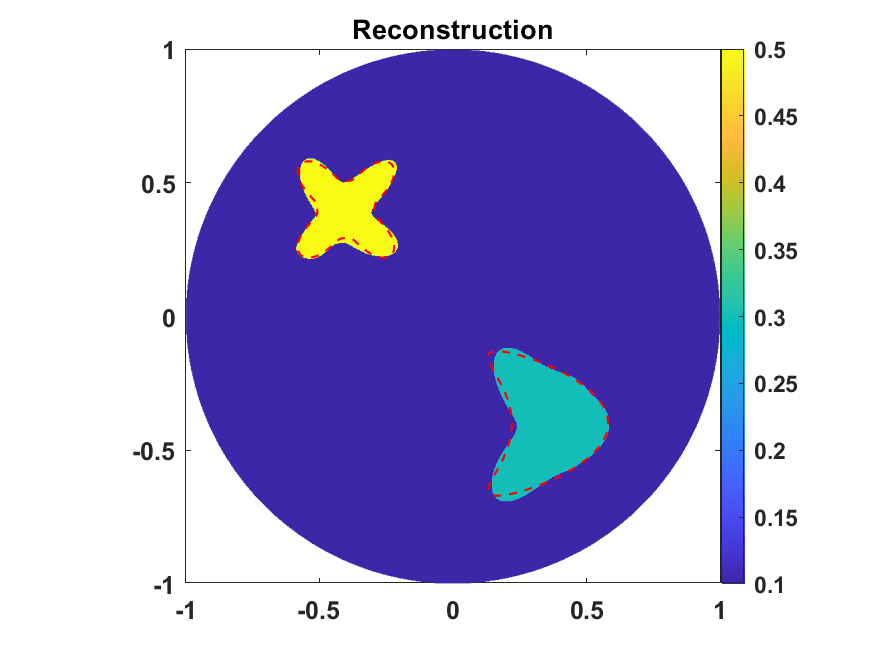}
        \caption{NL=4, NN=64}
        \label{gamma_QPAT_nL1_L4_N64_EP100_LD20_MB8}
    \end{subfigure}

    \begin{subfigure}[b]{0.28\textwidth}
        \includegraphics[scale=0.35]{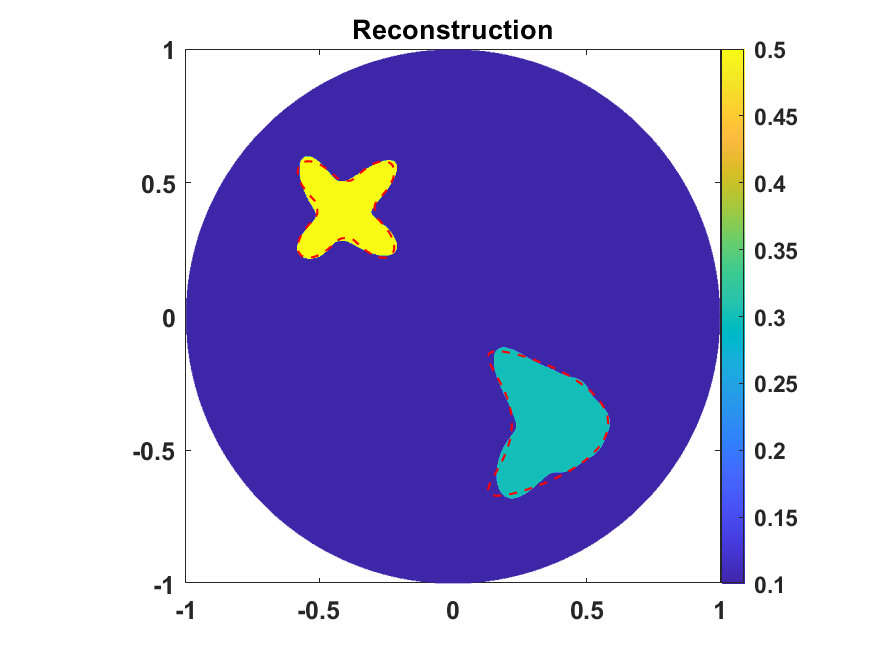}
        \caption{NL=6, NN=16}
        \label{gamma_QPAT_nL1_L6_N16_EP100_LD20_MB8}
    \end{subfigure}
    \begin{subfigure}[b]{0.28\textwidth}
        \includegraphics[scale=0.35]{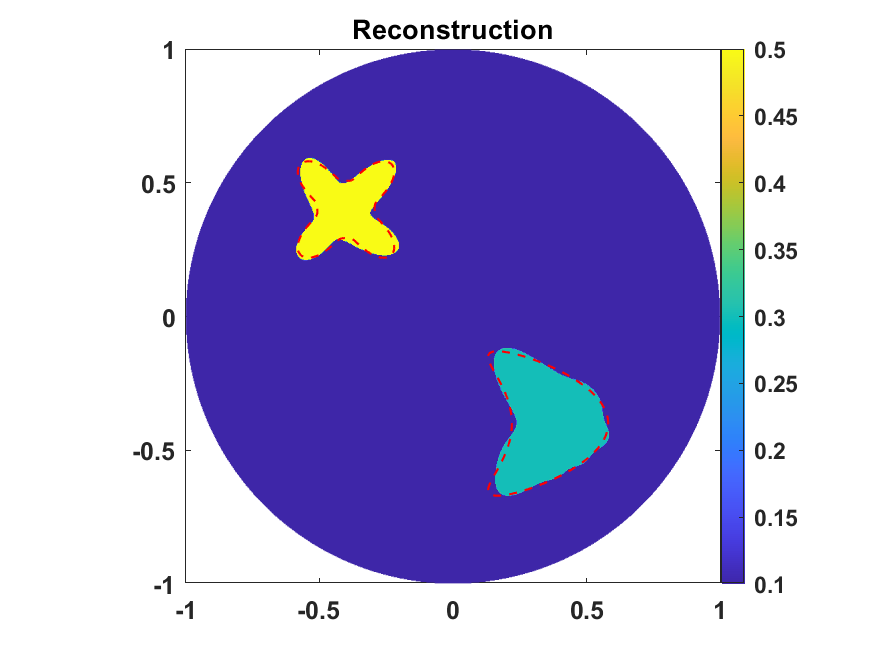}
        \caption{NL=6, NN=32}
        \label{gamma_QPAT_nL1_L6_N32_EP100_LD20_MB8}
    \end{subfigure}
        \begin{subfigure}[b]{0.28\textwidth}
        \includegraphics[scale=0.35]{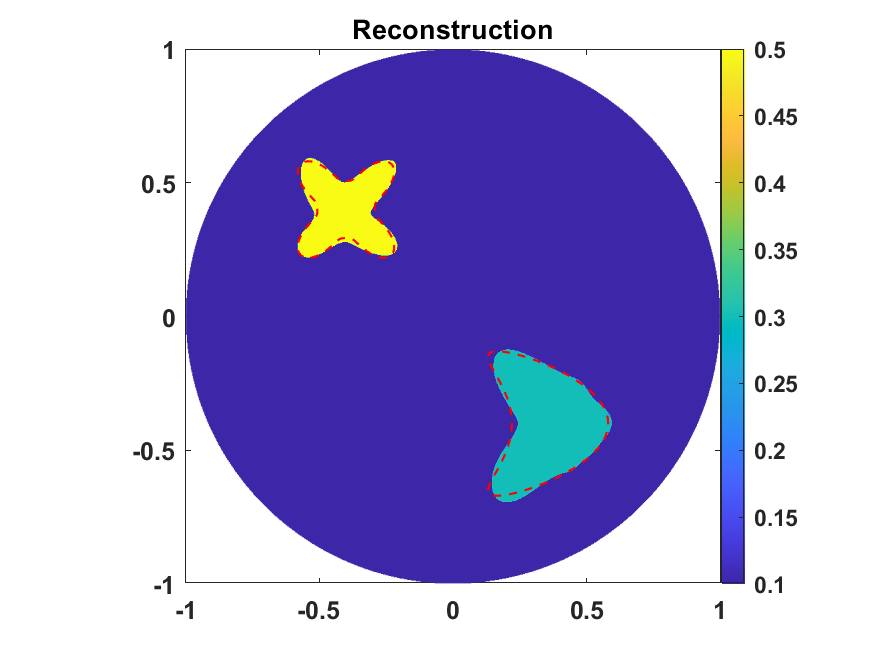}
        \caption{NL=6, NN=64}
        \label{gamma_QPAT_nL1_L6_N64_EP100_LD20_MB8}
    \end{subfigure}

    \begin{subfigure}[b]{0.28\textwidth}
        \includegraphics[scale=0.35]{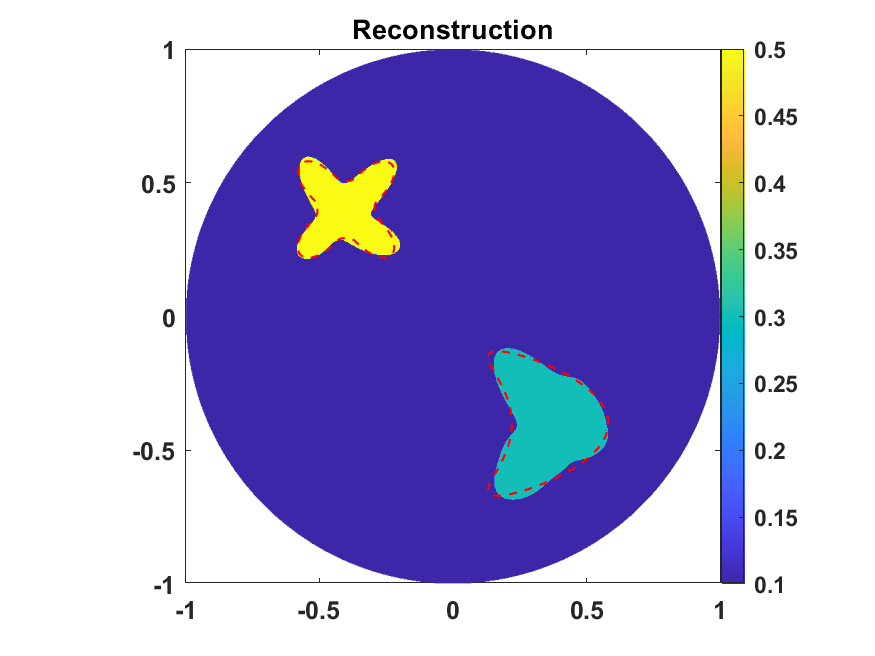}
        \caption{NL=8, NN=16}
        \label{gamma_QPAT_nL1_L8_N16_EP100_LD20_MB8}
    \end{subfigure}
    \begin{subfigure}[b]{0.28\textwidth}
        \includegraphics[scale=0.35]{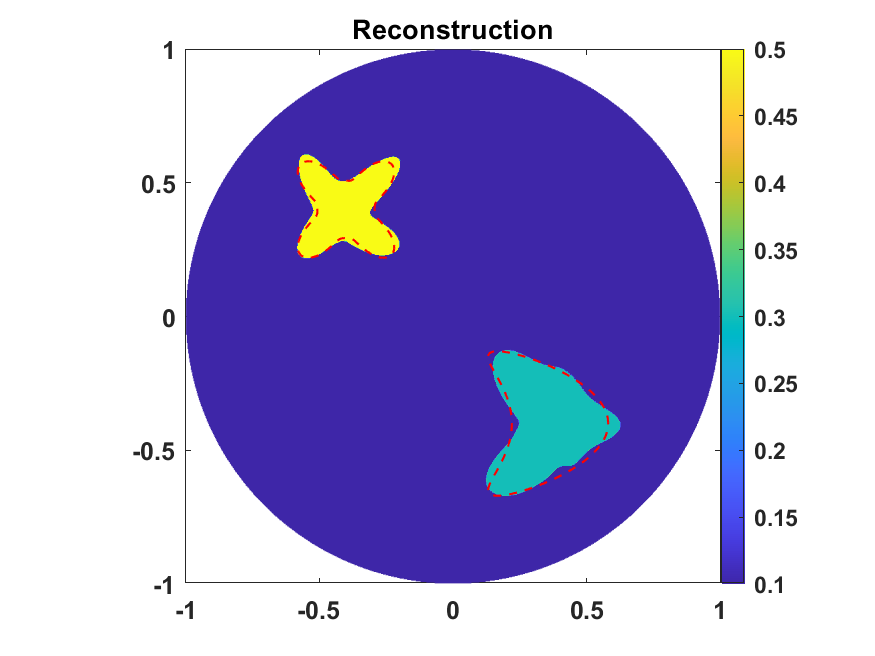}
        \caption{NL=8, NN=32}
        \label{gamma_QPAT_nL1_L8_N32_EP100_LD20_MB8}
    \end{subfigure}
        \begin{subfigure}[b]{0.28\textwidth}
        \includegraphics[scale=0.35]{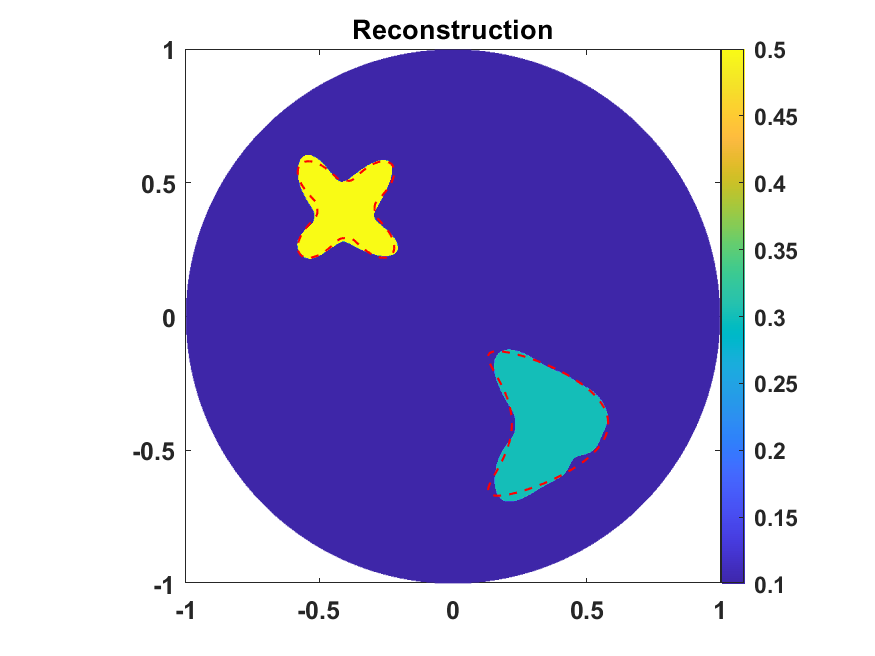}
        \caption{NL=8, NN=64}
        \label{gamma_QPAT_nL1_L8_N64_EP100_LD20_MB8}
    \end{subfigure}

    \begin{subfigure}[b]{0.28\textwidth}
        \includegraphics[scale=0.35]{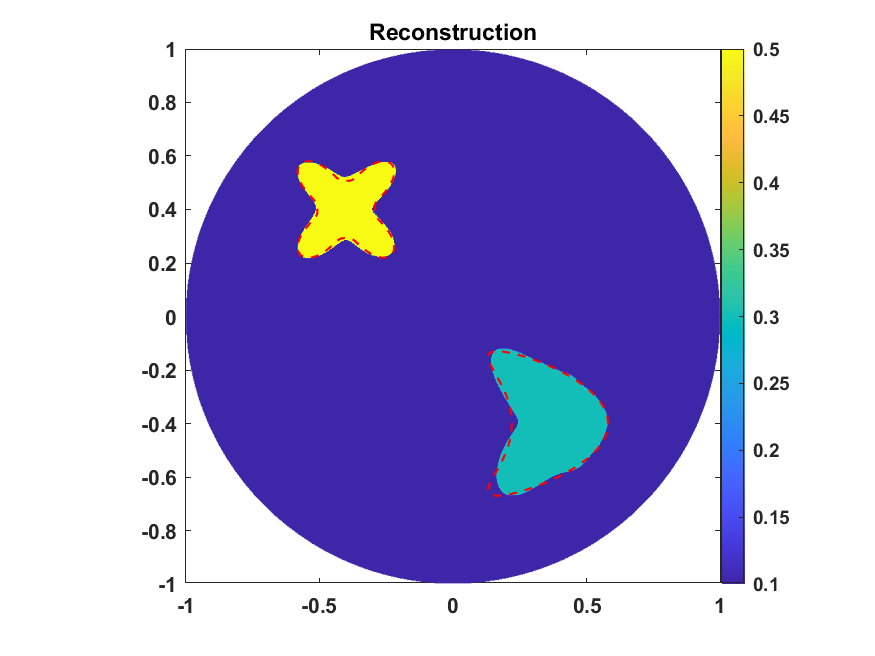}
        \caption{MCMC-FEM}
        \label{gamma_QPAT_FEM_nL1_ITER_100k}
    \end{subfigure}
        \begin{subfigure}[b]{0.28\textwidth}
        \includegraphics[scale=0.35]{gammatrue_qPAT}
        \caption{True}
        \label{gamma_QPAT_true1}
    \end{subfigure}

\caption{This figure shows qualitative results corresponding to the quantitative results in Table~\ref{tab:qpat_experiments_layers_neurons_mb_ep}, arranged top to bottom.
Figures~\ref{gamma_QPAT_nL1_L4_N16_EP100_LD20_MB8}--\ref{gamma_QPAT_nL1_L4_N64_EP100_LD20_MB8} show results for the same number of layers (NL = 4) with different numbers of neurons (NN = 16, 32, 64), trained for 100 epochs with a minibatch size of 8. Figures~\ref{gamma_QPAT_nL1_L6_N16_EP100_LD20_MB8}--\ref{gamma_QPAT_nL1_L6_N64_EP100_LD20_MB8} show results for the same number of layers (NL = 4) with different numbers of neurons (NN = 16, 32, 64), trained for 100 epochs with a minibatch size of 8. Figures~\ref{gamma_QPAT_nL1_L8_N16_EP100_LD20_MB8}--\ref{gamma_QPAT_nL1_L8_N64_EP100_LD20_MB8} show results for the same number of layers (NL = 4) with different numbers of neurons (NN = 16, 32, 64), trained for 100 epochs with a minibatch size of 8.
Figure~\ref{gamma_QPAT_FEM_nL1_ITER_100k} illustrates the result using traditional FEM-based MCMC inversion without using neural networks. Figure~\ref{gamma_QPAT_true1} shows the simulated ground truth conductivity. Reconstructions are performed on a finite element mesh consisting of 1515 grid points and 2904 triangles.}
\label{QPATResultsDifferentHyperparameter}
\end{figure}

Earlier, we tabulated the results for how sensitive the Bayesian inversion in the EIT is to the chosen network architecture and saw the method to be quite robust. In the tables \ref{tab:dot_experiments_layers_neurons_mb_ep} and \ref{tab:qpat_experiments_layers_neurons_mb_ep} below, we produce the results of a similar analysis carried out for the DOT and QPAT inversions. Example reconstructions are given in Figures \ref{DOTResultsDifferentHyperparameter} and \ref{QPATResultsDifferentHyperparameter}, respectively.

\begin{table}[!h]
\caption{Quantitative results of conductivity reconstructions using MCMC-Net with varying numbers of layers (NL) and neurons (NN), trained under different minibatch sizes (MB) and epochs (EP) for the DOT inverse problem. Reported metrics include training time (Train. Time) and inversion time (Inv. Time), both in minutes, along with the $L^{\infty}$ -loss, mean absolute error (MAE), and mean squared error (MSE). All reconstructions are based on MCMC sampling with 100,000 samples. The bottom row shows results from a traditional FEM-based MCMC inversion. The reconstructions correspond to a circular-shaped single-anomaly DOT inverse problem with simulated data corrupted by 1\% relative noise.}
\label{tab:dot_experiments_layers_neurons_mb_ep}
\centering
\renewcommand{\arraystretch}{1.2} 

\scalebox{0.94}{
\begin{tabular}{|c|c|c|c|c|c|c|c|c|}
\hline
NL & NN & MB & EP & Train. Time (m) & Inv. Time (m) & $L^{\infty}$ -loss & MAE & MSE \\
\hline
\textbf{4} & \textbf{16} & \textbf{8} & \textbf{100} & \textbf{14.36} & \textbf{3.59} & \textbf{2.343760} & \textbf{0.258924} & \textbf{0.241630} \\
4 & 16 & 16 & 100 & 7.68 & 3.63 & 2.337403   & 0.269474 & 0.241718 \\
6 & 16 & 8 & 100 & 13.90 & 3.97 & 2.180491   & 0.259781 & 0.231118 \\
8 & 16 & 8 & 100 & 16.20 & 4.37 & 2.216273   & 0.255994 & 0.230677 \\
\hline
4 & 16 & 8 & 200 & 26.62 & 3.66 & 2.301730  & 0.248475 & 0.228772 \\
4 & 16 & 16 & 200 & 12.75 & 3.62 & 2.324823 & 0.261933 & 0.241714 \\
6 & 16 & 8 & 200 & 28.71 & 3.92 & 2.329720  & 0.242877 & 0.223699 \\
8 & 16 & 8 & 200 & 31.41 & 4.31 & 2.321498  & 0.265998 & 0.239659 \\
\hline
4 & 32 & 8 & 100 & 12.43 & 3.55 & 2.292029   & 0.238404 & 0.222464 \\
6 & 32 & 8 & 100 & 14.57 & 3.99 & 2.244161   & 0.265378 & 0.247470 \\
8 & 32 & 8 & 100 & 16.12 & 4.30 & 2.450117   & 0.245358 & 0.224046 \\
\hline
4 & 64 & 8 & 100 & 18.92 & 4.08 & 2.259012   & 0.264918 & 0.243508 \\
6 & 64 & 8 & 100 & 21.59 & 4.45 & 2.287532   & 0.248400 & 0.235426 \\
8 & 64 & 8 & 100 & 23.40 & 4.76 & 2.388119   & 0.242490 & 0.230856 \\
\hline
\multicolumn{5}{|c|}{\textbf{MCMC-FEM (no neural network)}} & 34.25 &2.990881 &  0.281856 & 0.363245 \\
\hline
\end{tabular}
}
\end{table}
\begin{table}[!h]
\caption{Quantitative results of conductivity reconstructions using MCMC-Net with varying numbers of layers (NL) and neurons (NN), trained under different minibatch sizes (MB) and epochs (EP) for the QPAT inverse problem. Reported metrics include training time (Train. Time) and inversion time (Inv. Time), both in minutes, along with the $L^{\infty}$ -loss, mean absolute error (MAE), and mean squared error (MSE). All reconstructions are based on MCMC sampling with 100,000 samples. The bottom row shows results from a traditional FEM-based MCMC inversion. The reconstructions correspond to QPAT inverse problem with simulated data corrupted by 1\% relative noise.}
\label{tab:qpat_experiments_layers_neurons_mb_ep}
\centering
\renewcommand{\arraystretch}{1.2} 

\scalebox{0.94}{
\begin{tabular}{|c|c|c|c|c|c|c|c|c|}
\hline
NL & NN & MB & EP & Train. Time (m) & Inv. Time (m) & $L^{\infty}$ -loss & MAE & MSE \\
\hline
\textbf{4} & \textbf{16} & \textbf{8} & \textbf{100} & \textbf{20.62} & \textbf{5.58} & \textbf{0.200000} & \textbf{0.002066} & \textbf{0.000182} \\
4 & 32 & 8 & 100 & 20.85 & 5.56 & 0.200000  & 0.002452 & 0.000248 \\
4 & 64 & 8 & 100 & 24.21 & 5.89 & 0.240000  & 0.001887 & 0.000191 \\
\hline
6 & 16 & 8 & 100 & 23.94 & 6.47 & 0.240000  & 0.001949 & 0.000183 \\
6 & 32 & 8 & 100 & 23.74 & 6.04 & 0.200000  & 0.001915 & 0.000177 \\
6 & 64 & 8 & 100 & 26.60 & 5.92 & 0.200000  & 0.001811 & 0.000167 \\
\hline
8 & 16 & 8 & 100 & 27.26 & 6.30 & 0.200000  & 0.002225 & 0.000210 \\
8 & 32 & 8 & 100 & 26.58 & 5.75 & 0.360000  & 0.002824 & 0.000330 \\
8 & 64 & 8 & 100 & 29.33 & 6.23 & 0.280000  & 0.002273 & 0.000231 \\
\hline
\multicolumn{5}{|c|}{\textbf{MCMC-FEM (no neural network)}} & 11.15 &0.200000 &  0.001474 & 0.000119 \\
\hline
\end{tabular}
}
\end{table}


\section{Noise sensitivity analysis for EIT and DOT}
Using the network architectures specified in subsection \ref{sec:5.3}, we carry out numerical experiments to contrast the performance of the neural operator based method for Bayesian inversion of the EIT and DOT inverse problems, respectively. We report here the error incurred in reconstructions using both methods at varying noise levels. Quantitative results for the EIT and DOT inverse problems are provided the Tables~\ref{tab:noise_inversion_mae_mse_grouped} and~\ref{tab:dot_noise_inversion_mae_mse_grouped}, respectively.

\begin{table}[!h]
\caption{Inversion time and reconstruction errors for conductivity reconstruction in EIT with a single anomaly, evaluated at different noise levels.}
\label{tab:noise_inversion_mae_mse_grouped}
\centering
\renewcommand{\arraystretch}{1.2} 
\scalebox{0.95}{
\begin{tabular}{|c|ccc|ccc|}
\hline
\multirow{2}{*}{Noise Level} & \multicolumn{3}{c|}{MCMC-FEM} & \multicolumn{3}{c|}{MCMC-Net} \\
\cline{2-7}
 & Inv. Time (m) & MAE & MSE & Inv. Time (m) & MAE & MSE \\
\hline
2\%  & 125.76 & 0.083293 & 0.333083 & 4.17 & 0.086267 & 0.344979 \\ 
4\%  & 124.62 & 0.121954 & 0.487732 & 3.95 & 0.116006 & 0.463941 \\ 
6\%  & 124.93 & 0.157641 & 0.630484 & 4.10 & 0.157641 & 0.630484 \\ 
8\%  & 124.40 & 0.178459 & 0.713758 & 4.05 & 0.172511 & 0.689964 \\ 
10\% & 125.23 & 0.217121 & 0.868407 & 4.04 & 0.214147 & 0.856511 \\
\hline
\end{tabular}
}
\end{table}


\begin{table}[h!]
\caption{Inversion time and reconstruction errors for absorption reconstruction in DOT with a single anomaly, evaluated at different noise levels.}
\label{tab:dot_noise_inversion_mae_mse_grouped}
\centering

\renewcommand{\arraystretch}{1.2} 

\scalebox{0.85}{
\begin{tabular}{|c|cccc|cccc|}
\hline
\multirow{2}{*}{\shortstack{Noise \\ Level}} & \multicolumn{4}{c|}{MCMC-FEM} & \multicolumn{4}{c|}{MCMC-Net} \\
\cline{2-9}
 & Inv. Time (m) & $L^{\infty}$ -loss & MAE & MSE & Inv. Time (m) & $L^{\infty}$ -loss & MAE & MSE \\
\hline
2\%  & 34.33 & 3.040900 & 0.273252 & 0.354816 & 4.70 & 2.418734 & 0.257732 & 0.252306 \\ 
4\%  & 34.03 & 3.052751 & 0.279802 & 0.378908 & 4.81 & 2.454142 & 0.263362 & 0.262916 \\ 
6\%  & 34.39 & 3.088769 & 0.276369 & 0.380220 & 4.73 & 2.484857 & 0.270487 & 0.282817 \\ 
8\%  & 34.28 & 3.174917 & 0.274359 & 0.382229 & 3.57 & 2.576837 & 0.274480 & 0.294167 \\ 
10\% & 32.43 & 3.183243 & 0.276013 & 0.389562 & 3.51 & 2.666927 & 0.273841 & 0.302945 \\
\hline
\end{tabular}
}
\end{table}




\section{Reconstructions with more MCMC iterations}
\begin{figure}[!h]
    \centering 
    \begin{subfigure}[b]{0.45\textwidth}
        \includegraphics[height=0.2\textheight]{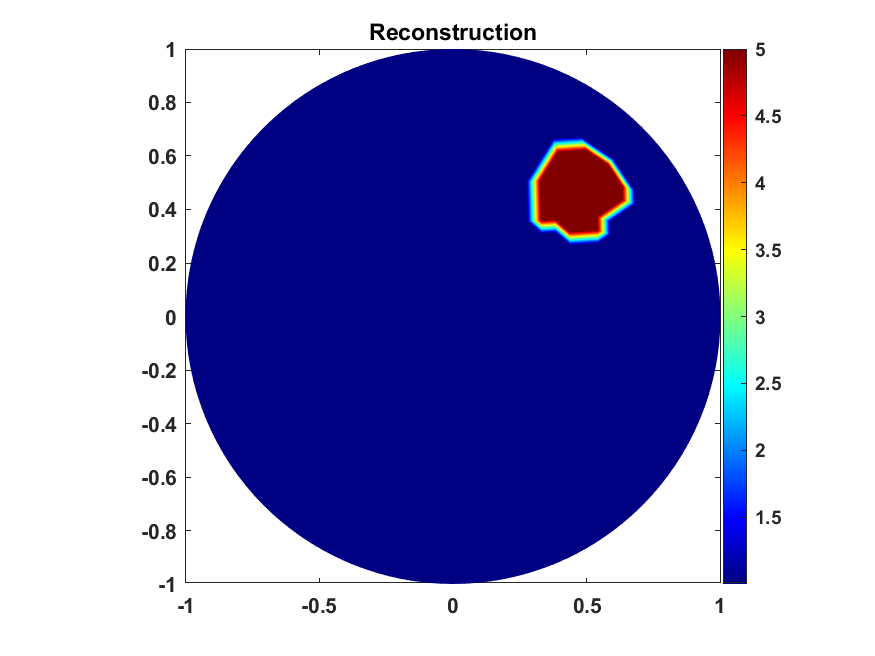}
        \caption{}
        \label{eit_SIGAN1_BRN_100K_TRS_100K}
    \end{subfigure}
    \hfill
    \begin{subfigure}[b]{0.45\textwidth}
        \includegraphics[height=0.2\textheight]{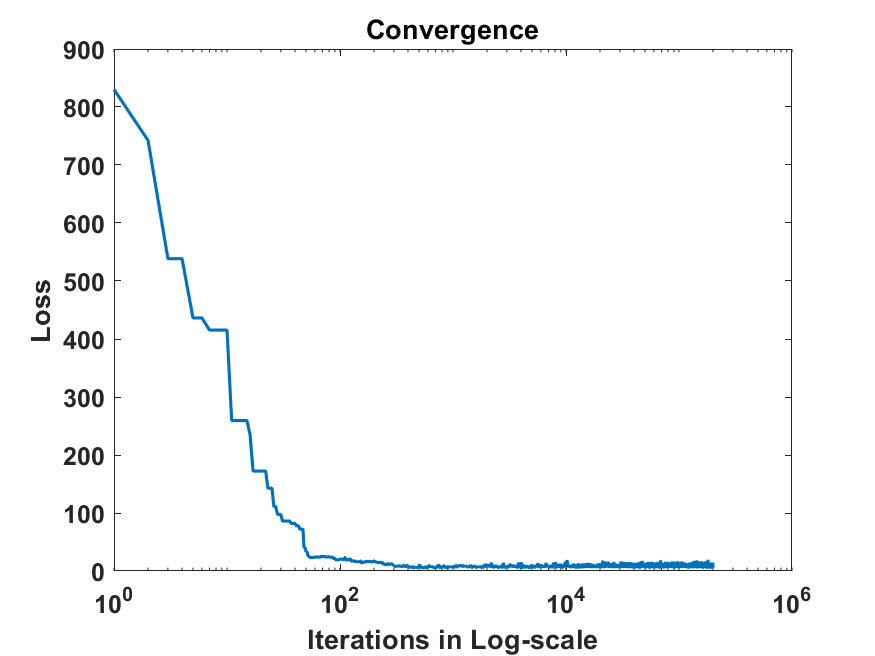}
        \caption{}
        \label{eit_SIGAN1_BRN_100K_TRS_100K_loss}
    \end{subfigure}
    \caption{Conductivity reconstruction using MCMC-Net in EIT after discarding 100,000 burn-in samples from a total of 200,000 samples. The reconstruction corresponds to a circular-shaped single anomaly with simulated data corrupted by 1\% relative noise. (a) Reconstructed conductivity. Inversion time: 8.08 minutes, MAE: 0.071386, MSE: 0.285500. (b) Convergence plot of the MCMC sampling.}
    \label{EITReconstructionAN1Loss}
\end{figure}

We mentioned in section \ref{NumericalExperiments} that to avoid the significant computational overhead in carrying out our extensive numerical experiments, we limited ourselves to 100,000 MCMC iterations for producing each reconstructed image. Here we produce an example reconstruction for the EIT inverse problem using MCMC-Net, where the number of MCMC iterations is 200,000. The resulting reconstruction and corresponding MCMC convergence behavior are shown in Figure~\ref{EITReconstructionAN1Loss}.

\section{Out-of-distribution tests} 
\begin{figure}[h!]
    \centering 

    \begin{subfigure}[b]{0.32\textwidth}
       \includegraphics[scale=0.35]{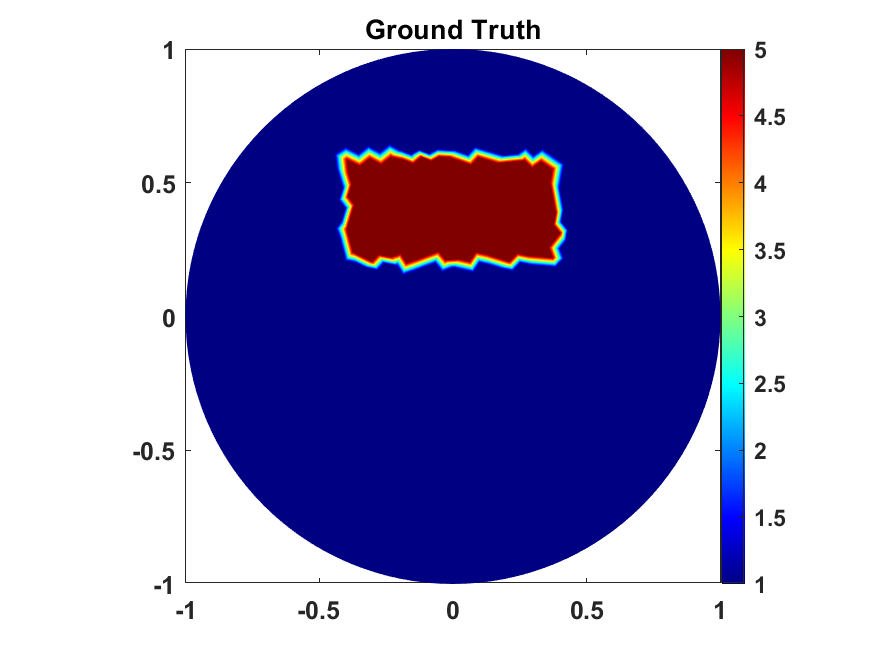}
       \caption{True}
        \label{dot_mutrue_sqbx_2129X1}
    \end{subfigure}
    \begin{subfigure}[b]{0.32\textwidth}
       \includegraphics[scale=0.35]{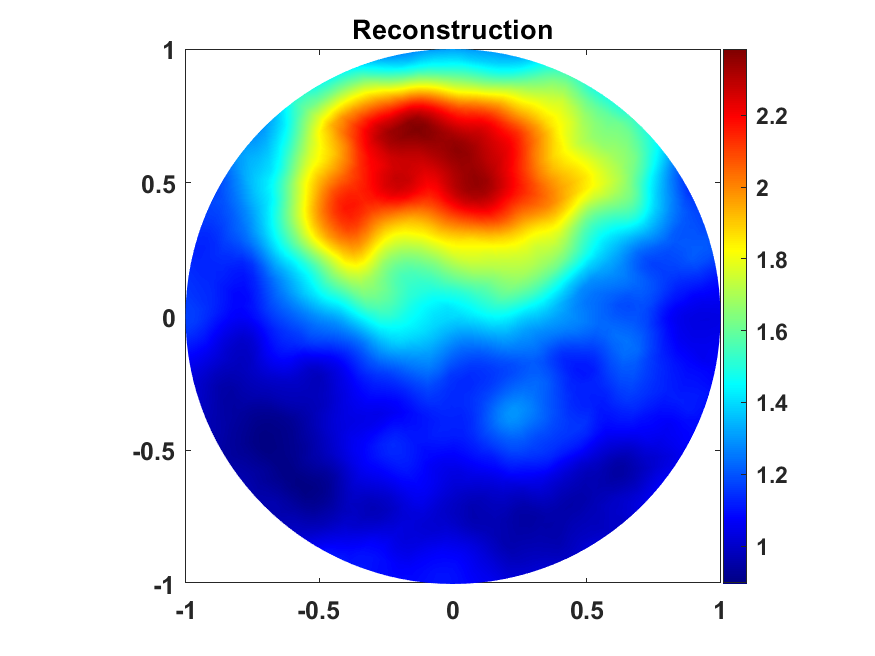}
        \caption{MCMC-Net}
        \label{dot_mu_CNN_sqbx_2129X1_nL1}
    \end{subfigure}
    \begin{subfigure}[b]{0.32\textwidth}
        \includegraphics[scale=0.35]{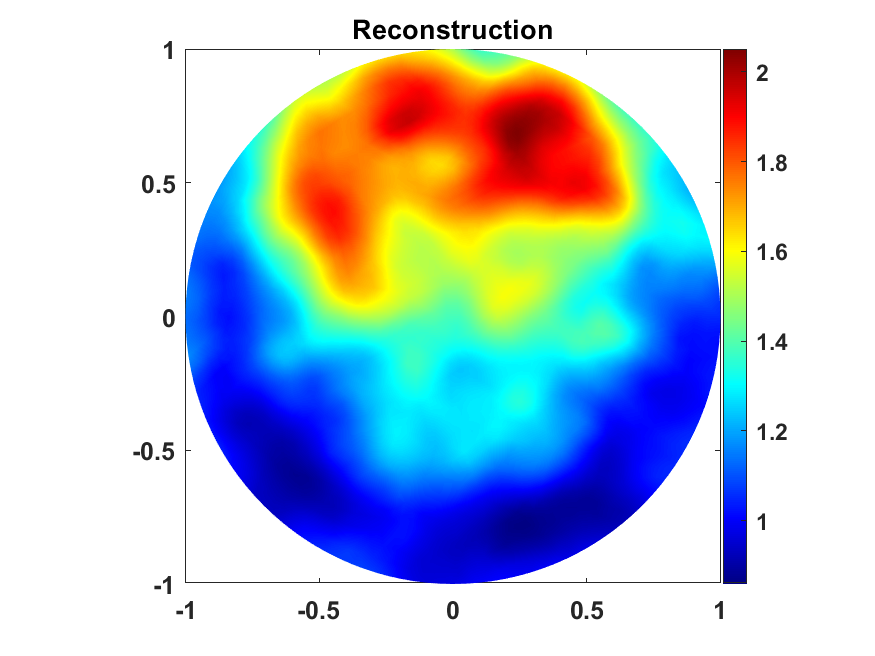}
       \caption{MCMC-FEM}
        \label{dot_mu_FEM_SQBX_2129X1_nL1}
    \end{subfigure}

    \begin{subfigure}[b]{0.32\textwidth}
        \includegraphics[scale=0.35]{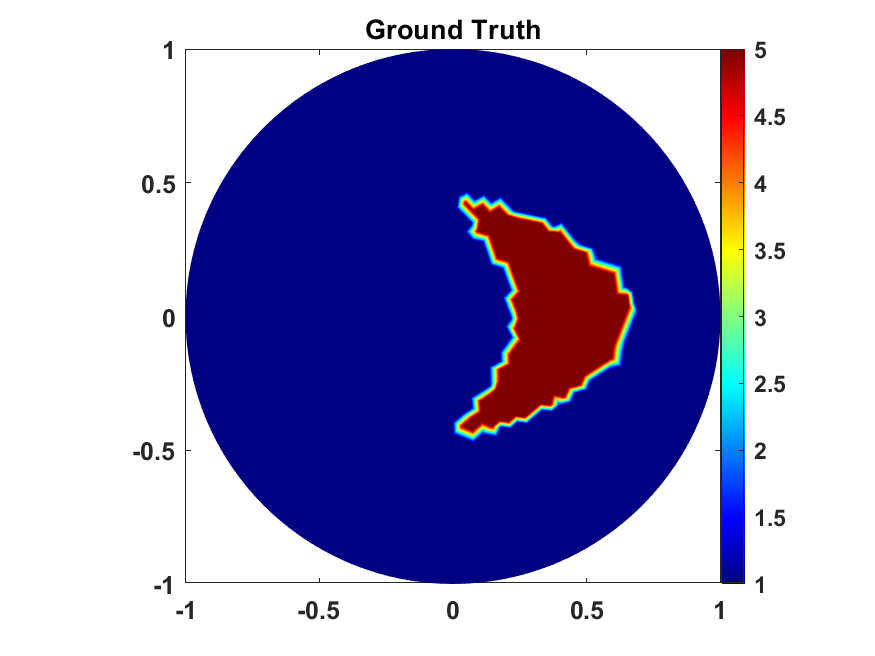}
        \caption{True}
        \label{dot_mutrue_moon_2129X1}
    \end{subfigure}
    \begin{subfigure}[b]{0.32\textwidth}
      \includegraphics[scale=0.35]{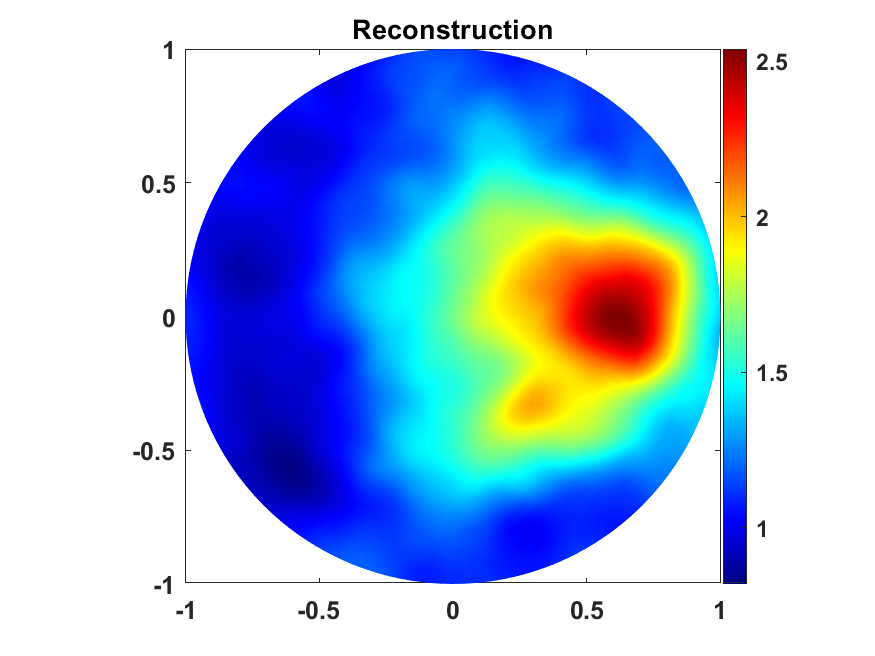}
        \caption{MCMC-Net}
        \label{dot_mu_CNN_moon_2129X1_nL1}
    \end{subfigure}
    \begin{subfigure}[b]{0.32\textwidth}
     \includegraphics[scale=0.35]{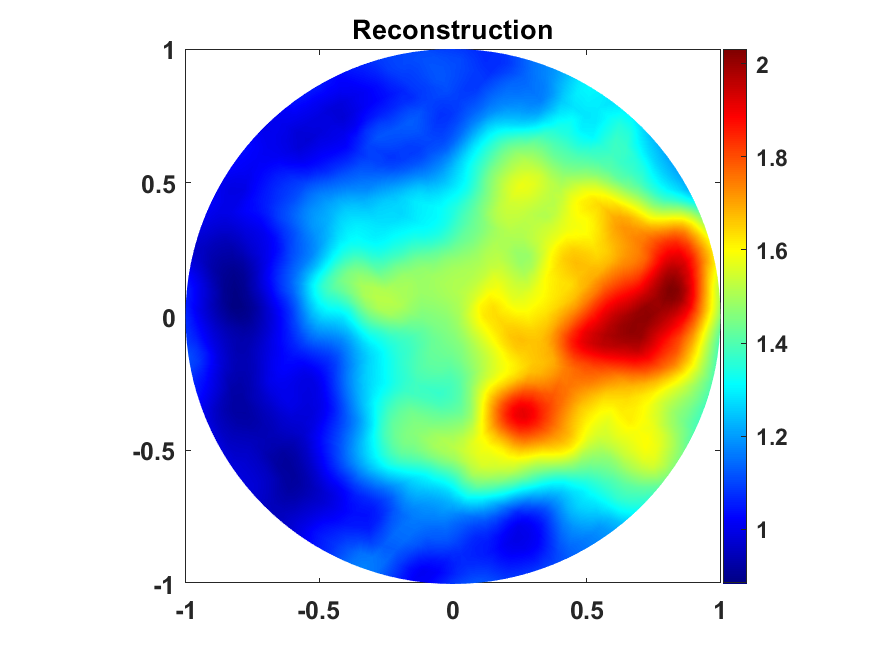}
       \caption{MCMC-FEM}
        \label{dot_mu_FEM_MOON_2129X1_nL1}
    \end{subfigure}

    \caption{Absorption reconstruction on unseen data for the DOT inverse problem. The CNN architecture is the same as used in Figure \ref{DOTan1Emu}. Data obtained at 1\% relative noise. 
    (a) Simulated absorption with a rectangle-shaped inclusion.
    (b) Reconstruction using MCMC-Net. Inversion time: 4.98 minutes, $L^{\infty}$ -loss: 3.399051, MAE: 0.595624, MSE: 1.118337.
    (c) Reconstruction using MCMC-FEM. Inversion time: 106.10 minutes, $L^{\infty}$ -loss: 3.511005, MAE: 0.641543, MSE: 1.323232.
    (d) Simulated absorption with a crescent-shaped inclusion.
    (e) Reconstruction using MCMC-Net. Inversion time: 4.95 minutes, $L^{\infty}$ -loss: 3.473424, MAE: 0.582329, MSE: 1.082405.
    (f) Reconstruction using MCMC-FEM. Inversion time: 106.34 minutes, $L^{\infty}$ -loss: 3.627973, MAE: 0.605065, MSE: 1.259771. Reconstructions are performed on a finite element mesh with 2129 grid points and 4128 triangles.}
    \label{DOTConductivityReconstructionsUnseenData}
\end{figure}

{{We also tested our networks for their generalization ability. Note that the forward operator networks were trained on phantoms with circular anomalous regions or were drawn from a Gaussian distribution with zero mean and a specified Matern kernel, see subsection \ref{sec:5.3} In Figure \ref{DOTConductivityReconstructionsUnseenData} we give examples of reconstructions for the DOT inverse problem for the case when the anomalous regions are non-circular and evidently not drawn from a Gaussian distribution. Note that MCMC-Net outperforms the traditional Bayesian inversion using the FEM solver.}}

\section{Reconstruction on a finer grid}
\begin{figure}[h!]
    \centering 
    \begin{subfigure}[b]{0.24\textwidth}
        \includegraphics[scale=0.31]{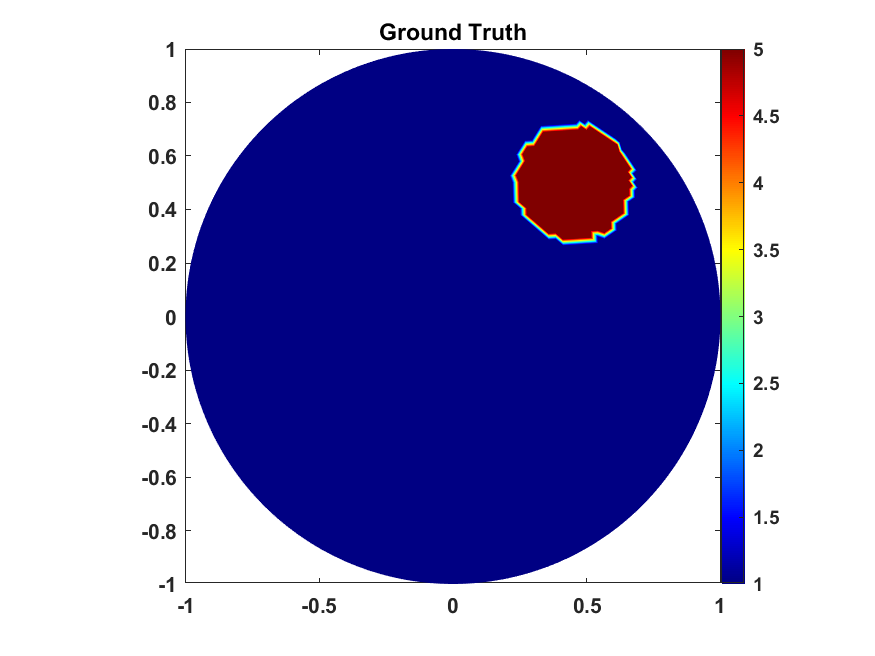}
        \caption{True}
        \label{EIT_sigmatrue_an1_5249X1}
    \end{subfigure}
    \begin{subfigure}[b]{0.24\textwidth}
        \includegraphics[scale=0.31]{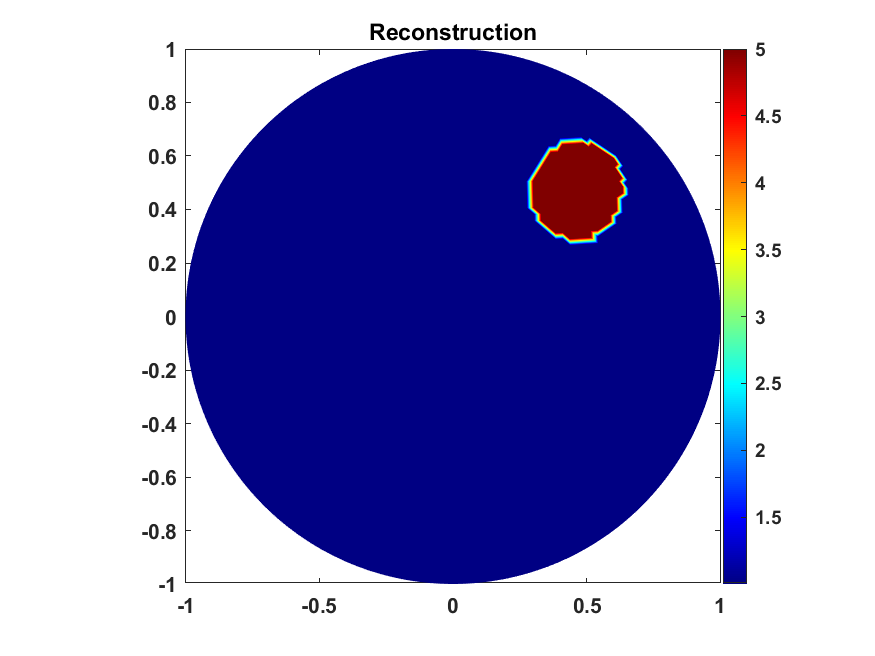}
        \caption{Reconstruction}
        \label{EIT_SIGAN1_5249X1_L4_N16_MB128_EP2K}
    \end{subfigure}
    \begin{subfigure}[b]{0.24\textwidth}
        \includegraphics[scale=0.31]{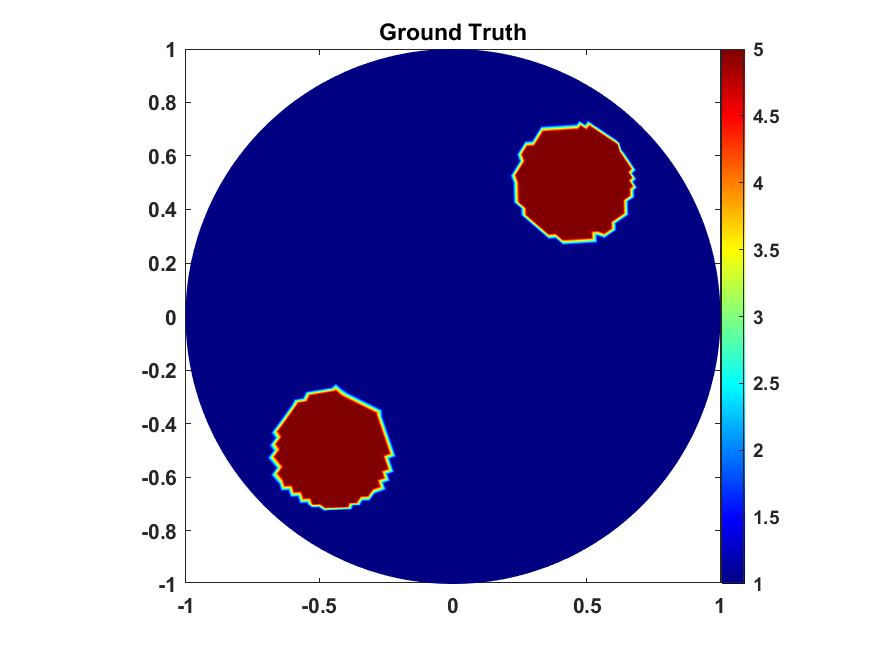}
        \caption{True}
        \label{EIT_sigmatrue_an2_5249X1}
    \end{subfigure}
    \begin{subfigure}[b]{0.24\textwidth}
        \includegraphics[scale=0.31]{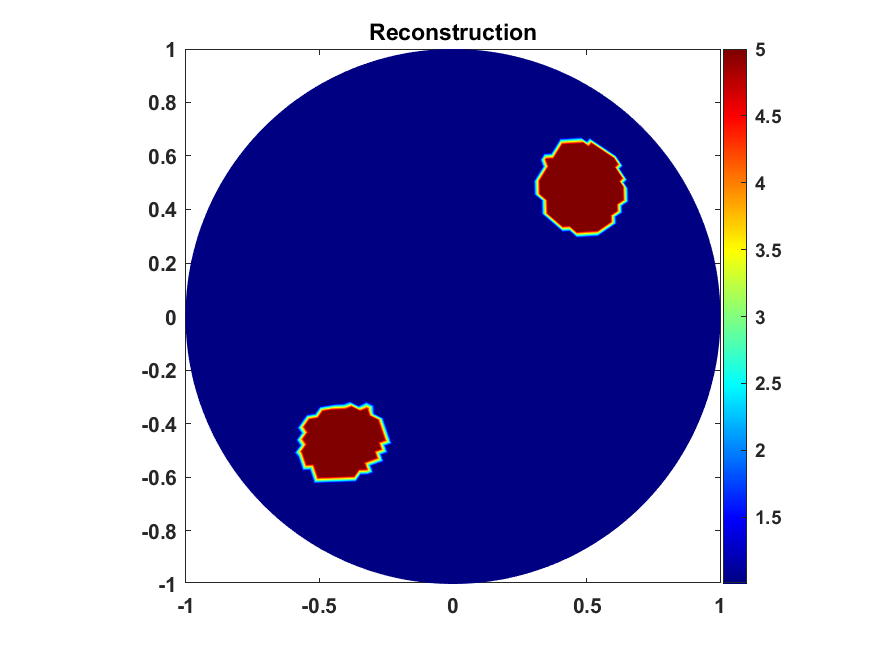}
        \caption{Reconstruction}
        \label{EIT_SIGAN2_5249X1_L4_N16_MB128_EP2K}
    \end{subfigure}

 \caption{Conductivity reconstructions on a finer grid. Reconstructions are performed on a finite element mesh with 5249 grid points and 10420 triangles. The CNN architecture consists of 4 layers with 16 neurons per layer, trained with a minibatch size of 128 and 2000 epochs. Data obtained at 1\% relative noise.
    (a) Simulated conductivity with a circular-shaped inclusion. 
    (b) Reconstruction corresponding to Figure~\ref{EIT_sigmatrue_an1_5249X1}. Inversion time: 11.52 minutes, MAE: 0.067082, MSE: 0.268239.
    (c) Simulated conductivity with two circular-shaped inclusions. 
    (d) Reconstruction corresponding to Figure~\ref{EIT_sigmatrue_an2_5249X1}. Inversion time: 11.91 minutes, MAE: 0.173008, MSE: 0.691934.
   }
    \label{EITConductivityReconstructionsSeenUnseenData}
\end{figure}

We also tried MCMC-Net on finer grids. The corresponding reconstructions are shown in Figure~\ref{EITConductivityReconstructionsSeenUnseenData}, which are obtained on a mesh that is four times finer than the one used in Figure~\ref{EITResultsDifferentHyperparameter}. We also report the inversion time using the MCMC-Net. 

\end{document}